\documentclass[11pt]{article}
\usepackage{amsmath}
\usepackage{amssymb}
\usepackage{theorem}
 \newtheorem{thm}{Theorem}[section]

 \theoremstyle{definition}
 
 \theoremstyle{remark}
 \newtheorem{rem}[thm]{Remark}
 \newtheorem{ex}[thm]{Example}
 \numberwithin{equation}{section}
\usepackage{mathtools}
\usepackage{mycommands}
\usepackage{xcolor}

\newcommand{\rank}{\mathrm{rank}}
\DeclareMathOperator*{\esssup}{ess\,sup}

\newcommand{\citechap}[1]{}
\newcommand{\refchapEIM}{Chapter 4}
\newcommand{\refchapRB}{Chapter 2}
\newcommand{\refchapRBtheory}{Chapter 3}
\newcommand{\refchapTENSORS}{Chapter 9}
\newcommand{\refchapPOD}{Chapter 1}

\begin{document}

%
%
%
%
%
%
%
%
%

\title
{\bfseries Low-rank methods for high-dimensional approximation and model order reduction}

\author{Anthony Nouy\footnotemark[2]}

\renewcommand{\thefootnote}{\fnsymbol{footnote}}
\footnotetext[2]{Ecole Centrale de Nantes, GeM, UMR CNRS 6183, France. e-mail: \texttt{anthony.nouy@ec-nantes.fr}
}

\renewcommand{\thefootnote}{\arabic{footnote}}

\date{}

\maketitle

%



\begin{abstract}
Tensor methods are among the most prominent tools for the numerical solution of high-dimensional problems where  functions of multiple variables have to be approximated. These methods exploit the tensor structure of function spaces and apply to many problems in computational science which are formulated in tensor spaces, such as problems arising in stochastic calculus, uncertainty quantification or parametric analyses.   
Here, we present complexity reduction methods  based on low-rank approximation methods. 
We analyze the problem of best approximation in subsets of low-rank tensors and discuss its connection with the problem of optimal model reduction in low-dimensional reduced spaces.  
We present different algorithms for computing approximations of a function in low-rank formats. In particular, we present constructive algorithms which are based either on a greedy construction of an approximation (with successive corrections in subsets of low-rank tensors) or on the greedy construction of tensor subspaces (for subspace-based low-rank formats). These algorithms can be applied for tensor compression, tensor completion or for the numerical solution of equations in low-rank tensor formats. 
A special emphasis is given to the solution of stochastic or parameter-dependent models. Different approaches are presented for the approximation of vector-valued or multivariate functions (identified with tensors), 
based on samples of the functions (black-box approaches) or on the models equations which are satisfied by the functions.
\bigskip

\end{abstract}

\maketitle

\section*{Introduction}

Low-rank approximation methods are among the most prominent complexity reduction methods for the solution of high-dimensional problems in computational science and engineering (see the surveys \cite{KOL09,CHI11,Khoromskij:2012fk,Grasedyck:2013} and monograph \cite{HAC12}). Typical problems include the solution of high-dimensional partial differential equations arising in physics or stochastic calculus, or the solution of parameter-dependent or stochastic equations using a functional approach, where functions of multiple (random) parameters have to be approximated. The construction of reduced order representations of the solution of complex parameter-dependent equations is of particular importance in parametric analyses (e.g. optimization, control, inverse problems) and uncertainty quantification (uncertainty propagation, sensitivity analyses, statistical inverse problems).

In practical applications, vector-valued or multivariate functions (as elements of tensor spaces) often present low-rank structures that can be efficiently exploited in order to reduce the complexity of their representation. 
In this chapter, we introduce the basic concepts on low-rank approximation, first for order-two tensors  and then for higher-order tensors. We present different methods for the approximation of a tensor, based either on a complete or incomplete information on the tensor, or on the knowledge of the equations satisfied by the tensor.
Particular emphasis is given to the solution of stochastic and parameter-dependent equations. 
\\

In Section \ref{sec:tensor_spaces}, we recall the definition and some useful properties of tensor Banach spaces.

In Section \ref{sec:low-rank-order-2}, we introduce the problem of the best rank-$r$ approximation of order-two tensors, with the Singular Value Decomposition (SVD) as a particular case (the case corresponding 
to tensor Hilbert spaces equipped with canonical norms). Emphasis is given to the case of Bochner spaces, which is of  particular interest for the analysis of parameter-dependent and stochastic equations.

In Section \ref{sec:Low-rank approximation of higher-order tensors}, we consider the case of higher-order tensors.  We first present different notions of rank and the associated low-rank approximation formats, with a special emphasis on subspace-based (or Tucker) formats. Then we discuss the problem of best approximation in subsets of low-rank tensors and its connection with the problem of finding optimal reduced spaces for the projection of a tensor (for subspace-based tensor formats). Then we present 
higher-order versions of the SVD which allow us to obtain quasi-best (and controlled) approximations in subspace-based tensor formats for the particular case of the approximation (compression) of a given tensor in a tensor Hilbert space equipped with a canonical norm. 

In Section \ref{sec:greedy}, we present constructive algorithms for the approximation in low-rank tensor formats.
These algorithms rely either on the greedy construction of the approximation, by defining successive corrections in a given subset of low-rank tensors (typically the set of rank-one tensors), or on the greedy construction of subspaces (for the approximation in subspace-based tensor formats). The latter approaches yield 
adaptive algorithms for projection-based model order reduction. For the case of parameter-dependent equations, these algorithms include the Empirical Interpolation Method (at the basis of Reduced Basis methods) and some variants of Proper Generalized Decomposition (PGD) methods.

In Section \ref{sec:low-rank-samples}, we present different approaches for the approximation of a function (identified with a tensor) in low-rank tensor formats, based on samples of the function. We present least-squares methods and interpolation methods, the latter ones being related to the problem of tensor completion.

In Section \ref{sec:parametric_problems}, we introduce a class of parameter-dependent (or stochastic) models and we show how these models can be formulated as tensor-structured equations, first by exploiting the order-two tensor structure of Bochner spaces, and then by exploiting higher-order tensor structures of Lebesgue spaces with product measures (e.g. induced by independent random parameters). 

Finally,  in Section \ref{sec:low-rank-equations}, we present low-rank methods for the solution 
of tensor-structured equations, relying either on the use of iterative solvers and standard low-rank compression methods, or on the minimization of a certain residual-based distance to the solution (using optimization algorithms in low-rank tensor manifolds or constructive algorithms). Particular emphasis is given to the case of parameter-dependent (or stochastic) equations. In this particular context, greedy algorithms provide adaptive 
methods for the construction of reduced-order models.

\section{Tensor spaces}\label{sec:tensor_spaces}

In this section, we introduce basic definitions on tensor Banach spaces and recall some useful properties.
For a detailed introduction to tensor analysis, we refer 
the reader  to the monographs 
 \cite{Defant:1993fk,light1985approximation,HAC12}. 

\subsection{Tensor Banach spaces}

Let us consider vector 
spaces $X_\nu$, $\nu \in \{1,\hdots,d\}$, equipped with norms ${\Vert\cdot\Vert_\nu}$.  For $(v^{(1)},\hdots,v^{(d)}) \in X_1\times \hdots\times X_d$, we denote by $\bigotimes_{\nu=1}^d v^{(\nu)}$ an elementary tensor. The 
\emph{algebraic tensor space} $X= \bigotimes_{\nu=1}^d X_\nu$ is defined as the linear span of elementary tensors:
$$
 {X}= \bigotimes_{\nu=1}^d X_\nu = \mathrm{span} \left\{ \bigotimes_{\nu=1}^d v^{(\nu)}  : v^{(\nu)} \in X_\nu, 1\le \nu\le d\right\},
$$
so that elements $v\in {X}$ can be written as finite linear combinations of elementary tensors, that means
\begin{align}
v = \sum_{i=1}^m v^{(1)}_i\otimes \hdots \otimes v^{(d)}_i\label{eq:canonical_decomp}
\end{align}
for some $m\in \mathbb{N}$ and some vectors $v_i^{(\nu)}\in X_\nu$, $1\le i\le m$, $1\le \nu\le d$.
A \emph{tensor Banach space} $X_{\Vert\cdot\Vert}$ equipped with a norm $\Vert\cdot\Vert$ 
is defined as the completion of an algebraic tensor space $X$
with respect to the norm $\Vert\cdot\Vert$, and we denote 
$
X_{\Vert\cdot\Vert} = \overline{X}^{\Vert\cdot\Vert} = {}_{\Vert\cdot\Vert} \bigotimes_{\nu=1}^d X_\nu .
$
If the norm $\Vert\cdot\Vert$ is associated with an inner product, the resulting space $X_{\Vert\cdot\Vert}$ is a \emph{tensor Hilbert space}. 
In the case of finite-dimensional spaces $X_\nu$, $X_{\Vert\cdot\Vert}$ does not depend on the choice of norm and it coincides with the normed 
algebraic tensor space $X$. 

\subsection{Tensor spaces of operators}
Let $X =  \bigotimes_{\nu=1}^d X_\nu$ and $Y =  \bigotimes_{\nu=1}^d Y_\nu$ be two normed algebraic tensor spaces. Let  $L(X_\nu, Y_\nu)$ (resp. $\Lc(X_\nu,Y_\nu)$) denote the set of linear operators (resp. continuous linear operators) from $X_\nu$ to $Y_\nu$. For $Y_\nu = \Rbb$, $L(X_\nu,\Rbb)=X_\nu^*$ is the algebraic dual space of $X_\nu$, while $\Lc(X_\nu,\Rbb)=X_\nu'$ is the continuous dual space of $X_\nu$.
For $A^{(\nu)} \in L(X_\nu, Y_\nu)$, $1\le \nu\le d $, the elementary tensor 
$A = \bigotimes_{\nu=1}^d A^{(\nu)} $
 is defined for elementary tensors $\bigotimes_{\nu=1}^d v^{(\nu)}\in X$ by 
$
A\left( \bigotimes_{\nu=1}^d v^{(\nu)}\right) = \bigotimes_{\nu=1}^d A^{\nu}(v^{(\nu)}),
$
and extended by linearity to the whole space $X$. The algebraic tensor space $ \bigotimes_{\nu=1}^d \Lc(X_\nu,Y_\nu)$ is defined in the same way. 
In the particular case where $Y=\Rbb$, with $Y_\nu=\Rbb$ for all $\nu$, an elementary tensor $  \bigotimes_{\nu=1}^d\varphi^{(\nu)} \in  \bigotimes_{\nu=1}^d X_\nu^*$ is such that for $v=\bigotimes_{\nu=1}^d v^{(\nu)} \in X$, 
$
\left(\bigotimes_{\nu=1}^d\varphi^{(\nu)} \right) (v) = \prod_{\nu=1}^d \varphi^{(\nu)}(v^{(\nu)}),
$
and we have $ \bigotimes_{\nu=1}^d X_\nu^* \subset X^*$.

\subsection{Minimal subspaces}
The \emph{minimal subspaces} of a tensor $v\in \bigotimes_{\nu=1}^d X_\nu$, denoted $U^{min}_\nu(v)$ for $1\le \nu \le d$, 
are defined by the property that 
$v\in \bigotimes_{\nu=1}^d U_\nu$ implies $U^{min}_\nu(v) \subset U_\nu$, while $v\in \bigotimes_{\nu=1}^d U^{min}_\nu(v)$. The minimal subspace
$U^{min}_\nu(v) \subset X_\nu $  can be equivalently characterized by 
$$
U^{min}_\nu(v) = \left\{ (id_{\nu} \otimes \varphi_{\nu^c}) (v) :  \varphi_{\nu^c} \in \bigotimes_{\beta\neq \nu} X^*_\beta \right\},
$$
where $id_\nu \in L( X_\nu , X_\nu)$ is the identity operator on $X_\nu$ and where we use the convention 
$id_{\nu} \otimes \left(\bigotimes_{\beta\neq \nu} \varphi_\beta\right) = \varphi_1\otimes \hdots \otimes \varphi_{\nu-1} \otimes 
id_\nu \otimes \varphi_{\nu+1} \otimes \hdots \otimes \varphi_d$. 
For $v$ having the representation \eqref{eq:canonical_decomp}, $U^{min}_\nu(v) \subset \mathrm{span} \{v_i^{(\nu)}\}_{i=1}^m,$ with an equality if the $m$ vectors $\{\otimes_{\beta\neq \nu} v_i^{(\beta)}\}_{i=1}^m$ are linearly independent. 
A minimal subspace can also be defined for any subset of dimensions $\alpha\subset  \{1,\hdots,d\}$ such that $1\le \#\alpha<d$. Letting $X_{\alpha} = \bigotimes_{\nu\in \alpha} X_\nu$, the minimal subspace $U_{\alpha}^{min}(v)\subset X_\alpha$  of $v$ is defined by
$$
U^{min}_\alpha(v) = \left\{ (id_{\alpha} \otimes \varphi_{\alpha^c}) (v) :  \varphi_{\alpha^c} \in \bigotimes_{\beta\notin \alpha} X^*_\beta \right\} .
$$
For $v$ having the representation \eqref{eq:canonical_decomp}, 
$U^{min}_\alpha(v) \subset \mathrm{span} \{ v_i^{(\alpha)}\}_{i=1}^m,$ with $v_i^{(\alpha)}= \bigotimes_{\nu\in \alpha} v_i^{(\nu)}$, and $U^{min}_\alpha(v) = \mathrm{span} \{ v_i^{(\alpha)}\}_{i=1}^m$ if the vectors $\{\otimes_{\nu\notin \alpha} v_i^{(\nu)}\}_{i=1}^m$ are linearly independent.
For $\alpha= \dot{\cup}_{k=1}^K \alpha_k $ being the disjoint union of non empty sets $\alpha_k\subset \{1,\hdots,d\}$, it holds 
$$
U^{min}_\alpha(v)  \subset \bigotimes_{k=1}^K U^{min}_{\alpha_k}(v).
$$
For a detailed introduction to minimal subspaces and their properties, see \cite{Falco:2012uq}.

\subsection{Tensor norms}
A norm $\Vert\cdot\Vert$ on $X$ is called a \emph{crossnorm} if 
$
\Vert \bigotimes_{\nu=1}^d v^{(\nu)}\Vert = \prod_{\nu=1}^d \Vert v^{(\nu)} \Vert_\nu
$
for all $(v^{(1)},\hdots,v^{(d)})\in X_1\times\hdots\times X_d$. For $\nu\in \{1,\hdots,d\}$, let $X_\nu' = \Lc(X_\nu,\Rbb)$ denote the continuous dual of 
$X_\nu$ equipped with the dual norm $\Vert \cdot\Vert'_\nu$ of $\Vert\cdot\Vert_\nu$. If $\Vert\cdot\Vert$ is a crossnorm and also the dual norm $\Vert\cdot\Vert'$ of $\Vert\cdot\Vert$ is a crossnorm on $\bigotimes_{\nu=1}^d X_\nu'$, that means 
$
\Vert \bigotimes_{\nu=1}^d \varphi^{(\nu)} \Vert' = \prod_{\nu=1}^d \Vert\varphi^{(\nu)} \Vert_\nu'
$
for all $\varphi^{(\nu)}\in X_\nu'$, then $\Vert\cdot\Vert$ is called a \emph{reasonable crossnorm}. For a reasonable crossnorm, the elementary tensor $\bigotimes_{\nu=1}^d \varphi^{(\nu)} $ is in  the space $X' = \Lc(X,\Rbb)$ equipped with the dual norm $\Vert\cdot\Vert'$, and it can be extended to an element in $(X_{\Vert\cdot\Vert})' = \Lc(X_{\Vert\cdot\Vert},\Rbb)$.
 A norm $\Vert\cdot\Vert$ on $X$ is said to be a  \emph{uniform crossnorm} if it is a reasonable crossnorm and 
if for any elementary operator $ \bigotimes_{\nu=1}^d A^{(\nu)}\in \bigotimes_{\nu=1}^d \Lc(X_\nu,X_\nu)$ and for any tensor  $v \in X$, it satisfies
$
\Vert \left(\bigotimes_{\nu=1}^d A^{(\nu)} \right) (v)\Vert \le \left( \prod_{\nu=1}^d \Vert A^{(\nu)} \Vert_{X_\nu\leftarrow X_\nu} \right) \Vert v\Vert,
$
where $\Vert A^{(\nu)} \Vert_{X_\nu\leftarrow X_\nu}$ denotes the operator norm of $A^{(\nu)}$. Therefore, when $X$ is equipped with a uniform crossnorm, 
$A = \bigotimes_{\nu=1}^d A^{(\nu)} $ belongs to the  space $\Lc(X,X)$ of continuous operators from $X$ to $X$, and the operator norm of $A$ is $\Vert A \Vert_{X\leftarrow X}= \prod_{\nu=1}^d \Vert A^{(\nu)}\Vert_{X_\nu\leftarrow X_\nu}$. The operator $A$ can then be uniquely extended to a continuous operator $\overline A \in \Lc(X_{\Vert\cdot\Vert},X_{\Vert\cdot\Vert})$.
\\
\par 
Some norms can be directly defined from the norms $\Vert\cdot\Vert_\nu$ on $X_\nu$, $1\le \nu\le d$.
The \emph{injective norm} $\Vert\cdot\Vert_{\vee}$ is a particular uniform crossnorm defined for an algebraic tensor $v$ as 
$$
\Vert v\Vert_{\vee} = \sup\{(\varphi^{(1)}\otimes \hdots \otimes \varphi^{(d)}) (v) : \varphi^{(\nu)}\in X_\nu' , \Vert \varphi^{(\nu)}\Vert_\nu'=1,1\le \nu \le d \}.
$$
The \emph{projective norm} $\Vert\cdot\Vert_\wedge$ is another particular uniform 
crossnorm defined for an algebraic tensor $v$ as 
$$
\Vert v\Vert_{\wedge} = \inf \left\{ \sum_{i=1}^m \prod_{\nu=1}^d \Vert v^{(\nu)}_i \Vert_\nu : v = \sum_{i=1}^m \bigotimes_{\nu=1}^d v_i^{(\nu)} \right\},
$$
where the infimum is taken over all possible representations of $v$. 
The injective and projective norms are respectively the weakest and strongest reasonable crossnorms in the sense that for any 
reasonable crossnorm $\Vert \cdot \Vert$, we have 
$ \Vert\cdot\Vert_\vee \lesssim \Vert\cdot\Vert \lesssim \Vert\cdot\Vert_\wedge$, therefore yielding the following inclusions between the corresponding tensor Banach spaces: $X_{\Vert \cdot\Vert_\wedge} \subset X_{\Vert \cdot\Vert} \subset X_{\Vert \cdot\Vert_\vee} $.

In the case where spaces $X_\nu$, $1\le \nu \le d$,   are Hilbert spaces 
associated with inner products $\langle \cdot,\cdot\rangle_\nu$, a natural inner product, called the \emph{induced} or \emph{canonical inner product},
 can be defined 
for elementary tensors as 
$$
\langle \bigotimes_{\nu=1}^d v^{(\nu)}, \bigotimes_{\nu=1}^d w^{(\nu)} \rangle = \prod_{\nu=1}^d \langle v^{(\nu)}, 
 w^{(\nu)} \rangle_\nu,
$$
and extended by linearity to the whole algebraic tensor space $X$. This yields the definition of a natural tensor Hilbert space $X_{\Vert\cdot\Vert}$. The associated norm, called the \emph{canonical norm}, is in fact the unique crossnorm associated with an inner product, and it is a uniform crossnorm. 
 
\subsection{Examples of tensor Banach spaces}

Here, we introduce examples of tensor Banach spaces that are of particular importance in parametric and stochastic analyses. 

\subsubsection{$L^p$ spaces with product measure}\label{sec:Lp_spaces}
Let $(\Xi,\Sigma,\mu)$ be a measure space with $\Xi\subset \Rbb^s$ and $\mu$ a finite measure supported on $\Xi$ (e.g. a probability measure). 
For $1\le p \le \infty$, the Lebesgue space $L^p_\mu(\Xi)$ is defined as the Banach space of (equivalence classes of) 
measurable functions $v:\Xi \rightarrow \Rbb$ with finite norm
\begin{align*}
&\Vert v \Vert_{p} = \left( \int_{\Xi} \vert v(y) \vert^p \mu(dy)\right)^{1/p} \quad \text{for }  1\le p <\infty, \quad \text{and}
\\&\Vert v \Vert_{\infty} = \esssup_{y\in \Xi} \vert v(y)\vert \quad \text{for } p=\infty.
\end{align*}
Now, let us assume that $(\Xi,\Sigma,\mu)$ is the product of measure spaces 
 $(\Xi_\nu,\Sigma_\nu,\mu_\nu)$ where $\Xi_\nu\subset \Rbb^{s_\nu}$   and $\mu_\nu$ is a finite measure, $1 \le \nu \le d$ ($s=\sum_{\nu=1}^d s_\nu$). That means 
$\Xi=\Xi_1\times \hdots \times \Xi_d$, $\Sigma = \Sigma_1\otimes \hdots \otimes \Sigma_d$ and $\mu=\mu_1\otimes\hdots \otimes \mu_d$. 
We can define the algebraic tensor space $X =  L^p_{\mu_1}(\Xi_1)\otimes \hdots \otimes L^p_{\mu_d}(\Xi_d)$. The  natural injection from $X$ to $L^p_\mu(\Xi)  $ is such that $(v^{(1)}\otimes \hdots \otimes v^{(d)})(y_1,\hdots,y_d) = v^{(1)}(y_1)\hdots v^{(d)}(y_d)$ for $(y_1,\hdots,y_d) \in \Xi$.  
$X$ is then the set of functions $v$ in $L^p_\mu(\Xi)$ that can be written 
$$
v(y_1,\hdots,y_d) = \sum_{i=1}^m v_i^{(1)}(y_1)\hdots v_i^{(d)}(y_d),
$$
for some $m\in \Nbb$ and some functions $ v_i^{(\nu)} \in L^p_{\mu_\nu} (\Xi_\nu)$.
 We have the property that the resulting tensor Banach  space $X_{\Vert\cdot\Vert_p} = {}_{\Vert\cdot\Vert_p}\bigotimes_{\nu=1}^d
 L^p_{\mu_\nu}(\Xi_\nu)$ is such that
\begin{align*}&X_{\Vert\cdot\Vert_p} = L^p_\mu(\Xi) \quad \text{for }  1\le p<\infty, \quad  \text{and }  \\
&X_{\Vert\cdot\Vert_\infty} \subset L^\infty_\mu(\Xi)  \quad \text{for  } p=\infty,\end{align*}
with equality $X_{\Vert\cdot\Vert_\infty} = L^\infty_\mu(\Xi)$ if $\Xi$ is a finite set (see \cite{Defant:1993fk}). For any $p$, the norm $\Vert \cdot\Vert_p$  is a reasonable crossnorm. 
In the case $p=2$,   $L^2_\mu(\Xi)$ is a Hilbert space which can be identified with the tensor Hilbert space $X_{\Vert\cdot\Vert_2} = {}_{\Vert\cdot\Vert_2}\bigotimes_{\nu=1}^d  L^2_{\mu_\nu}(\Xi_\nu)$. The norm $\Vert\cdot\Vert_2$ is the canonical inner product norm, which is a uniform crossnorm. For $1< p<\infty$, $X_{\Vert\cdot\Vert_p}$ is reflexive and separable. 

\subsubsection{Bochner spaces}\label{sec:Bochner}
Bochner spaces are of particular importance in the analysis 
 of  parameter-dependent and stochastic equations.
Let $V$ denote a Banach space equipped with a norm $\Vert \cdot\Vert_V$ and let $(\Xi,\Sigma,\mu)$ denote a measure space, where $\Xi\subset \Rbb^s$ and $\mu$ is a finite measure (e.g. a probability measure). For $1\le p\le \infty$, the Bochner space $ L^p_\mu(\Xi;V)$ is the 
Banach space of all (equivalence classes of) Bochner measurable functions\footnote{See e.g. \cite[Section 1.5]{roubivcek2013nonlinear} for the definition of Bochner measurability and Bochner integrability.} $v:\Xi\rightarrow V$ with bounded norm 
\begin{align*}
&\Vert v \Vert_{p} = \left(\int_\Xi \Vert v(y) \Vert_{V}^p \mu(dy) \right)^{1/p}
\quad \text{
for } 1\le p<\infty, \quad \text{and}
\\
 &  
\Vert v \Vert_{\infty} = \esssup_{y\in \Xi}  \Vert v(y)\Vert_{V}\quad \text{for }
p=\infty.\end{align*}
Let us note that $L^p_\mu(\Xi) = L^p_{\mu}(\Xi;\Rbb)$. We can define the algebraic tensor space $X =   L^p_\mu(\Xi)\otimes V$  and the natural injection from $X$ to $L^p_\mu(\Xi;V)$ by $ \lambda \otimes v \mapsto \lambda(\cdot)v$, such that $(\lambda\otimes v)(y) = \lambda(y)v$ 
for $y\in \Xi$. The space $X$ is composed by functions that can be written 
$$
v(y)= \sum_{i=1}^m s_i(y) v_i,
$$
for some $m\in \Nbb$ and some vectors $v_i\in V$ and functions $s_i\in L^p_\mu(\Xi)$, $1\le i\le m$.
We have the property that the resulting tensor Banach space $X_{\Vert\cdot\Vert_p} =  L^p_\mu(\Xi) \otimes_{\Vert\cdot\Vert_p} V$ is such that
\begin{align*}
&X_{\Vert\cdot\Vert_p} = L^p_\mu(\Xi;V) \quad \text{for } 
1\le p<\infty, \quad \text{and} 
\\
&X_{\Vert\cdot\Vert_{\infty}} \subset L^\infty_\mu(\Xi;V) \quad \text{for } p=\infty,
\end{align*}
with an equality $X_{\Vert\cdot\Vert_{\infty}}= L^\infty_\mu(\Xi;V)$ if $V$ is a finite-dimensional space or if $\Xi$ is a finite set\footnote{Note that if $\Xi$ is a finite set, then for any $1\le p,q \le \infty$, the norms $\Vert\cdot\Vert_p$ and $\Vert \cdot\Vert_q$ are equivalent and therefore,  the topological tensor spaces $X_{\Vert\cdot\Vert_{p}}$ and $X_{\Vert\cdot\Vert_{q}}$ coincide.}. For any $1\le p \le \infty$, the norm $\Vert \cdot\Vert_p$  is a reasonable crossnorm. For $p=2$ and if $V$ is a Hilbert space, then $\Vert \cdot\Vert_2$ is an inner product norm which makes $L^2_\mu(\Xi;V)$ a Hilbert space. 
Then $X_{\Vert\cdot\Vert_2} = L^2_\mu(\Xi)  \otimes_{\Vert\cdot\Vert_2} V$ is a tensor Hilbert space and  $\Vert \cdot\Vert_2$ is the canonical norm, which is a uniform crossnorm. For $1< p <\infty$, if  $V$ is reflexive and separable, then the Bochner tensor space $L^p_\mu(\Xi)\otimes_{\Vert\cdot\Vert_p} V$ is reflexive (see  \cite[Proposition 1.38]{roubivcek2013nonlinear}).
 
\subsection{Approximation in finite-dimensional tensor spaces}\label{sec:approx-fi-dim}

Let $X_{\Vert\cdot \Vert} = {}_{\Vert\cdot\Vert} \bigotimes_{\nu=1}^d X_\nu$ be a tensor Banach space. Approximations of elements 
of $X_{\Vert\cdot\Vert}$ are typically searched in finite-dimensional subspaces of $X$ that can be constructed as follows. 
Let $\{\psi_{k_\nu}^{(\nu)}\}_{k_\nu\in I_\nu}$ be a set of linearly independent elements in $X_\nu$, with $I_\nu$ such that $\#I_\nu = n_\nu$. Let $X_{\nu,I_\nu} = \mathrm{span}\{\psi_{k_\nu}^{(\nu)}\}_{k_\nu\in I_\nu} \subset X_\nu$. Let $I=I_1\times \hdots \times I_d$. Then $$X_I = X_{1,I_1} \otimes \hdots \otimes X_{d,I_d}$$ is a finite-dimensional subspace of $X$ with dimension  $\# I = \prod_{\nu=1}^d  n_\nu$ and with a basis $\{\psi_k\}_{k\in I}$ defined by 
$
\psi_k = \psi^{(1)}_{k_1} \otimes \hdots \otimes \psi_{k_d}^{(d)}, $ $k=(k_1,\hdots,k_d) \in  I.
$
An element $u \in X_I$ can be written 
\begin{align}
u = \sum_{k\in I} a_k \psi_k  = \sum_{k_1\in I_1} \hdots \sum_{  k_d \in I_d} a_{k_1,\hdots,k_d} 
 \psi_{k_1}^{(1)} \otimes \hdots \otimes \psi_{k_d}^{(d)},\label{eq:approxu}
\end{align}
where the set of coefficients $a =(a_k)_{k\in I} \in \mathbb{R}^{I}$ can be identified with a tensor $a\in 
 \mathbb{R}^{n_1} \otimes \hdots \otimes \Rbb^{n_d}$.   
If $X$ is a Hilbert space with inner product $\langle\cdot,\cdot\rangle$ and  if the basis $\{\psi_k\}_{k\in I}$ is orthonormal, then the coefficients in \eqref{eq:approxu} are given by
$a_{k_1,\hdots,k_d} = \langle \bigotimes_{\nu=1}^d \psi_{k_\nu}^{(\nu)} , u\rangle . $\\

\paragraph{Complexity reduction using sparse and low-rank tensor methods.}
The approximation space $X_I$ has a dimension which grows exponentially with the dimension $d$, which makes unpractical standard linear approximation methods in $X_I$ for a large $d$.  We can distinguish two main families of complexity reduction methods in tensor spaces: 
low-rank approximation methods and sparse approximation methods. Sparse approximation methods aim at defining 
suitable index sets $J \subset I$ with 
small cardinality for the approximation of a tensor in the corresponding low-dimensional space $X_J = \mathrm{span}\{\psi_k\}_{k\in J}$.
 The construction of index sets $J$ can be based on a priori analyses 
 \cite{Nobile:2008b,Beck:2014} or on adaptive algorithms \cite{Bieri:2010,SCH11,Cohen:2010,Cohen:2011,chkifa2015breaking,Chkifa:2013,Chkifa:2014fk}. Sparse and low-rank methods exploit different low-dimensional structures of tensors and these two complexity reduction methods can also be combined \cite{Bachmayr:2015fk,2013arXiv1305.0030C}. In this chapter, we only focus on low-rank approximation methods. 
Note that in practical applications, complexity reduction methods are most often used for the approximation of tensors in a fixed finite-dimensional space $X_I$, possibly adapted afterwards using a posteriori error estimates (see e.g. \cite{Bachmayr:2015fk}). Thus, low-rank and sparse tensor methods aim at finding a representation of the form \eqref{eq:approxu}
with a low-dimensional representation of the tensor of coefficients $a$.

\subsection{About best approximation problems}\label{sec:bestapproxproblem}
Here we recall definitions and classical  results on the problem of best approximation of an element $u \in X_{\Vert\cdot\Vert}$ from a subset $M$ in $X_{\Vert\cdot\Vert}$,
\begin{align}
\min_{v \in M} \Vert u - v \Vert.\label{bestapproxproblem}
\end{align}
A subset $M$ is \emph{proximinal} if for any $u$, there exists an element of best approximation in $M$. 
Any finite-dimensional linear subspace of $X_{\Vert\cdot\Vert}$ is proximinal. When  $X_{\Vert\cdot\Vert}$ is reflexive, a sufficient condition for a subset $M$ to be proximinal is that 
$M$ is weakly closed. 
In particular, any 
 closed convex set of a normed space is weakly closed. When 
$X_{\Vert\cdot\Vert}$ is finite-dimensional or when $M$ is a subset of a finite-dimensional subspace in $X_{\Vert\cdot\Vert}$, then a sufficient condition for $M$ to be proximinal is that $M$ is closed.
\par
A subset $M$ is a \emph{unicity set} if for any $u$, there exists at most one element of best approximation of $u$ in $M$. A subset $M$ is a \emph{Chebyshev set} if it is a proximinal unicity set, that means if for any $u$, there exists a unique element of best approximation of $u$ from $M$. Any convex subset $M$ of a strictly convex normed space is a unicity set.

\section{Low-rank approximation of order-two tensors}\label{sec:low-rank-order-2}

In this section, we consider the problem of the low-rank approximation of order-two tensors. 
We denote by $ S \otimes V$ an algebraic tensor space, where $S$ and $V$ are Banach spaces, 
and by $ S\otimes_{\Vert\cdot \Vert} V$ the corresponding tensor Banach space equipped with a norm $\Vert\cdot\Vert$. 
 
\subsection{Best rank-$r$ approximation} \label{sec:bestlowrank}
The \emph{rank} of $u \in S\otimes V$, denoted $\rank(u)$,  
 is  the minimal $r \in \Nbb$  such that
\begin{align}
u = \sum_{i=1}^r s_i \otimes v_i  , \label{eq:urankr}
\end{align}
for some vectors $\{v_i\}_{i=1}^r \in V^r$ and $\{s_i\}_{i=1}^r \in  S^r$.  
 We denote by $\Rc_r$  
  the set of tensors in $S\otimes V$ with a rank bounded by $r$, 
$$
\Rc_r = \left\{ \sum_{i=1}^r s_i \otimes v_i:   \{s_i\}_{i=1}^r \in  S^r,  \{v_i\}_{i=1}^r \in V^r\right\},
$$
or equivalently 
$$
\Rc_r = \left\{\sum_{i=1}^r \sum_{j=1}^r a_{ij} s_i \otimes v_j:   
 a=(a_{ij}) \in \Rbb^{r\times r}, \{s_i\}_{i=1}^r \in  S^r,\{v_i\}_{i=1}^r \in V^r\right\}.
$$
Let $u \in S\otimes_{\Vert \cdot\Vert} V$. An element $u_r$ of best approximation  of $u$ in $\Rc_r$ with respect to the norm $\Vert\cdot\Vert$ is such that 
\begin{align}
\Vert u-u_r \Vert = \min_{v\in \Rc_r} \Vert u -v \Vert.\label{eq:bestrankr}
\end{align}
If the norm $\Vert\cdot\Vert$ is not weaker than the injective norm, then $\Rc_r$ is weakly closed in $S\otimes_{\Vert \cdot\Vert} V$  (see Lemma 8.6 in \cite{HAC12}), and therefore proximinal if $S\otimes_{\Vert \cdot\Vert} V$ is reflexive. However, $\Rc_r$ is not a convex set and there is no guaranty of uniqueness of an element of best approximation. 
\begin{ex}\label{ex:bochnerbestlowrankapprox}
As an example, for $1< p <\infty$ and $V$ a reflexive and separable Banach space, the Bochner tensor space $L^p_\mu(\Xi)\otimes_{\Vert\cdot\Vert_p} V$ is reflexive and $\Vert\cdot \Vert_p$ is not weaker than the injective norm (see Section \ref{sec:Bochner}). Therefore, $\Rc_r$ is proximinal in $L^p_\mu(\Xi)\otimes_{\Vert\cdot\Vert_p} V$ if $1< p <\infty$ and $V$ is a reflexive and separable Banach space. 
\end{ex}

\subsection{Optimal subspaces}\label{sec:optimal_subspaces}
Now, we introduce equivalent reformulations of the best rank-$r$ approximation problem \eqref{eq:bestrankr} by 
using subspace-based parametrizations of $\Rc_r$.  We first note that $\Rc_r$ has a simple characterization using minimal subspaces. Indeed,
\begin{align*}
 \Rc_r &= \left\{ u \in S\otimes V :  \dim(U_1^{min}(u)) =  \dim(U_2^{min}(u)) \le r\right\},
  \end{align*}
 where the left and right minimal subspaces are respectively  
$U^{min}_{1}(u) = \left\{ (id_S \otimes  \varphi)  (u):  \varphi\in V^*\right\},$ 
 $U^{min}_{2}(u) = \left\{ (\psi \otimes id_V )  (u):  \psi\in S^*\right\}.$
Let $\Gbb_r(E)$ denote the \emph{Grassmann manifold} of $r$-dimensional subspaces in the vector space $E$.
First, we have 
 \begin{align}
 \Rc_r = &\left\{ u \in S_r \otimes V_r :  S_r \in \Gbb_r(S),V_r \in \Gbb_r(V)
 \right\},\label{eq:Rrsubspace1}
 \end{align}
 and the best rank-$r$ approximation problem \eqref{eq:bestrankr} can be  equivalently written 
  \begin{align}
  \min_{S_r\in \Gbb_r(S)} \min_{V_r\in \Gbb_r(V)}\min_{v\in S_r \otimes V_r}  \Vert u-v \Vert. \label{eq:bestrankrsubspace1}
 \end{align}
The solution of \eqref{eq:bestrankrsubspace1} yields optimal $r$-dimensional spaces $V_r$ and $S_r$ for the approximation of $u$ in the ``reduced'' tensor space $S_r\otimes V_r$.
Also, we have the following parametrization which only involves subspaces in $V$:
  \begin{align}
 \Rc_r =  &\left\{ u \in S \otimes V_r: V_r \in \Gbb_r(V) 
 \right\}, \label{eq:Rrsubspace2}
\end{align}
which yields the following reformulation of the best rank-$r$ approximation problem \eqref{eq:bestrankr}:
  \begin{align}
 \min_{V_r\in \Gbb_r(V)} \min_{v\in S \otimes V_r} \Vert u-v \Vert. \label{eq:bestrankrsubspace2}
 \end{align}
 The solution of \eqref{eq:bestrankrsubspace2} yields an optimal $r$-dimensional subspace $V_r$ for the approximation of $u$ in the ``reduced'' tensor space $S\otimes V_r$. 
 
\paragraph{Hilbert case.} Suppose that $S$ and $V$ are Hilbert spaces and that $S\otimes_{\Vert\cdot\Vert} V$ is a Hilbert space with a norm $\Vert\cdot\Vert$ associated with an inner product $\langle \cdot,\cdot\rangle$. For a finite-dimensional linear subspace $V_r\subset V$, let $P_{S\otimes V_r}$ denote the orthogonal projection from $S\otimes_{\Vert\cdot\Vert} V$ onto $S\otimes V_r$ such that 
\begin{align}\min_{v\in S\otimes V_r} \Vert u-v\Vert^2 = \Vert u - P_{S\otimes V_r} u\Vert^2 = \Vert u \Vert^2 - \Vert P_{S\otimes V_r} u \Vert^2.\label{eq:optimal_space_Vr_Hilbert}
\end{align}
The optimal subspace $V_r$ is  the solution of 
\begin{align}
\max_{V_r \in \Gbb_r(V)} \mathcal{R}_u(V_r) \quad \text{with} \quad \mathcal{R}_u(V_r)= \Vert P_{S\otimes V_r} u\Vert^2,\label{eq:maxrayleighVr}
\end{align}
which is an optimization problem on the Grassmann manifold $ \Gbb_r(V)$. 
The application 
$$\Vert \cdot\Vert_r : u\mapsto \Vert u \Vert_r = \max_{V_r\in \Gbb_r(V)} \Vert P_{S\otimes V_r} u\Vert= \max_{V_r\in \Gbb_r(V)} \max_{\substack{w\in S\otimes V_r \\ \Vert w\Vert=1}} \langle w, u \rangle$$ defines a norm on $S \otimes_{\Vert\cdot\Vert} V$ and the best rank-$r$ approximation $u_r$ satisfies 
$$
\Vert u-u_r\Vert^2 = \Vert u\Vert^2 - \Vert u_r\Vert^2 = \Vert u\Vert^2 - \Vert u\Vert_r^2.
$$
\begin{rem}
In the case where $\langle\cdot,\cdot\rangle$ is the canonical inner product,  $P_{S\otimes V_r} = id_S \otimes P_{V_r}$ where $P_{V_r}$ is the orthogonal projection from $V$ to $V_r$. Then, finding the optimal subspace $V_r$ is equivalent to finding the dominant eigenspace of an operator (see Section \ref{sec:svd}).
\end{rem}

\subsection{Tensors as operators}

The following results are taken from 
\cite[Section 4.2.13]{HAC12}. We restrict the presentation to the case where $V$ is a Hilbert space with inner product $\langle \cdot ,\cdot\rangle_V$.
An element 
$ u = \sum_{i=1}^r s_i\otimes v_i \in S\otimes V$ with rank $r$ can be identified with a rank-$r$ linear operator from 
$V$ to $S$ such that for $v\in V$, 
$$u(v) =\sum_{i=1}^r s_i \langle v_i,v \rangle_V .$$
Then, the  algebraic tensor space $S\otimes V$ can be identified 
with the set $\Fc(V,S)$ of finite rank operators from $V$ to $S$. The injective norm $\Vert\cdot\Vert_\vee$ coincides 
with the operator norm, so that the tensor Banach space
$S\otimes_{\Vert\cdot\Vert_\vee} V$ can be identified with 
the closure $\overline{\Fc(V,S)}$ of $\Fc(V,S)$ with respect to the operator norm, which coincides with the set of compact operators\footnote{$\overline{\Fc(V,S)}$ coincides with $\Kc(V,S)$ if the Banach space $V$ has the approximation property, which is the case for $V$ a Hilbert space.} 
$\Kc(V,S)$ from $V$ to $S$. Therefore, for any norm $\Vert \cdot\Vert$ not weaker than the injective norm, we have 
\begin{align}
 S\otimes_{\Vert\cdot\Vert} V
 \subset S\otimes_{\Vert\cdot\Vert_\vee} V = \Kc(V,S). \label{eq:tensor_space_compact_operators}
\end{align} 
Also, the tensor Banach space $S\otimes_{\Vert\cdot\Vert_\wedge} V$ equipped with the projective norm  $\Vert\cdot\Vert_\wedge$ can be 
identified with the space of nuclear operators 
$\Nc(V,S)$ from $V$ to $S$.  
Therefore, for any norm $\Vert \cdot\Vert$ not stronger than 
the projective norm $\Vert\cdot\Vert_\wedge$, we have 
\begin{align}
\Nc(V,S) = S\otimes_{\Vert\cdot\Vert_\wedge} V \subset S\otimes_{\Vert\cdot\Vert} V
 . \label{eq:tensor_space_nuclear_operators}
\end{align} 

\subsection{Singular value decomposition}\label{sec:svd}
In this section, we consider the case where 
$V$ and $S$ are Hilbert spaces. The spaces $V$ and $S$ are identified with their dual spaces $V'$ and $S'$ respectively. 
Let $\Vert\cdot\Vert$ denote the  canonical inner product norm. Let $n=\min\{\dim(V),\dim(S)\}$, with $n<\infty$ or $n=\infty$. 
Let $u$ in $ S\otimes_{\Vert\cdot\Vert_\vee} V = \Kc(V,S)$, the set of compact operators\footnote{Note that for $n<\infty$, $S\otimes_{\Vert\cdot\Vert_\vee} V = S\otimes V = \Fc(V,S) = \Kc(V,S).$}. Then, 
there exists a 
decreasing sequence of non-negative numbers $\sigma = \{\sigma_i\}_{i = 1}^n$ and two orthonormal systems 
$\{v_i\}_{i=1}^n \subset V$ and $\{s_i\}_{i=1}^n\subset S$ such that
\begin{align}
u = \sum_{i =1}^n \sigma_i s_i \otimes v_i, \label{eq:svd_compact}
\end{align}
where in the case $n=\infty$, the only accumulation point of the sequence $\sigma$ is zero and the series converges with respect to the injective norm $\Vert\cdot\Vert_\vee$ which coincides with the operator norm
 (see Theorem 4.114 in \cite{HAC12}). 
 The expression \eqref{eq:svd_compact} is  the so-called  \emph{Singular Value Decomposition} (SVD) of $u$, where $(s_i,v_i) \in S\times V$  is a couple 
of left and right singular vectors of $u$ 
associated with a singular value $\sigma_i$, verifying 
$$
u(v_i) = \sigma_i s_i \quad \text{and} \quad u^*(s_i) = \sigma_i v_i,
$$
where $u^* \in \Kc(S,V)$ is the adjoint operator of $u$ defined by $\langle s , u(v) \rangle_{S} = \langle u^*(s) , v\rangle_{V}$ for all $(v,s)\in V\times S$.
Let $u_r$ be the rank-$r$ truncated SVD defined by
$$
u_r = \sum_{i=1}^r \sigma_i s_i \otimes v_i.
$$
We have 
$$
\Vert u \Vert_\vee = \Vert \sigma\Vert_{\ell_\infty} = \sigma_1 \quad \text{and} \quad  \Vert u-u_r \Vert_\vee = \sigma_{r+1}.
$$
If we assume that $\sigma\in \ell_2$, then
$u\in S\otimes_{\Vert\cdot\Vert} V$ and 
$$
\Vert u\Vert=\Vert \sigma\Vert_{\ell_2}= \big(\sum_{i=1}^n \sigma_i^2\big)^{1/2},\quad \Vert u-u_r \Vert =  
\big(\sum_{i=r+1}^n \sigma_i^2\big)^{1/2} .
$$
The  canonical norm $\Vert \cdot\Vert$  coincides with the Hilbert-Schmidt norm of operators. 
We have the important property that 
$$
 \Vert u -u_r \Vert = \min_{w\in \Rc_r} \Vert u-w\Vert ,
$$
which means that an optimal rank-$r$ approximation of $u$ in the norm $\Vert\cdot\Vert$ can be obtained by retaining the first $r$ terms of the SVD. Moreover, 
\begin{align}
 \Vert u -u_r \Vert = \min_{w\in \Rc_1} \Vert u-u_{r-1} -w\Vert \label{eq:svdgreedy},
\end{align}
and
$$
 \Vert u -u_r \Vert^2 =  \Vert u -u_{r-1} \Vert^2 - \sigma_r^2 = \Vert u\Vert^2 - \sum_{i=1}^r \sigma_i^2 = \Vert u\Vert^2 - \Vert u_r\Vert^2.
$$
The $r$-dimensional  subspaces \begin{align*}&S_r = U^{min}_1(u_r) = \mathrm{span}\{s_i\}_{i=1}^r \in \Gbb_r(S) \text{ and } \\
&  
V_r = U^{min}_2(u_r) = \mathrm{span}\{v_i\}_{i=1}^r\in \Gbb_r(V)\end{align*} are 
respectively left and right dominant singular spaces of $u$. These  subspaces are solutions of 
problems
\eqref{eq:bestrankrsubspace1} and \eqref{eq:bestrankrsubspace2}, which means that they 
are optimal $r$-dimensional subspaces with respect to the canonical norm. Therefore, the SVD defines increasing sequences of 
optimal subspaces $\{V_r\}_{r\ge 1}$ and $\{S_r\}_{r\ge 1}$, such that 
$$
V_{r}\subset V_{r+1} \quad \text{and} \quad  S_{r} \subset S_{r+1}.
$$
Note that the optimal subspaces $V_r$ and $S_r$ are uniquely defined if $\sigma_{r}>\sigma_{r+1}$.
 Denoting $C_u: V\rightarrow V$ the compact operator defined by $C_u= u^*\circ u$, we have that $(v_i,\sigma_i^2) \in V\times \Rbb^+$ is an eigenpair of $C_u$, i.e.
$
C_u v_i = \sigma_i^2 v_i.
$
An optimal subspace $V_r$ is a dominant $r$-dimensional  eigenspace of $C_u.$ It is a solution 
of \eqref{eq:maxrayleighVr}. Here, 
the orthogonal projection
 $P_{S\otimes V_r}$ from $S\otimes_{\Vert\cdot\Vert} V$ to $S\otimes V_r$ is such that $P_{S\otimes V_r} = \overline{ id_S\otimes P_{V_r}}$, and
we have that 
$
\mathcal{R}_u(V_r) =  {R}_u(\mathbf{V}) = \mathrm{Trace}(\left\{ C_u \mathbf{V} , \mathbf{V} \right\}_V\left\{\mathbf{V} , \mathbf{V}\right\}_V^{-1}),
$
where $\mathbf{V} = \{v_i\}_{i=1}^r\in (V)^r$ is any basis of $V_r$, $C_u\mathbf{V}=\{C_uv_i\}_{i=1}^r$, and where 
$ \left\{ \{w_i\}_{i=1}^r , \{v_i\}_{i=1}^r \right\}_V = (\langle w_i , v_j \rangle_V )_{1\le i,j \le r} \in \Rbb^{r\times r}$. ${R}_u(\mathbf{V})$ is the Rayleigh quotient of $C_u$.

\subsection{Low-rank approximations in Bochner spaces}\label{sec:low-rank-bochner}

Here, we consider the particular case of low-rank approximations in 
Bochner spaces $L^p_\mu(\Xi ; V)$, $1\le p\le \infty$, where $\mu$ is a finite measure. This case is of particular interest for subspace-based model order reduction of  parameter-dependent (or stochastic) problems. 
Here we consider $V$ as a Hilbert space with norm $\Vert\cdot\Vert_V$. 
The considered algebraic tensor space is $L^p_\mu(\Xi)\otimes V$, and
the set $\Rc_r$ of   elements in $L^p_\mu(\Xi)\otimes V$ with rank at most $r$ is identified  with the set of functions $u_r:\Xi \rightarrow V$ of the form 
$$
u_r(y) = \sum_{i=1}^r s_i(y) v_i, \quad y\in \Xi.
$$
For a given $u\in L^p_\mu(\Xi ; V)$, let $\rho^{(p)}_r(u)$ denote the error of best rank-$r$ approximation in  $L^p_\mu(\Xi)\otimes V$, defined by
$$ \rho^{(p)}_r(u) = \inf_{w\in \Rc_r} \Vert u-w\Vert_p,$$
or equivalently by
$$
\rho^{(p)}_r(u) = \inf_{V_r\in \Gbb_r(V)} \inf_{w\in L^p_\mu(\Xi)\otimes V_r} \Vert u-w\Vert_p = \inf_{V_r \in \Gbb_r(V)} \Vert u-P_{V_r}u\Vert_p,$$
where $P_{V_r}$ is the orthogonal projection from $V$ to $V_r$ and $(P_{V_r}u)(y) = P_{V_r} u(y)$. 
For $1\le p<\infty$,  
$$
\rho^{(p)}_r(u) = \inf_{V_r \in \Gbb_r(V)} \left( \int_\Xi \Vert u(y) - P_{V_r} u(y) \Vert_V^p \mu(dy)\right)^{1/p},
$$
and for $p=\infty$,
$$
\rho^{(\infty)}_r(u) = \inf_{V_r \in \Gbb_r(V)}  \esssup_{y\in \Xi} \Vert u(y) - P_{V_r} u(y) \Vert_V.
$$
If we assume that $\mu$ is a probability measure, we have 
for all $1\le p\le q \le \infty$,
\begin{align*}
\rho^{(1)}_r(u) \le  \rho^{(p)}_r(u) \le \rho^{(q)}_r(u) \le \rho^{(\infty)}_r(u).
\end{align*}
There are two cases of practical importance. The first case is $p=2$, where $L^2_\mu(\Xi;V) = L^2_\mu(\Xi)\otimes_{\Vert \cdot\Vert_2} V$ is a Hilbert space and $\Vert \cdot\Vert_2$ is the canonical norm, so that we are in the situation where the best rank-$r$ approximation is the $r$-term truncated singular value decompositon of $u$ (see Section \ref{sec:svd}), called in this context Karhunen-Lo\`eve decomposition\footnote{Karhunen-Lo\`eve decomposition usually corresponds to the singular value decomposition of a centered second-order stochastic process $u$, that means of $u-\Ebb_\mu(u) = u-\int_\Xi u(y) \mu(dy)$.}. Then 
$
\rho^{(2)}_r(u)  = (\sum_{i\ge r+1}\sigma_i^2)^{1/2},$ 
where $\{\sigma_i\}_{i\ge 1}$ is the sequence of decreasing singular values of $u$.
The other important case is $p=\infty$. If we assume that $\Xi= \mathrm{support}(\mu)$ is compact and that $u$ is continuous from $\Xi$ to $V$, then the set of solutions $u(\Xi) = \{u(y) : y\in \Xi\}$ is a compact subset of $V$ and $\rho^{(\infty)}_r(u)$ coincides with the \emph{Kolmogorov $r$-width} $d_r(u(\Xi))_V$ of $u(\Xi)\subset V$,
\begin{align*}
\rho^{(\infty)}_r(u) &= \inf_{V_r \in \Gbb_r(V)}  \sup_{y\in \Xi} \Vert u(y) - P_{V_r} u(y) \Vert_V \\ &=
\inf_{V_r \in \Gbb_r(V)}  \sup_{v \in u(\Xi)} \Vert v - P_{V_r} v \Vert_V 
:= d_r(u(\Xi))_V.
\end{align*}
\begin{rem}
In the case $p=2$, there exists a sequence of nested optimal spaces $V_r$ associated with $\rho^{(2)}_r(u)$. In the case $p\neq 2$, up to the knowledge of the author, it remains an open question to prove whether or not there exists a sequence of nested optimal spaces.
\end{rem}

\section{Low-rank approximation of higher-order tensors}\label{sec:Low-rank approximation of higher-order tensors}

In this section, we consider the problem of the low-rank approximation of higher-order tensors and we will see how to extend the principles of Section \ref{sec:low-rank-order-2}.  Although several concepts apply to  general tensor Banach spaces (see \cite{Falco:2012fk,Falco:2012uq,falco2015geometric}), we restrict the presentation to the case of tensor Hilbert spaces. 

Let $X_\nu$, $\nu\in D:=\{1,\hdots,d\}$, denote Hilbert spaces equipped with norms $\Vert\cdot\Vert_\nu$ and associated inner products $\langle \cdot,\cdot \rangle_\nu$. 
We denote by $X= \bigotimes_{\nu \in D} X_\nu$ the algebraic tensor space, equipped with a norm $\Vert\cdot\Vert$ associated with an inner product $\langle\cdot,\cdot\rangle$, and by $X_{\Vert\cdot\Vert} $ the corresponding tensor Hilbert space.

\subsection{Low-rank tensor formats}

A subset $\Sc_r$ of low-rank tensors tensors in $X$ can be formally defined as a set $\Sc_r=\{v\in X : \rank(v) \le r\}$. There is no ambiguity in the case of order-two tensors, for which there is a unique notion of rank and $\Sc_r=\Rc_r$, with $r\in \Nbb$. However, there are several notions of rank for higher-order tensors, thus leading to different  subsets $\Sc_r$. For a detailed introduction to higher-order low-rank tensor formats, see \cite{KOL09,HAC12}\citechap{ and \refchapTENSORS{}}. Here, we briefly recall the main tensor formats, namely the canonical format and the subspace-based (or Tucker) formats. The approximation in the latter formats is closely related to subspace-based model order reduction.

\subsubsection{Canonical rank and canonical format}
The \emph{canonical rank} of a tensor $v \in X$ is the minimal integer $r\in\Nbb$ such that 
\begin{align}
v= \sum_{i=1}^r v_i^{(1)} \otimes \hdots \otimes v_i^{(d)}\label{eq:u_decomp_canonique}
\end{align}
for some vectors $v_i^{(\nu)}$, $1\le i\le r$, $1\le \nu \le d$. The set of  tensors with a canonical rank bounded by $r$ is denoted by $\Rc_r$. 
\begin{rem}\label{rem:param-Rr}
The elements of $ \Rc_r$ can be written $v= F_{\Rc_r}(\{v_{i}^{(\nu)}:1\le i \le r, 1\le \nu\le d \})$, where $F_{\Rc_r}$ is a multilinear map that parametrizes the subset $\Rc_r$ with 
$M(\Rc_r) =  r (\sum_{\nu=1}^d \dim(X_\nu))$ real parameters. We have $M(\Rc_r) \le d N r$, with $N=\max_\nu \dim(X_\nu)$.
\end{rem}

\subsubsection{$\alpha$-rank}
A natural notion of rank can be defined for a subset of dimensions, based on the notion of minimal subspaces. Let $\alpha \subset D$ be a subset of dimensions and $\alpha^c = D\setminus \alpha$, with $\alpha$ and $\alpha^c$ non empty. The \emph{$\alpha$-rank} of $v$, denoted $\rank_\alpha(v)$, is defined by
\begin{align}
\rank_\alpha(v) = \dim(U^{min}_\alpha(v)) .\label{eq:alpha-rank}
\end{align}
The $\alpha$-rank coincides with the classical notion of rank for order-two tensors.
A tensor $v\in X$  can be identified with a tensor $\Mc_\alpha(v) \in X_\alpha \otimes X_{\alpha^c}$, where $X_\alpha = \bigotimes_{\nu\in \alpha} X_\nu$ and $X_{\alpha^c} = \bigotimes_{\nu\in \alpha^c} X_\nu$, such that for $v$ of the form \eqref{eq:u_decomp_canonique},
$\Mc_\alpha(v) = \sum_{i=1}^r v_i^{(\alpha)} \otimes v_i^{(\alpha^c)}$, with $v_i^{(\alpha)} = \bigotimes_{\nu\in \alpha} v_i^{(\nu)}$ and 
$v_i^{(\alpha^{c})} = \bigotimes_{\nu\in \alpha^c} v_i^{(\nu)}$. $\Mc_\alpha : \bigotimes_{\nu \in D} X_\nu \rightarrow X_\alpha \otimes X_{\alpha^c}$ is a so-called ``\emph{matricisation}'' (or ``\emph{unfolding}'') operator.  The $\alpha$-rank of $v$ then coincides with the classical rank 
of the order-two tensor $\Mc_\alpha(v)$, i.e.
$
\rank_\alpha(v) = \rank(\Mc_\alpha(v)).
$
Subsets of low-rank tensors can now be defined by imposing the $\alpha$-rank for a collection of subsets $\alpha \in 2^D$. 
\begin{rem}\label{rem:matricisation_infinite}
The definition $\eqref{eq:alpha-rank}$  of the $\alpha$-rank also holds for elements  $v \in X_{\Vert\cdot\Vert}$. In this case, the interpretation as the rank of an order-two tensor requires the extension of the  matricisation operator to the topological tensor space $X_{\Vert\cdot\Vert}$.
\end{rem}

\subsubsection{Tucker rank and Tucker format}
The \emph{Tucker rank} (or \emph{multilinear rank}) of a tensor $v \in X$ is defined as the tuple $(\rank_\nu(v))_{\nu\in D}  \in \Nbb^d$. The set of tensors with a Tucker rank bounded by $r=(r_\nu)_{\nu\in D}$ is the set of \emph{Tucker tensors}
\begin{align*}
\Tc_r = \left\{ v\in X : \rank_{\nu}(v) = \dim(U_\nu^{min}(v)) \le r_\nu, \nu\in D \right\},
\end{align*} 
which can be equivalently characterized by 
\begin{align}
\Tc_r = \left\{ v\in U_1 \otimes \hdots \otimes U_d : U_\nu \in  \Gbb_{r_\nu}(X_\nu), \nu\in D\right\}. \label{eq:tucker_subspace}
\end{align} 
An element $v\in \Tc_r$ can be written 
$$ v = \sum_{i_1=1}^{r_1} \hdots \sum_{i_d=1}^{r_d} 
C_{i_1,\hdots,i_d} v_{i_1}^{(1)} \otimes \hdots \otimes v_{i_d}^{(d)} 
$$
for some $C\in \Rbb^{r_1\times \hdots\times r_d}$ (the \emph{core tensor}) and some $ v_{i_\nu}^{(\nu)} \in X_\nu$, $1\le i_\nu\le r_\nu$, $\nu\in D$. 
\begin{rem}\label{rem:param-Tr}
The elements of $ \Tc_r$ can be written $v= F_{\Tc_r}(C,\{v_{i_\nu}^{(\nu)}:1\le i_\nu \le r_\nu, 1\le \nu\le d \})$, where $F_{\Tc_r}$ is a multilinear map that parametrizes the subset $\Tc_r$ with 
$M(\Tc_r) = \prod_{\nu=1}^d r_\nu + \sum_{\nu=1}^d r_\nu \dim(X_\nu)$ real parameters. We have $M(\Tc_r) \le R^d + d N R$ with $R=\max_\nu r_\nu$ and $N=\max_\nu \dim(X_\nu)$.
\end{rem}

\subsubsection{Tree-based rank and tree-based Tucker format}
A more general notion of rank can be associated with a tree of dimensions.  
Let $T_D$ denote a \emph{dimension partition tree} of $D$, which is a subset of $2^D$ such that all vertices $\alpha \in T_D$ are non empty subsets of $D$, $D$ is the root of $T_D$, every vertex $\alpha\in T_D$ with $\#\alpha \ge 2$ has at least two sons, and the sons of a vertex $\alpha \in T_D$ form a partition of $\alpha$. The set of sons of $\alpha \in T_D$ is denoted $S(\alpha)$. A vertex $\alpha$ with $\#\alpha=1$ is called a leaf of the tree and is such that $S(\alpha)=\emptyset$. The set of leaves of $T_D$ is denoted $\Lc(T_D)$.   The \emph{tree-based  Tucker rank} of a tensor $u$ associated with a dimension tree $T_D$, denoted $T_D\text{-}\rank(u)$, is a tuple $(\rank_\alpha(u))_{\alpha\in T_D} \in \Nbb^{\#T_D}$. 
Letting $r = (r_\alpha)_{\alpha\in T_D} \in \Nbb^{\#T_D}$ be a tuple of integers,  the subset of \emph{tree-based Tucker tensors} with tree-based Tucker rank bounded by $r$ is defined by
\begin{align}
\Bc\Tc_r &= \left\{v\in X : \rank_\alpha(v) =  \dim(U^{min}_\alpha(v)) \le r_\alpha, \alpha\in T_D\right\}.\label{eq:treebasedsubset}
\end{align}
A tuple $r=(r_\alpha)_{\alpha\in T_D}$ is said admissible for $T_D$ if there exists an element $v \in X \setminus \{0\}$ such that $\dim(U^{min}_\alpha(v))=r_\alpha$ for all $\alpha \in T_D$. Here we use the convention $U^{min}_D(v) = \mathrm{span}\{v\}$, so that $r_D=1$ for $r$ admissible.  The set $\Bc\Tc_r$ can be  equivalently defined 
 by
\begin{align}
\Bc\Tc_r &= \left\{v\in \bigotimes_{\alpha \in S(D)} U_\alpha: \begin{array}{l} \displaystyle
U_\alpha \subset \bigotimes_{\beta \in S(\alpha)} U_\beta \text{ for all } \alpha\in T_D\setminus \{\Lc(T_D)\cup D\} \\
\text{and } \dim(U_\alpha)=r_\alpha \text{ for all } \alpha\in T_D\setminus D \end{array} \right\}.\label{eq:TBr_subspace}
\end{align}
For an element $v\in \Bc\Tc_r$ with an admissible tuple $r$, if $\{v^{(\alpha)}_{i_\alpha}\}_{i_\alpha =1}^{r_\alpha}$ denotes a basis of 
$U^{min}_\alpha(v)$ for $\alpha \in T_D$, with $v^{(D)}_1=v$, then for all $\alpha \in T_D\setminus \Lc(T_D)$,  
$$
v^{(\alpha)}_{i_\alpha} = \sum_{\substack{1\le i_\beta \le r_\beta \\\beta\in S(\alpha)}} C^{(\alpha)}_{i_\alpha,(i_\beta)_{\beta\in S(\alpha)}} \bigotimes_{\beta\in S(\alpha)} v^{(\beta)}_{i_\beta},
$$
for $1\le i_\alpha\le r_\alpha$, where the $C^{(\alpha)} \in \Rbb^{r_\alpha \times (\times_{\beta\in S(\alpha)} r_\beta )}$ are the so-called \emph{transfer tensors}. Then, proceeding recursively, we obtain the following  representation of $v$:
$$
v = \sum_{\substack{1\le i_\nu\le r_\nu \\ \nu \in D}} \left(\sum_{\substack{1\le i_\alpha \le r_\alpha \\ \alpha \in T_D \setminus \Lc(T_D) 
}}   \prod_{\substack{\mu \in T_D \setminus \Lc(T_D)}} C^{(\mu)}_{i_\mu,(i_\beta)_{\beta\in S(\mu)}}
\right)  \bigotimes_{\nu\in D} v^{(\nu)}_{i_\nu}.
$$
\begin{rem}\label{rem:param-BTr}
The elements of $ \Bc\Tc_r$ can be written $v= F_{\Bc\Tc_r}(\{v_{i_\nu}^{(\nu)}:1\le i_\nu \le r_\nu, 1\le \nu\le d \}, \{C^{(\alpha)}:\alpha\in T_D\setminus \Lc(T_D)\})$, where $F_{\Bc\Tc_r}$ is a multilinear map that parametrizes the subset $\Bc\Tc_r$ with $M(\Bc\Tc_r) = \sum_{\nu=1}^d r_\nu \dim(X_\nu) +  \sum_{\alpha \in T_D \backslash \Lc(T_D)} r_\alpha \prod_{\beta\in S(\alpha)} r_\beta  $ real parameters. We have $M(\Bc\Tc_r) \le d N R + R^{\#S(D)} + \sum_{T_D \backslash \{\Lc(T_D)\cup D\}} R^{\#S(\alpha)+1}  \le dNR +  R^S + (d-2) R^{S+1} $,  with $R=\max_\alpha r_\alpha$, $S=\max_{\alpha\notin \Lc(T_D)} \#S(\alpha)$, and  
$N=\max_\nu \dim(X_\nu)$.
\end{rem}
\begin{rem}
For a tree $T_D$ such that $S(D) = \Lc(T_D) = \{\{1\},\hdots,\{d\}\}$, 
the set $\Bc\Tc_{(1,r_1,\hdots,r_d)}$ coincides with the set of Tucker tensors $\Tc_{(r_1,\hdots,r_d)}$. 
For a binary tree $T_D$, i.e. such that $\# S(\alpha) = 2$ for all $\alpha \notin \Lc(T_D)$, the set $\Bc\Tc_{r}$ coincides with the set of Hierarchical Tucker (HT) tensors introduced in \cite{HAC09}.
\end{rem}

The reader is referred to \cite{HAC09,HAC12,falco2015geometric} for a detailed presentation of tree-based Tucker formats and their properties.

\subsubsection{Tensor-Train rank and Tensor-Train format}\label{TT-format}
The Tensor-Train (TT) format (see \cite{OSE11}) is a particular (degenerate) case of tree-based Tucker format which is associated with a particular binary dimension tree  $$T_D = \{\{k\}:1\le k \le d\} \cup \{\{k,\hdots,d\}: 1\le k \le d-1\}$$ such that 
$S(\{k,\hdots,d\}) = \{\{k\},\{k+1,\hdots,d\}\}$ for $1\le k \le d-1$.
The TT-rank of a tensor $u$, denoted $\rank_{TT}(u)$, is the tuple 
$(\rank_{\{k+1,\hdots,d\}}(u))_{k=1}^{d-1} $. For a tuple $r=(r_1,\hdots,r_d)\in \Nbb^{d-1}$, the set of tensors with TT-rank bounded by $r$ is defined by
\begin{align}
\Tc\Tc_r &= \left\{v\in X : \rank_{\{k+1,\hdots,d\}}(v) \le r_k
\right\},
\end{align}
which corresponds to the definition of a subset of tree-based Tucker tensors with inactive constraints on the ranks  $\rank_{\{k\}}(v)$ for $2\le k\le d-1$.
\begin{rem}
More precisely, $\Tc\Tc_r$ coincides with the subset $\Bc\Tc_{m}$ of tree-based Tucker tensors with a tree-based Tucker rank bounded by $m = (m_\alpha)_{\alpha\in T_D}$ if  $m$ is such that $m_{\{k+1,\hdots,d\}}=r_k$ for $1\le k \le d-1$ and $m_{\{k\}} \ge r_{k}r_{k+1}$ for $2\le k \le d-1$, the latter conditions implying that the constraints $\rank_{\{k\}}(v)\le m_{\{k\}}$ are  inactive for $2\le k \le d-1$.
\end{rem}
An element $v\in \Tc\Tc_r$ admits the following representation
$$
v = \sum_{i_1=1}^{r_1}\sum_{i_2=1}^{r_2}\hdots \sum_{i_{d-1}}^{r_{d-1}} v_{1,i_1}^{(1)} \otimes v_{i_1,i_2}^{(2)} \hdots \otimes v_{i_{d-1},1}^{(d)},
$$
where $v_{i_{\nu-1},i_{\nu}}^{(\nu)} \in X_\nu$.
\begin{rem}\label{rem:param-TTr}
The elements of $ \Tc\Tc_r$ can be written $v= F_{\Tc\Tc_r}(\{v^{(\nu)}:1\le \nu\le d \})$, with $v^{(\nu)}\in (X_\nu)^{r_{\nu-1}\times r_\nu}$ (using the convention $r_0=r_d=1$),  where $F_{\Tc\Tc_r}$ is a multilinear map that parametrizes the subset $\Tc\Tc_r$ with $M(\Tc\Tc_r) = \sum_{\nu=1}^d  r_{\nu-1}r_\nu \dim(X_\nu) $ real parameters. We have $M(\Tc\Tc_r) \le d N R^2  $,  with $R=\max_k r_k$ and 
$N=\max_\nu \dim(X_\nu)$.
\end{rem}

 \subsection{Best approximations in subspace-based low-rank tensor formats}\label{sec:best-approx-higher-order}
 
 \subsubsection{Tucker format}\label{sec:best-approx-Tucker}
Let us first consider the best approximation problem in Tucker format.
A best approximation of $u\in X_{\Vert\cdot\Vert}$ in the subset of Tucker tensors $\Tc_r$ with a rank bounded by $r=(r_1,\hdots,r_d)$ is defined by
\begin{align}
\Vert u-u_r\Vert = \min_{v\in \Tc_r} \Vert u-v\Vert. \label{eq:bestapproxTr}
\end{align} 
Based on the definition \eqref{eq:tucker_subspace} of $\Tc_r$, problem \eqref{eq:bestapproxTr}
can be equivalently written 
\begin{align}
\Vert u-u_r \Vert = \min_{U_1\in \Gbb_{r_1}(X_1)} \hdots \min_{U_d \in \Gbb_{r_d}(X_d)}
 \min_{v\in U_1\otimes \hdots \otimes U_d} \Vert u-v\Vert. \label{eq:bestapproxTr_subspace}
\end{align} 
A solution $u_r$ to problem \eqref{eq:bestapproxTr_subspace}
 yields optimal subspaces $U_\nu =  U^{min}_\nu(u_r)$ with dimension less than $r_\nu$, for $1\le \nu \le d$.  

\par 
Different conditions ensure that the set $\Tc_r$ is proximinal, that means that there exists a solution to the best approximation problem \eqref{eq:bestapproxTr} for any $u$ (see Section \ref{sec:bestapproxproblem}). 
If the norm $\Vert \cdot \Vert$ is not weaker than the injective norm, then $\Tc_r$ is weakly closed (see \cite{Falco:2012uq}), and therefore proximinal if $X_{\Vert\cdot\Vert}$ is reflexive (e.g. for  $X=\bigotimes_{\nu\in D} L^p_{\mu_\nu}(\Xi_\nu)$  for any $1<p<\infty$, see Section \ref{sec:Lp_spaces}). In particular, if $X$ is finite-dimensional, $\Tc_r$ is closed and therefore proximinal.
 
\subsubsection{Tree-based Tucker format}\label{sec:best-approx-treebased}
Let us now consider the best approximation problem in the more general tree-based Tucker format. 
The best approximation of $u\in X_{\Vert\cdot\Vert}$ in the subset of tree-based Tucker tensors $\Bc\Tc_r$ with $T_D$-rank bounded by $r=(r_\alpha)_{\alpha\in T_D}$ is defined by
\begin{align}
\Vert u-u_r\Vert = \min_{v\in \Bc\Tc_r} \Vert u-v\Vert. \label{eq:bestapproxBTr}
\end{align} 
Based on the definition \eqref{eq:TBr_subspace} of $\Bc\Tc_r$, Problem \eqref{eq:bestapproxBTr} can be equivalently written
\begin{align}
\Vert u-u_r\Vert = \min_{(U_\alpha)_{\alpha\in T_D \setminus D} \in \Gc_r(T_D)} \min_{v \in \bigotimes_{\alpha \in S(D)} U_\alpha}\Vert u-v\Vert, \label{eq:bestapproxBTr_subspace}
\end{align} 
where $\Gc_r({T_D})$ is a set of subspaces defined by
\begin{align*}
\Gc_r({T_D})= 
\Big\{ (U_{\alpha})_{\alpha \in T_D\setminus D} :  U_\alpha \in \Gbb_{r_\alpha}(X_\alpha) \text{ for all $\alpha \in T_D\setminus D$}, &\\
\text{ and } U_\alpha \subset \bigotimes_{\beta \in S(\alpha)}  U_\beta \text{ for all } \alpha \in T_D\setminus \{D \cup \Lc(T_D)\}    &\Big\}.
\end{align*}
Therefore, a best approximation $u_r \in \Bc\Tc_r$ yields a collection of optimal subspaces $U_\alpha$ with dimension  $r_\alpha$, $\alpha \in T_D\setminus D$, with a hierarchical structure. 

The proof of the existence of a best approximation in $\Bc\Tc_r$ requires some technical conditions involving norms defined for all the vertices of the tree (see \cite{falco2015geometric}). In particular, these conditions are satisfied in the case of tensor Hilbert spaces equipped with a canonical norm, and also for $L^p$ spaces.

\subsection{Optimization problems in subsets of low-rank tensors}\label{sec:optimization-low-rank-subsets}

Standard subsets of low-rank tensors $\Sc_r$ (such as $\Rc_r$, $\Tc_r$, $\Bc\Tc_r$ or $\Tc\Tc_r$) are not vector spaces nor convex sets. Therefore, the solution of a best approximation problem in $\Sc_r$,
or more generally of an optimization problem
\begin{align}
\min_{v\in \Sc_r} J(v),\label{optimSr}
\end{align}
with $J: X_{\Vert \cdot\Vert} \rightarrow \Rbb$, requires ad-hoc minimization algorithms.
Standard subsets of low-rank tensors admit a parametrization of the form
\begin{align}
\Sc_r = \{ v = F_{\Sc_r}(p_1,\hdots,p_M): p_i \in P_i, 1\le i\le M\},\label{multilinear-param-Sr}
\end{align}
where $F_{\Sc_r} : P_1 \times \hdots \times P_M \to X$ is a multilinear map and the $P_i$ are vector spaces or standard submanifolds of vector spaces (e.g. Stiefel manifolds) (see Remarks \ref{rem:param-Rr}, \ref{rem:param-Tr}, \ref{rem:param-BTr} and \ref{rem:param-TTr} respectively for $\Rc_r$, $\Tc_r$,  $\Bc\Tc_r$ and $\Tc\Tc_r$). The optimization problem \eqref{optimSr} is then rewritten as an optimization problem on the parameters
$$
\min_{p_1 \in P_1,\hdots,p_M\in P_M} J(F_{\Sc_r}(p_1,\hdots,p_M)),
$$
which allows the use of more or less standard optimization algorithms (e.g. Newton, steepest descent, block coordinate descent), possibly exploiting the manifold structure of $P_1 \times \hdots \times P_M$ (see e.g. \cite{Espig:2012kx,vandereycken2013low,Uschmajew:2012fk}).
Alternating minimization algorithms (or block coordinate descent algorithms)  transform the initial optimization problem into a succession of simpler optimization problems. They
consist in solving successively the minimization problems
$$
\min_{p_i \in P_i} J(F_{\Sc_r}(p_1,\hdots,p_M)),
$$
each problem being a minimization problem in a linear space (or standard manifold)  $P_i$ of a functional $p_i \mapsto J(F_{\Sc_r}(p_1,\hdots,p_M))$ which inherits from some properties of the initial functional $J$ (due to the linearity of the partial map $p_i\mapsto F_{\Sc_r}(p_1,\hdots,p_M)$ from $P_i$ to $X$). The available convergence results  for these optimization algorithms in a general setting only ensure local convergence or global convergence to critical points (see e.g. \cite{Rohwedder:2013fk,2015arXiv150600062E}).

\subsection{Higher-order singular value decomposition}\label{sec:hosvd}
The \emph{Higher-Order Singular Value Decomposition} (HOSVD), introduced in \cite{DEL00} for the Tucker format, in \cite{GRA10} for the Hierarchical Tucker format, and in \cite{OSE11} for the TT-format, constitutes a possible generalization of the SVD for tensors of order $d \ge 3$ which allows us to obtain quasi-best approximations (but not necessarily best approximations) in subsets of low-rank tensors (for tree-based Tucker formats). It relies on the use of the SVD for order-two tensors applied to matricisations of a tensor. 
Here, we consider a tensor Hilbert space $X$ equipped 
with the canonical norm $\Vert\cdot\Vert$.
For each nonempty subset $\alpha\subset D$, $X_\alpha = \bigotimes_{\nu\in \alpha} X_\nu$ is also equipped with the canonical norm, denoted $\Vert \cdot\Vert_\alpha$. 
\par
Let us consider an element $u$ in the algebraic tensor space\footnote{The case where $u\in X_{\Vert\cdot\Vert}\setminus X$ introduces some technical difficulties related to the definition of tree-based topological tensor spaces (see \cite{falco2015geometric}).}  $X$.
For $\alpha\subset D$, let  $u_{\alpha,r_\alpha} \in X$ 
denote the best approximation of $u$ with $\alpha$-rank bounded by $r_\alpha$, i.e.
\begin{align*}
\Vert u- u_{\alpha,r_\alpha} \Vert &= \min_{\rank_\alpha(v) \le r_\alpha}  \Vert u-v\Vert.
\end{align*}
$u_{\alpha,r_\alpha} $ is such that $\Mc_\alpha(u_{\alpha,r_\alpha} )$ is the rank-$r_\alpha$ truncated SVD of $\Mc_\alpha(u) \in X_\alpha \otimes X_{\alpha^c}$, which can be written
$$
u_{\alpha,r_\alpha}   = \sum_{i=1}^{r_\alpha} \sigma_{i}^{(\alpha)} u_i^{(\alpha)} \otimes u_i^{(\alpha^c)},
$$
where $\sigma_i^{(\alpha)}$ are the dominant singular values and $u_i^{(\alpha)} $ and $u_i^{(\alpha^c)} $ the corresponding left and right singular vectors of $\Mc_\alpha(u)$. Let $U^{(\alpha)}_{r_\alpha}= U^{min}_{\alpha}(u_{\alpha,r_\alpha} ) = \mathrm{span}\{u_i^{(\alpha)}\}_{i=1}^{r_\alpha} $ denote the resulting optimal subspace in $X_\alpha$ and $P_{U^{(\alpha)}_{r_\alpha}}$ the corresponding orthogonal projection from  $X_{\alpha} $ to $U^{(\alpha)}_{r_\alpha}$
(associated with the canonical inner product in $X_{\alpha}$). The projection is such that $u_{\alpha,r_\alpha}  = (P_{U^{(\alpha)}_{r_\alpha}} \otimes id_{\alpha^c}) (u)$. We note that $\{U_{r_\alpha}^{(\alpha)}\}_{r_\alpha \ge 1}$ is an increasing sequence of subspaces. We have the  orthogonal decomposition $U_{r_\alpha}^{(\alpha)} = \bigoplus_{i_\alpha=1}^{r_\alpha} W^{(\alpha)}_{i_\alpha}$ with $W^{(\alpha)}_{i_\alpha} = \mathrm{span}\{u_{i_\alpha}^{(\alpha)}\}$, and $P_{U_{r_\alpha}^{(\alpha)}} = \sum_{i_\alpha=1}^{r_\alpha} P_{W_{i_\alpha}^{(\alpha)}}.$

\subsubsection{HOSVD in Tucker format}
Let $r=(r_1,\hdots,r_d) \in \Nbb^d$ such that $r_\nu \le \rank_\nu(u)$ for $1\le \nu\le d$.
For each dimension  $\nu \in D$, we define the optimal $r_\nu$-dimensional space $U_{r_\nu}^{(\nu)}$ and the corresponding  orthogonal projection  $P_{U_{r_\nu}^{(\nu)}}$. Then, we define the space $U_r = \bigotimes_{\nu=1}^d U^{(\nu)}_{r_\nu}$ and the associated orthogonal projection $$P_{U_r} = P_{U_{r_1}^{(1)}}\otimes \hdots \otimes P_{U_{r_d}^{(d)}}.$$
Then, the truncated HOSVD of $u$ with multilinear rank $r$ is defined by 
$$
u_{r} = P_{U_r} (u) \in \Tc_r.
$$
We note that subspaces $\{U_r\}_{r\in \Nbb^d}$ are nested: for $s,r\in \Nbb^d$ such that $s\ge r$, we have $U_{r}\subset U_s$.
The approximation $u_r$ can be obtained by truncating a decomposition of $u$. Indeed, noting that $U_r = \bigoplus_{i\le r} W_i$, with $W_i = \bigotimes_{\nu \in D} W_{i_\nu}^{(\nu)}$, we have
$$
u_r = \sum_{i\le r} w_i, \quad w_i=  P_{W_i} (u),
$$
which converges to 
$
u 
$ when $r_\nu \to \rank_\nu(u)$ for all $\nu$.
We have that $u_r$ is a quasi-optimal approximation of $u$ in $\Tc_r$ (see \cite[Theorem 10.3]{HAC12}), such that
\begin{align*}
\Vert u-u_r\Vert \le \sqrt{d} \min_{v\in \Tc_r} \Vert u-v \Vert.
\end{align*} 

\begin{rem}
Another version of the HOSVD can be found in \cite[Section 10.1.2]{HAC12}, where the spaces $U_{r_\nu}^{(\nu)}$, $1\le \nu \le d$, are computed successively. The space $U_{r_\nu}^{(\nu)}$ is defined as the dominant singular space of  
$\Mc_{\{\nu\}}(u^{(\nu-1)})$, with $u^{(\nu-1)} = P_{U_{r_1}^{(1)}}\otimes \hdots \otimes P_{U_{r_{\nu-1}}^{(\nu-1)}} u.$
\end{rem}

\subsubsection{HOSVD in tree-based Tucker format}
Let $T_D$ be a dimension tree and $r = (r_\alpha)_{\alpha\in T_D}$ be an admissible set of ranks, with $r_\alpha \le \rank_\alpha(u)$ for all $\alpha\in T_D$. For each vertex $\alpha\in T_D$, we define the optimal $r_\alpha$-dimensional subspace $U_{r_\alpha}^{(\alpha)} \subset X_\alpha$ and the associated projection $P_{U_{r_\alpha}^{(\alpha)}}$. Let $P_{r_\alpha}^{(\alpha)} = P_{U^{(\alpha)}_{r_\alpha}} \otimes id_{\alpha^c}$. Then the truncated HOSVD of $u$ with $T_D$-rank $r$ is defined by  
$$
u_{r} = P^{T_D}_r (u) \in \Bc\Tc_r,
$$
with 
$$
 P^{T_D}_r (u) = P^{T_D,(L)}_{r}  P^{T_D,(L-1)}_{r} \hdots  P^{T_D,(1)}_r,\quad  P^{T_D,(\ell)}_{r} = 
\prod_{\substack{\alpha \in T_D \\ level(\alpha)=\ell}} P_{r_\alpha}^{(\alpha)},
$$
where $level(\alpha)$ is the {level} of a vertex in the tree, with $level(D)=0$, and where $L=\max_{\alpha	\in T_D} level(\alpha)$.
 We have that $u_r$ is a quasi-optimal approximation of $u$ in $\Bc\Tc_r$ (see \cite[Theorem 11.58]{HAC12}), such that
\begin{align*}
\Vert u-u_r\Vert \le \sqrt{2d-2-s} \min_{v\in \Bc\Tc_r} \Vert u-v \Vert,
\end{align*} 
with $s=1$ if $\#S(D)=2$ and $s=0$ if $\#S(D)>2$. 
\begin{rem}
For the TT-format, the truncated HOSVD of $u$ with TT-rank $r=(r_1,\hdots,r_{d-1})$ is defined by 
$
u_r = P_{r_{d-1}}^{(\{d\})} \hdots P_{r_1}^{(\{2,\hdots,d\})} (u)$, where $ P_{r_{k}}^{(\{k+1,\hdots,d\})}$ is the orthogonal projection associated with the optimal $r_k$-dimensional subspace $U^{\{k+1,\hdots,d\}}_{r_k}$ in  $X_{\{k+1,\hdots,d\}}$, $1\le k\le d-1$ (no projection associated with vertices $\{k\}$, $1\le k \le d-1$). We have that $\Vert u-u_r \Vert\le \sqrt{d-1} \min_{v\in \Tc\Tc_r} \Vert u-v \Vert $.
\end{rem}
\begin{rem}
Other versions of HOSVD for tree-based formats can be found in \cite[Sections 11.4.2.2 and 11.4.2.3]{HAC12}, where the spaces $U_{r_\alpha}^{(\alpha)}$, $\alpha\in T_D$, are computed successively.\end{rem}

\section{Greedy algorithms for low-rank approximation}\label{sec:greedy}

It can be observed in many practical applications that best approximations in low-rank tensor formats present good convergence properties (with respect to the rank). However, the computational complexity for computing best approximations drastically increases with the rank. Also, in general, the sequence of best approximations of a tensor  is not associated with a decomposition of the tensor, that means that best approximations can not be obtained by truncating a decomposition of the tensor. In Sections \ref{sec:svd} and \ref{sec:hosvd}, we have seen that the SVD or one of its extensions for higher-order tensors
allows recovering such a notion of decomposition. However, it is restricted to the approximation of a tensor in a tensor Hilbert space equipped with the canonical norm, and it requires an explicit representation of the tensor.

Greedy algorithms (sometimes called Proper Generalized Decomposition methods) aim at recovering a notion of decomposition,  by relying either   on greedy constructions of the approximation (by computing successive corrections in subsets of low-rank tensors) or on greedy constructions of subspaces (for subspace-based low-rank formats).  These algorithms are applied in a more general setting where one is interested in constructing low-rank approximations $w$ that minimize some distance $\Ec(u,w)
$ to a tensor $u$.   These constructions, although they are suboptimal, allow reducing the computational complexity for the computation of high rank  approximations and they sometimes achieve quasi-optimal convergence (with the rank). These quasi-optimality properties are observed in some practical applications but they still require a theoretical justification. 

\begin{rem}\label{rem:greedy_svd}
Note that in the particular case where $X=S\otimes V$ with $V$ and $S$ Hilbert spaces, and $ \Ec(u,w) = \Vert u-w\Vert$ with $\Vert\cdot\Vert$ the canonical norm, all the algorithms presented in this section yield  the singular value decomposition of $u$ (provided that successive minimization problems are solved exactly). In general, when deviating from this particular case,  the presented algorithms  yield  different decompositions. 
\end{rem}

\subsection{Greedy construction of the approximation}\label{sec:greedy-higher-order}

A natural way to recover a notion of tensor decomposition is to define a sequence of approximations with increasing canonical rank obtained by successive rank-one corrections. This algorithm constitutes the most prominent version of so-called Proper Generalized Decomposition and it has been used in many applications (see the review \cite{CHI11} and the monograph \cite{Chinesta_book:2014}). 
 Starting from $u_0=0$, a rank-$r$ approximation $u_r \in \Rc_r$ is defined by 
$$
u_r = u_{r-1} + w_r,
$$
where $w_r = \otimes_{\nu=1}^d w^{(\nu)}_r \in \Rc_1$ is the optimal rank-one correction of $u_{r-1}$ such that 
\begin{align}
\Ec( u,u_{r-1} + w_r)  = \min_{w\in \Rc_1}  \Ec( u,u_{r-1} +w).\label{eq:rankonecorrection_orderd} 
\end{align}
This algorithm can be interpreted as a greedy algorithm in the dictionary of rank-one tensors $\Rc_1$ in $X_{\Vert\cdot\Vert}$ and it allows recovering a notion of decomposition, even for higher-order tensors. Indeed, assuming that the sequence $\{u_r\}_{r\ge 1}$ strongly converges to $u$, then $u$ admits the decomposition 
\begin{align}
u = \sum_{i\ge 1} w_i^{(1)} \otimes \hdots \otimes w_i^{(d)},\label{eq:decompucanonical}
\end{align}
and the approximation $u_r$ with canonical rank $r$ can be obtained by truncating this series after $r$ terms, therefore justifying the notion of decomposition.   
When $\Ec(u,w)=\Vert u-w\Vert$, conditions for the convergence of greedy algorithms in a general setting can be found in  \cite{TEM11}. In the case of the minimization of convex functionals, convergence results can be found in \cite{LEB09,CAN11,Falco:2012fk,temlyakov2012greedy}. 
 Note that this greedy construction is not specific to the particular setting of tensor approximation. The available convergence results do not take into account any particular structure of the tensor $u$ and are usually pessimistic. 
However, except for very particular cases (see Remark \ref{rem:greedy_svd}), 
 this algorithm only provides a suboptimal sequence of rank-$r$ approximations. Depending on the properties of $\Ec(u,w)$, the convergence with the rank $r$ may be strongly deteriorated by this greedy construction, compared with the best approximation error in canonical format, that is
 $\sigma(u;\Rc_r) = \inf_{v\in \Rc_r}  \Ec(u ,v )$ (which corresponds to the error of best $r$-term approximation in the dictionary $\Rc_1$).

 A classical improvement of the above construction (known as orthogonal greedy algorithm) consists in first computing a rank-one correction $w_r = \bigotimes_{\nu=1}^d w_r^{(\nu)}$ by solving \eqref{eq:rankonecorrection_orderd}  and then (after a possible normalization of $w_r$) in defining 
$$
u_r = \sum_{i=1}^r \sigma^{(r)}_i w_i,
$$
where the set of coefficients $(\sigma_{i}^{(r)})_{i=1}^r$ is solution of
$$
\Ec( u , u_r) = \min_{(\sigma_i^{(r)})_{i=1}^r \in \Rbb^r} \Ec( u, \sum_{i=1}^r \sigma^{(r)}_i  w_i ).
$$
In many applications, it is observed that this additional  step does not significantly improve the convergence of the sequence $u_r$. 
\begin{rem}
In the orthogonal greedy construction, the approximation $u_r$ cannot be obtained by truncating a decomposition of the form \eqref{eq:decompucanonical}, and therefore, the sequence $u_r$ has to be interpreted as a decomposition in a general sense. 
\end{rem}
The orthogonal greedy algorithm has been analyzed in 
\cite{Falco:2012fk} as a particular case of a family of algorithms using more general dictionaries of low-rank tensors, and using improvement strategies that are specific to the context of low-rank tensor approximation. In fact, improvements that seem to be efficient in practice do not rely anymore on greedy approximations, but rather adopt a subspace point of view in which low-rank corrections are only used for the greedy construction of subspaces. This requires to move to other tensor formats, as presented in the next section. 
\begin{rem}
Note that the above algorithms define sequences of spaces $U_r^{(\nu)} = \mathrm{span}\{w_1^{(\nu)},\hdots,w_r^{(\nu)}\}$ in $X_\nu$ verifying the nestedness property $U_{r}^{(\nu)} \subset U_{r+1}^{(\nu)} $ and such that $u_r\in U_r^{(1)}\otimes \hdots \otimes U_r^{(d)}$. However, the algorithms do not exploit this subspace point of view. 
\end{rem}

\subsection{Greedy construction of subspaces for order-two tensors}\label{sec:greedy-subspaces}
 
 When using subspace-based tensor formats, 
other notions of decompositions can be obtained by defining a sequence of approximations in increasing tensor spaces. Here, we present algorithms for the approximation of an order-two tensor $u$ in $X_{\Vert\cdot\Vert}$ with $X=S\otimes V$. Their extensions to the case higher-order tensors are presented in Section \ref{sec:greedy-subspaces-higher-order}.

 \subsubsection{Fully greedy construction of subspaces}\label{sec:greedy-subspace-order-two}
For an  order-two tensor $u$, the best rank-$r$ approximation problems \eqref{eq:bestrankrsubspace1}, $r\ge 1$,  yield a sequence of rank-$r$ approximations $$
u_r = \sum_{i=1}^r  s_i^{(r)} \otimes v_i^{(r)}.
$$ 
 The associated sequences of reduced approximation spaces $S_r = U^{min}_1(u_r) = \mathrm{span}\{s_i^{(r)}\}_{i=1}^r$ and $V_r = U^{min}_2(u_r) = \mathrm{span}\{v_i^{(r)}\}_{i=1}^r$, such that 
 \begin{align}\label{eq:urinSrVr}
u_r \in S_r \otimes V_r, 
\end{align}
do not necessarily satisfy \begin{align}\label{eq:nested_spaces}
S_r \subset S_{r+1} \quad \text{and} \quad V_r \subset V_{r+1}.\end{align}

A notion of decomposition can be obtained by defining a sequence of rank-$r$ approximations $u_r$ in an increasing sequence of subspaces $S_r\otimes V_r$, which means such that minimal subspaces 
$S_r = U_{1}^{min}(u_r)$ and $V_r = U_{2}^{min}(u_r)$ verify the nestedness property \eqref{eq:nested_spaces}. The resulting approximation $u_r$ is defined as the best approximation in $S_r \otimes V_r$, i.e.
\begin{align}
\Ec( u ,u_r )  = \min_{w\in S_r\otimes V_r} \Ec(u,w), \label{eq:approxSrotimesVr}
\end{align}
and can be written under the form 
$$
u_r = \sum_{i=1}^r \sum_{j=1}^r \sigma_{ij}^{(r)} s_i \otimes v_j,
$$
where $\{s_i\}_{i=1}^r$ and $\{v_i\}_{i=1}^r$ are bases of $S_r$ and $V_r$ respectively, and where 
 $\sigma^{(r)}\in \Rbb^{r\times r}$ is the solution of 
 $$
 \min_{\sigma^{(r)}\in \Rbb^{r\times r}} \Ec( u ,   \sum_{i=1}^r \sum_{j=1}^r \sigma_{ij}^{(r)} s_i \otimes v_j).
 $$
Different constructions of nested subspaces can be proposed.
\\

\paragraph{Optimal  construction with nested minimal subspaces.}
A first and natural definition of $u_r$ is such that 
\begin{align*}
\Ec( u, u_r)  = \min_{\substack{S_r \in \Gbb_r(S) \\ S_r \supset S_{r-1}}} \min_{ \substack{ V_r\in \Gbb_r(V)\\ V_r\supset V_{r-1}}} \min_{w\in S_r\otimes V_r} \Ec(u,w),
\end{align*}
which corresponds to the definition \eqref{eq:bestrankrsubspace1} of optimal rank-$r$ approximations with the only additional constraint that minimal subspaces of successive approximations are nested.  This definition can be equivalently written in terms of the new elements $s_r\in S$ and $v_r \in V$ and of the matrix of coefficients $\sigma^{(r)}\in \Rbb^{r\times r}$:
\begin{align}
\Ec( u, u_r) = \min_{s_r \in S} \min_{v_r\in V} \min_{\sigma^{(r)} \in \Rbb^{r\times r}} \Ec(  u  ,  \sum_{i=1}^r \sum_{j=1}^r \sigma_{ij}^{(r)} s_i \otimes v_j).\label{eq:optimalgreedysubspaces}
\end{align}
 
\paragraph{Suboptimal construction.}
A simpler but suboptimal construction (compared to \eqref{eq:optimalgreedysubspaces}) consists in defining the new elements $s_r\in S$ and $v_r\in V$ by computing an optimal rank-one correction of the previous approximation $u_{r-1}$. More precisely, given  $u_{r-1} = \sum_{i=1}^{r-1} \sum_{j=1}^{r-1} \sigma_{ij}^{(r-1)} s_i \otimes v_j$, $s_r\in S$ and $v_r\in V$ are defined by
\begin{align*}
 \min_{s_r \in S} \min_{v_r\in V}  \Ec( u ,  u_{r-1}  + s_r\otimes v_r),
\end{align*}
and then the approximation $u_r$ is obtained by solving \eqref{eq:approxSrotimesVr} with spaces 
$S_r = S_{r-1} + \mathrm{span}\{s_r\}$ and $V_r =  V_{r-1} + \mathrm{span}\{v_r\}$.

\subsubsection{Partially greedy construction of subspaces}\label{sec:partial_greedy_subspace}

Another notion of decomposition can be obtained by  imposing the nestedness property for only one of the minimal subspaces, say $V_r\subset V$, which results in a sequence $u_r$ of the form 
$$
u_r = \sum_{i=1}^r s_i^{(r)} \otimes v_i.
$$
This is a non-symmetric point of view which focuses on the construction of reduced spaces in $V$. This point of view is of particular interest in the case of Bochner spaces (see Section \ref{sec:low-rank-bochner}), for the model reduction of parameter-dependent or stochastic equations (see Section {\ref{sec:Proper Generalized Decompositions for parametric equations in tensor format}} and references \cite{NOU07,NOU08b,NOU09b,CHE12,Tamellini:2014}), and also for the model reduction of time-dependent evolution equations (see \cite{LAD99b,LAD10,NOU10}).

\paragraph{Optimal construction with nested minimal subspaces.} The sequence of rank-$r$ approximations $u_r$ can  be defined by
\begin{align}
\Ec( u, u_r)  =  \min_{\substack{V_r\in \Gbb_r(V) \\ V_r\supset V_{r-1}}} \min_{w\in S\otimes V_r} \Ec(u,w). \label{eq:optimal-partial-greedy}
\end{align} 
This definition corresponds to the definition \eqref{eq:bestrankrsubspace1} of optimal rank-$r$ approximations with the only additional constraint that the minimal subspaces $V_r = U^{min}_2(u_r)$ are nested. It is equivalent to the minimization problem
\begin{align}
 \min_{v_r\in V} \min_{\{s_i^{(r)}\}_{i=1}^r \in S^r}  \Ec( u,\sum_{i=1}^r s_i^{(r)} \otimes v_i ),\label{eq:optimal-partial-greedy-functions}
\end{align}
which can be solved  by an alternating minimization algorithm {(see Section \ref{sec:Proper Generalized Decompositions for parametric equations in tensor format} for the application to parameter-dependent equations)}. 

\paragraph{Suboptimal construction.}
Suboptimal constructions can also be introduced in order to reduce the computational complexity, e.g. by computing a rank-one correction of $u_{r-1}$ defined by $\min_{v_r\in V} \min_{s_r \in S } \Ec( u, u_{r-1} + s_r\otimes v_r )$, and then by solving 
$\Ec( u, u_r)  = \min_{w\in S\otimes V_r} \Ec( u,w)$ with $V_r= V_{r-1} + \mathrm{span}\{v_r\}$. 

\subsubsection{Partially greedy construction of subspaces in Bochner spaces}\label{sec:partial_greedy_parametric}
Let $X=L^p_\mu(\Xi)\otimes V$ and let $\Vert \cdot \Vert_p$ denote the Bochner norm.
For $1<p<\infty$ (and in particular for $p=2$), we can consider the algorithm presented in Section \ref{sec:partial_greedy_subspace}, with $\Ec(u,v)=\Vert u-v \Vert_p$. It defines the rank-$r$ approximation $u_r$ by
$$
\Vert u-u_r\Vert_p = \min_{\substack{V_r\in \Gbb_r(V) \\ V_r \supset V_{r-1}}}  \min_{w\in S\otimes V_r} \Vert  u-w\Vert_p
=  \min_{\substack{V_r\in \Gbb_r(V) \\ V_r \supset V_{r-1}}} \Vert  u-P_{V_r} u\Vert_p,
$$ 
which is a well-posed optimization problem (as a best approximation problem in a weakly closed subset of the reflexive Banach space 
$L^p_\mu(\Xi;V)$, see Sections \ref{sec:bestapproxproblem} and \ref{sec:Bochner}).
This algorithm generates an increasing sequence of reduced approximation spaces $V_r$ that are optimal in a ``$L^p$ sense''. 
For $p=\infty$, an ideal greedy construction would define $u_r = P_{V_r} u$ with $V_r$ solution of 
$$
 \inf_{\substack{V_r\in \Gbb_r(V) \\ V_r \supset V_{r-1}}} \Vert  u-P_{V_r} u\Vert_\infty = \inf_{\substack{V_r\in \Gbb_r(V) \\ V_r \supset V_{r-1}}} \esssup_{y\in \Xi}\Vert  u(y)-P_{V_r} u(y)\Vert_V.
 $$
 Suboptimal constructions can be proposed in order to avoid computational issues related to the optimization with respect to the $L^\infty$-norm. 
Suppose that $u:\Xi\rightarrow V$ is continuous and $\Xi= \mathrm{support}(\mu)$ is compact. Then, starting from $V_0=0$, one defines  $V_r = V_{r-1} + \mathrm{span}\{v_r\}$ with $v_r\in V$ such that 
\begin{equation}
\sup_{y\in \Xi}\Vert u(y)-P_{V_{r-1}} u(y) \Vert_V = \Vert v_r - P_{V_{r-1}} v_r \Vert_V.\label{greedy_linfini}
\end{equation}
This is the greedy construction used in the Empirical Interpolation Method \cite{MAD09}.
Convergence results for this algorithm can be found in \cite{Binev:2011,buffa:2012,DeVore:2013fk}, where 
the error $\Vert u-u_r \Vert_\infty = \sup_{y\in \Xi} \Vert u(y) - P_{V_r} u(y)\Vert_V$ is compared with the best rank-$r$ approximation error 
$\rho^{(\infty)}_r(u) = d_r(u(\Xi))_V$.

\subsection{Greedy construction of subspaces for higher-order tensors}\label{sec:greedy-subspaces-higher-order}

Here we extend the constructive algorithms presented in Section \ref{sec:greedy-subspaces} to the case of higher-order subspace-based tensor formats. 

\subsubsection{Greedy construction of subspaces for the Tucker format}\label{sec:greedy-subspace-tucker}

The algorithms presented in Section \ref{sec:greedy-subspace-order-two} can be naturally generalized in order to provide constructive algorithms for the approximation of tensors in Tucker format. These algorithms construct a sequence of approximations $u_m$ in nested tensor spaces
$U_m = U_{m}^{(1)}\otimes \hdots \otimes U_m^{(d)}$, with $U^{(\nu)}_{m} \subset U^{(\nu)}_{m+1}$, therefore allowing a notion of decomposition to be recovered. 

\paragraph{Construction of subspaces based on rank-one corrections.}
A first strategy, introduced in \cite{Giraldi2013154}, consists in progressively enriching the spaces by the factors of rank-one corrections. More precisely, we start with $u_0=0$. Then, for $m\ge 1$, we compute a rank-one correction $w_m = \bigotimes_{\nu=1}^d w_m^{(\nu)} \in \Rc_1$ of $u_{m-1}$, which is solution of $$
\Ec (u,u_{m-1} + w_m) = \min_{w\in\Rc_1} \Ec( u,u_{m-1} + w ),
$$
and then define $U_{m}^{(\nu)}  = \mathrm{span}\{w_{i}^{(\nu)}\}_{i=1}^{m}$, for all $\nu\in D$. Then, $u_m \in U_m = \bigotimes_{\nu=1}^d U^{(\nu)}_{m} $ is defined by
$$
\Ec( u,u_m ) = \min_{v\in U_m} \Ec (u,v),
$$
and can be written
$$
u_m = \sum_{i_1=1}^{m} \hdots \sum_{i_d=1}^m \sigma_{i_1,\hdots,i_d}^{(m)} \bigotimes_{\nu=1}^d w_{i_\nu}^{(\nu)}.
$$
This construction has also been applied for the construction of an approximate inverse of an operator in low-rank format \cite{Giraldi:2014}. 
  For some 
  applications (see \cite{Giraldi2013154,Giraldi:2014}), when $\Ec(u,v)\sim \Vert u-v\Vert$, we observe an error $\Ec( u ,u_m )$ which behaves as the best approximation 
 error in Tucker format $\sigma(u;\Tc_{r^{(m)}}) = \inf_{v\in \Tc_{r^{(m)}}} \Vert u-v\Vert$ with $r^{(m)}=(m,\hdots,m)$. The theoretical justification of these observations remains an open problem. The above construction is isotropic in the sense that subspaces are enriched in all directions $\nu\in D$ simultaneously. This does not allow us to take advantage of possible anisotropic structures of the tensor $u$. 
\begin{rem}
Of course,  computing an approximation in the tensor product space $U_m$ is not tractable in high dimension $d$ without additional complexity reduction techniques. In \cite{Giraldi:2014}, it is proposed to approximate  $u_m$ in a low-rank hierarchical tensor format in the tensor space $U_m$.
\end{rem}

\paragraph{Optimal greedy construction of subspaces.}
Another natural algorithm consists in simply adding the nestedness property of subspaces in the definition \eqref{eq:bestapproxTr_subspace} of best approximations. 
Starting from $u_0=0$, we let $u_{m-1}$ denote the approximation at step $m-1$ of the construction and $U_{m-1}^{(\nu)} = U^{min}_\nu(u_{m-1})$ 
, for $\nu \in D$. At step $m$,  we select a set of dimensions $D_m \subset D$ to be enriched, we let $U_m^{(\nu)} = U_{m-1}^{(\nu)}$ for $\nu\notin D_m$ and we define $u_m$ by
\begin{align}
\Ec(u,u_m)
= \min_{\substack{(U^{(\nu)}_{m})_{\nu\in D_m} \\ 
\dim(U^{(\nu)}_m) = \dim(U^{(\nu)}_{m-1}) + \Delta r_\nu^{(m)} \\
U^{(\nu)}_m \supset U^{(\nu)}_{m-1}}}
 \min_{v\in U_m^{(1)}\otimes \hdots \otimes U_m^{(d)}} \Ec(u,v). \label{eq:approxTr_subspace_greedy_anisotropic}
\end{align} 
Choosing $D_m=D$ and $\Delta r_\nu^{(m)}=1$ for all $\nu \in D$ at each step corresponds to an isotropic enrichment (similar to the previous construction based on rank-one corrections).
However, this isotropic construction does not allow any particular structure of the tensor $u$ to be exploited. Choosing $D_m \neq D$ or different values for the $\Delta r_\nu^{(m)}$, $\nu \in D$, yields anisotropic constructions but the selection of $D_m$ and $\Delta r_\nu^{(m)}$, $\nu\in D$, requires the introduction of some error indicators.
This type of construction seems to provide good convergence properties with respect to the rank $r^{(m)} = (r_\nu^{(m)})_{\nu\in D}$, with $r_\nu^{(m)} = \dim(U^{(\nu)}_{m})$. However, it remains an open and challenging question to prove that this type of construction can achieve quasi-optimality compared to the best rank-$r^{(m)}$ approximation for certain classes of functions (e.g. associated with a certain decay of the best rank-$r^{(m)}$ approximation error).

 \subsubsection{Greedy construction of subspaces for the tree-based tensor format}

The construction presented in Section \ref{sec:greedy-subspace-tucker} can be extended to more general tree-based Tucker formats, these formats being related to the notion of subspaces. The idea is again to start from the subspace-based formulation of the best approximation 
problem in tree-based Tucker format \eqref{eq:bestapproxBTr_subspace}, and to propose a suboptimal greedy construction of subspaces which consists in adding a nestedness property for the successive  minimal subspaces.  
We start from $u_0=0$. Then we let $u_{m-1}$ denote the approximation at step $m-1$ of the construction and we let $U_{m-1}^{(\alpha)} = U^{min}_\alpha(u_{m-1})$ denote the current minimal subspaces of dimensions $r^{(m-1)}_\alpha = \dim(U^{min}_\alpha(u_{m-1}))$, $\alpha \in T_D$. Then, at step $m$, we select a set of vertices $T_m \subset T_D$ and we define $r^{(m)} = (r^{(m)}_\alpha)_{\alpha \in T_D}$ with $r^{(m)}_\alpha = r_\alpha^{(m-1)} + \Delta r_\alpha^{(m)} $ for $\alpha\in T_m$ and $r^{(m)}_\alpha = r_\alpha^{(m-1)} $ for $\alpha \in T_D\setminus T_m$. Then, we define $u_m$ as the solution of 
 \begin{align}
\Ec(u,u_m)
= \min_{\substack{(U^{(\alpha)}_{m})_{\alpha\in T_D\setminus D}  \in \Gc_{r^{(m)}}({T_D})\\ U^{(\alpha)}_{m}\supset U^{(\alpha)}_{m-1}}}  \min_{v \in    \bigotimes_{\alpha\in S(D)} U^{(\alpha)}_m}\Ec(u,v).  
\end{align} 
The selection of vertices $T_m$ and of the $\Delta r_\alpha^{(m)}$, $\alpha\in T_m$, requires the introduction of error indicators and strategies of enrichment (preserving admissibility of $T_D$-rank). This type of construction seems to be a good candidate for really exploiting specific tensor structures but the analysis and the implementation of this type of strategy remain open and challenging issues.

\subsection{Remarks on the solution of minimization problems}
Constructive algorithms presented in this section require the solution of successive minimization problems in  subsets which are not vector spaces nor convex sets. In practice, one can rely on standard optimization algorithms by exploiting a multilinear parametrization of these approximation subsets (see Section \ref{sec:optimization-low-rank-subsets} for the optimization in standard low-rank manifolds, e.g. the set of rank-one tensors $\Rc_1$). As an illustration for a non standard subset introduced in the present section, let us consider the solution of  
\eqref{eq:optimalgreedysubspaces}, which is written 
$$
\min_{s_r\in S}\min_{v_r\in V}\min_{\sigma^{(r)}\in \Rbb^{r\times r}} J(s_r,v_r,\sigma^{(r)}), 
$$ 
with $J(s_r,v_r,\sigma^{(r)}) =  \Ec( u , \sum_{i=1}^r \sum_{j=1}^r \sigma_{ij}^{(r)} s_i\otimes v_j)$.
A natural alternating minimization algorithm then consists in successively solving minimization problems 
$$
\min_{s_r\in S} J(s_r,v_r,\sigma^{(r)}), \quad 
\min_{v_r\in V} J(s_r,v_r,\sigma^{(r)}) \quad \text{and} \quad 
\min_{\sigma^{(r)} \in \Rbb^{r\times r}} J(s_r,v_r,\sigma^{(r)}).
$$
Note that in practice, algorithms do not yield exact solutions of optimization problems. 
The analysis of constructive algorithms presented in this section should therefore take into account  these approximations and quantify their impact. Several convergence results are available for weak greedy algorithms \cite{TEM11,temlyakov2012greedy}, which are perturbations of ideal greedy algorithms presented in Section \ref{sec:greedy-higher-order}. 

\section{Low-rank approximation using samples}\label{sec:low-rank-samples}

In this section, we present methods for the practical construction of low-rank approximations of a vector-valued or multivariate function (identified with a tensor) from sample evaluations of the function.  

\subsection{Low-rank approximation of vector-valued functions}
\label{sec:low-rank-sample-norm}
Let $u : \Xi \to V$ be a vector-valued function, with $V$ a Banach space and $\Xi$ a set equipped with a measure $\mu$, and let us assume that $u \in L^p_\mu(\Xi;V)$.
A low-rank approximation of  $u$ can be defined from sample evaluations of $u$.
Let $\Xi_K = \{y^k\}_{k=1}^K$ be a set of sample points in $\Xi$ (e.g. samples drawn according to a probability measure $\mu$ on $\Xi$). Then, for $w\in L^p_\mu(\Xi;V)$,
we define 
\begin{align}
&\Vert w\Vert_{\infty,K} = \sup_{1\le k \le K} \Vert w(y^k) \Vert_V, \; \text{and}\label{infinite-K-norm} \\
&\Vert w \Vert_{p,K} = \big(
\sum_{k=1}^K \omega^k \Vert  w(y^k) \Vert_V^p\big)^{1/p}\quad 
\text{for } p<\infty,\label{p-K-norm}
\end{align}
where $\{\omega^k\}_{k=1}^K$ is a set of positive weights. For $p<\infty$, if the $y^k$ are i.i.d. samples drawn according the probability measure $\mu$ and if $\omega^k = K^{-1}$ for all $k$, then $\Vert w\Vert_{p,K}$ is a statistical estimate of the Bochner norm $\Vert w\Vert_p$.
For $1\le p\le \infty$, the application $v \mapsto \Vert v\Vert_{p,K}$ defines a semi-norm on 
 $L^p_\mu(\Xi;V)$. An optimal rank-$r$ approximation $u_r$ of $u$ with respect to 
 the semi-norm  $\Vert \cdot\Vert_{p,K}$ is defined by
 $$
  \Vert u-u_r \Vert_{p,K} = \min_{w\in \Rc_r} \Vert u-w \Vert_{p,K} := \rho^{(p,K)}_r(u),
 $$
 or equivalently by
  \begin{equation}
  \Vert u-u_r \Vert_{p,K}  = \min_{V_r\in \Gbb_r(V)} \Vert u-P_{V_r} u\Vert_{p,K},\label{bestrankrseminorm}
 \end{equation}
where $(P_{V_r}u)(y^k) = P_{V_r} u(y^k)$.  
 The restriction of a function $w \in L^p_\mu(\Xi;V)$ to the subset $\Xi_K$, which is the tuple $\{w(y^k)\}_{k=1}^K \in V^K $, can be identified   with a tensor 
 $\mathbf{w}$ in the tensor space $\Rbb^K \otimes V$ equipped with a norm $\Vert \cdot\Vert_p$ such that $\Vert \mathbf{w} \Vert_p = \Vert w \Vert_{p,K}$. 
  The restriction to $\Xi_K$ of the best rank-$r$ approximation $u_r$ of $u$
  is then identified with the 
 best rank-$r$ approximation $\mathbf{u}_r$ of $\mathbf{u}$ in $\Rbb^K\otimes V$, and can be written
 $$
\mathbf{u}_r = \sum_{i=1}^r \mathbf{s}_i \otimes v_i \in \Rbb^K \otimes V,
$$
where $\mathbf{s}_i \in \Rbb^K$ can be identified with sample evaluations $\{s_i(y^k)\}_{k=1}^K$ of a certain function $s_i\in L^p_\mu(\Xi)$ such that 
\begin{align}
u_r(y^k) =  \sum_{i=1}^r s_i(y^k) v_i.\label{eq:rank-u-approx-samples}
\end{align}
Any rank-$r$ function $u_r$ whose restriction to $\Xi_K$ is identified with $\mathbf{u}_r$ is a solution of the best approximation 
problem \eqref{bestrankrseminorm}. The selection of a particular solution $u_r$ requires an additional approximation step. 
Such a particular solution can be obtained by interpolation of functions $s_i$ on the set of points $\Xi_K$ (e.g. using polynomial interpolation on structured grids $\Xi_K$, or nearest neighbors, Shepard or Radial Basis interpolations for unstructured samples). 
\begin{rem}
Other approximation methods (e.g. least-squares) can be used for 
the approximation of functions $s_i$ from their evaluations at sample points $\Xi_K$. However, if the interpolation property is not satisfied, then the resulting function $u_r$ is not necessarily a solution of \eqref{bestrankrseminorm}.
\end{rem}

\paragraph{Case $p=2$.}
For $p=2$ and $V$ a Hilbert space, the norm $\Vert \cdot\Vert_{2}$ on $\Rbb^K\otimes V$ such that 
$\Vert \mathbf{w}\Vert_{2} = \Vert w \Vert_{2,K}$ coincides 
with the canonical inner product norm (when $\Rbb^K$ is equipped with the weighted $2$-norm $\Vert a\Vert_{2} = (\sum_{k=1}^K \omega^k \vert a_k \vert^2)^{1/2}$). Therefore, $\mathbf{u}_r$ coincides with the truncated rank-$r$ singular value decomposition of $\mathbf{u}$, where vectors $\{v_i\}_{i=1}^r$ are the $r$ dominant eigenvectors of the operator $C_u^K :v \mapsto  \sum_{k=1}^K \omega^k u(y^k) \langle u(y^k) , v \rangle_V.$ 
 The best rank-$r$ approximation error is such that 
$
\rho_r^{(2,K)}(u) = (\sum_{i=r+1}^K \sigma_i)^{1/2} ,
$
 where $\{\sigma_i\}_{i=1}^K$ is the set of singular values of $\mathbf{u}$ (eigenvalues of $C_u^K$) sorted in decreasing order. In a probabilistic setting, when $\{y^k\}_{k=1}^K$ are i.i.d. samples (drawn according to probability measure $\mu$) and $\omega^k=K^{-1}$ for all $k$, 
$C_u^K$ is the so-called empirical correlation operator of the $V$-valued random variable $u$. Its $r$-dimensional dominant eigenspace $V_r = \mathrm{span}\{v_i\}_{i=1}^r$ is a statistical estimate of the optimal $r$-dimensional subspace associated with the best rank-$r$ approximation of $u$ in $L^2_\mu(\Xi;V)$. This corresponds to the standard \emph{Principal Component Analysis}. The obtained reduced approximation space $V_r$ can then be used for computing an approximation of $u(\xi)$ in $V_r$ for all $\xi \in \Xi$. This approach is at the basis of Galerkin \emph{Proper Orthogonal Decomposition} methods for parameter-dependent equations (see e.g. \cite{kahlbacher2006galerkin}). 

\paragraph{Case $p=\infty$.}  
For $p=\infty$, the best rank-$r$ approximation is well defined and the corresponding error is 
$$
\rho^{(\infty,K)}_r(u) = \min_{V_r\in \Gbb_r(V)} \sup_{1 \le k \le K} \Vert u(y^k) - P_{V_r} u(y^k)\Vert_V = d_r(u(\Xi_K))_V,
$$
where $ d_r(u(\Xi_K))_V$ is the Kolmogorov $r$-width of the finite subset $u(\Xi_K)=\{u(y^k)\}_{k=1}^K$ of $V$. 
Suboptimal constructions of low-rank approximations can be proposed. 
In particular, on can rely on the greedy algorithm \eqref{greedy_linfini} with $\Xi$  replaced by $\Xi_K$, which results in a sequence of nested spaces $V_r$.
This algorithm coincides with the Empirical Interpolation Method (EIM) for finite parameter sets (see \cite{MAD09,Chaturantabut:2010fk}\citechap{ and \refchapEIM{}}), sometimes called Discrete Empirical Interpolation Method (DEIM). Here also, the reduced approximation space
$V_r$ can then be used for the computation of an approximation of $u(\xi)$ in $V_r$ for all $\xi \in \Xi$.

\subsection{Higher-order low-rank approximation of multivariate functions}\label{sec:higher-order-low-rank-sample}

Here, we consider the approximation of a real-valued multivariate function $g: \Xi\rightarrow \Rbb $
 from a set of evaluations $\{g(y^k)\}_{k=1}^K$ of $g$ on a set of points $ \Xi_K = \{y^k\}_{k=1}^K $ in $\Xi$. The function $g$ can be a variable of interest that is a function of a solution $u: \Xi\rightarrow V$ of a parameter-dependent equation  (i.e. $g(\xi) = Q(u(\xi);\xi)$ with 
$Q(\cdot;\xi):V\rightarrow \Rbb$). It can also be the coefficient of the approximation of a function $u:\Xi\rightarrow V$
on a certain basis of a subspace of $V$ (e.g. one of the functions $s_i$ in the representation \eqref{eq:rank-u-approx-samples}). 

Let us assume that $\mu = \mu_1\otimes \hdots \otimes \mu_d$ is a product measure on $\Xi=\Xi_1\times \hdots \times \Xi_d$, with  $\mu_\nu$ being a measure on $\Xi_\nu \subset \Rbb$, $1\le \nu\le d$.
Using the notations of Section \ref{sec:approx-fi-dim}, we consider the approximation of $g$ in a finite-dimensional subspace $X_I = X_{1,I_1}\otimes \hdots \otimes X_{d,I_d}$ in $X = L^2_{\mu_1}(\Xi_1)\otimes \hdots \otimes L^2_{\mu_d}(\Xi_d)$, where $X_{\nu,I_\nu}$ is a $K_\nu$-dimensional subspace of $L^2_{\mu_\nu}(\Xi_\nu)$ with basis $\Psi^{(\nu)} = \{\psi_{k_\nu}^{(\nu)}\}_{k_\nu\in I_\nu}$.  

\subsubsection{Least-squares}
The standard discrete least-squares method for the approximation 
of $g$ in a subset $\Sc_r$ of low-rank tensors in $X_I$
  (see e.g. \cite{BEL11,Doostan:2013,2013arXiv1305.0030C}) consists in solving
  $$
\min_{h\in \Sc_r} \Vert g-h\Vert_{2,K}^2, \quad  \text{with} \quad 
 \Vert g-h\Vert_{2,K}^2  = {\frac{1}{K}}
 \sum_{k=1}^K (g(y^k)-h(y^k))^2,
$$
which is a quadratic convex optimization problem on a nonlinear set. 
Algorithms presented in Section \ref{sec:optimization-low-rank-subsets} can be used for the solution of this optimization problem. Assuming that $\Sc_r$ admits a simple parametrization of the form 
$\Sc_r= \{ v = F_{\Sc_r}(p_1,\hdots,p_M) : p_i \in \Rbb^{N_i}, 1\le i\le M\},$
where  $F_{\Sc_r} : \bigtimes_{i=1}^M \Rbb^{N_i} \rightarrow X_I$ is a multilinear map, the 
discrete least-squares minimization problem then takes the form 
$$
 \min_{p_1\in \Rbb^{N_1}, \hdots,p_M\in \Rbb^{N_M}}  \Vert g - F_{\Sc_r}(p_1,\hdots,p_M) \Vert_{2,K}^2,
$$
where the function to minimize is quadratic and convex with respect to each argument $p_i$, $1\le i\le M$. Also, greedy algorithms presented in Section \ref{sec:greedy} can be used for the construction of approximations in low-rank formats (see \cite{2013arXiv1305.0030C} for the construction in canonical format).

When the number of available samples is not sufficient to get a stable estimate of the $\sum_{i=1}^M N_i$ real parameters, regularization techniques can be used in a quite straightforward way (see e.g. \cite{Doostan:2013} for the use of  $\ell_2$ regularization, or \cite{2013arXiv1305.0030C} for the use of sparsity-inducing regularizations). 
However, these approaches are still heuristic and for a given low-rank format, some challenging questions remain open: how much samples are required to get 
a stable approximation in this format ?  are there sampling strategies (not random)  that are optimal with respect to this format ?  

\subsubsection{Interpolation}\label{sec:interpolation}
Here we present interpolation methods for the approximation of $g$ in $ X_I$.
 
 If $\Psi^{(\nu)}$ is a set of interpolation functions associated with a set of points  $\Xi_{\nu,K_\nu} = \{y_\nu^{k_\nu}\}_{k_\nu\in I_{\nu}}$ in $\Xi_\nu$,  then $\{\psi_k(y) = \psi_{k_1}^1(y_1)\hdots \psi_{k_d}^d(y_d)\}_{k\in I}$  is a set of interpolation functions associated with 
the tensorized grid $\Xi_K = \Xi_{1,K_1}\times \hdots \times \Xi_{d,K_d}$ composed of $K=\prod_{\nu=1}^d K_\nu$ points. An interpolation $\Ic_K(u)$ of $u$ is then given by 
 $$
\Ic_K(u)(y) = \sum_{k\in I} u(y^k) \psi_k(y),
 $$
 so that $\Ic_K(u)$ is completely characterized by the order-$d$ tensor $a \in \Rbb^{K_1}\otimes \hdots \otimes \Rbb^{K_d}$ whose components $a_{k_1,\hdots,k_d}  = u(y_{1}^{k_1},\hdots,y_d^{k_d})$ are the evaluations of $u$ on the interpolation grid $\Xi_K$.

Then, low-rank approximation methods can be used in order to obtain an approximation of the tensor $a\in \Rbb^{K_1}\otimes \hdots \otimes \Rbb^{K_d}$ using only a few entries of the tensor (i.e. a few evaluations of the function $u$). 
This is related to the problem of tensor completion. A possible approach consists in evaluating some entries of the tensor taken at random and then in reconstructing the tensor by the minimization of a least-squares functional (which is an algebraic version of the least-squares method described in the previous section) or dual approaches using regularizations of rank minimization problems (see \cite{2014arXiv1404.3905R}). 
An algorithm has been introduced in \cite{Espig:2009yf} for the approximation in canonical format, using least-squares minimization with a structured set of entries selected adaptively.  
Algorithms have also been proposed for an adaptive construction of low-rank approximations of $a$ in  Tensor Train format \cite{OSE10} or  Hierarchical Tucker format \cite{Ballani2013639}.   These algorithms are extensions of Adaptive Cross Approximation (ACA) algorithm  to high-order tensors and provide approximations that interpolate the tensor $a$ at some adaptively chosen entries.

\section{Tensor-structured parameter-dependent or stochastic equations}\label{sec:parametric_problems}

In this section, we  consider a general class of linear parameter-dependent or stochastic equations and 
we formulate these equations as tensor-structured 
 equations.

\subsection{A class of linear parameter-dependent equations}\label{sec:weakly_coercive_problems}

Let $\xi$ denote some parameters taking values in a set $\Xi \subset \Rbb^s$. 
$\Xi$ is equipped with a finite measure $\mu$ (when $\xi$ are random parameters, $\mu$ is the probability measure induced by $\xi$). 
For an integrable function $g:\Xi\rightarrow \Rbb$, we denote by 
$
\int_\Xi g(y) \mu(dy)$ the integral with respect to the measure $\mu$, which is the mathematical expectation $\Ebb_\mu(g(\xi))$ for $\mu$ being a probability measure. 
Let $V$ and $W$ be Hilbert spaces and let $V'$ and $W'$ be their respective continuous dual spaces. We denote by 
$\langle \cdot,\cdot\rangle$  the duality pairing.
We consider the problem of finding $u: \Xi \rightarrow V$ such that it holds
\begin{align}
b(u(\xi),w;\xi) =    \langle f(\xi) ,w\rangle   , \quad \forall w\in W, \label{eq:problem}
\end{align}
 for almost all $\xi \in \Xi$,
where $b(\cdot,\cdot;\xi) : V\times W \rightarrow \Rbb$ is a parameter-dependent continuous bilinear form and $f(\xi)\in W'$ 
is a parameter-dependent continuous linear form. We suppose that $\xi \mapsto b(\cdot,\cdot;\xi)$ is Bochner measurable and we 
 suppose that $b(\cdot,\cdot;\xi)$ is uniformly continuous and uniformly weakly coercive, that means that there exist constants $\alpha$ and $\beta$ independent of $\xi$ such that it holds (for almost all $\xi \in \Xi$) 
\begin{align}
&\sup_{v\in V}\sup_{w\in W} \frac{b(v,w;\xi)}{\Vert v\Vert_V \Vert w\Vert_W} \le \beta<\infty, \label{eq:continuity}
\\
& \adjustlimits \inf_{v\in V} \sup_{w\in W} \frac{b(v,w;\xi)}{\Vert v\Vert_V \Vert w\Vert_W} \ge \alpha >0. \label{eq:BNB1}
\end{align} 
Also, we assume that for all $w\neq 0\in W$, we have 
\begin{align}
\sup_{v\in V} b(v,w;\xi) >0. \label{eq:BNB2}
\end{align} 
Note that condition \eqref{eq:BNB2} is deduced from \eqref{eq:BNB1}
when $\dim(V)=\dim(W)<\infty$.
When $V=W$, a parametrized family of bilinear forms $b(\cdot,\cdot;\xi) : V\times V\rightarrow \Rbb$, $\xi\in\Xi$, is uniformly coercive if there exists a constant independent of $\xi$ such that 
\begin{align}
\inf_{v\in V} \frac{b(v,v;\xi)}{\Vert v\Vert_V^2} \ge \alpha >0 \label{eq:uniformly_coercive},
\end{align}
which implies both conditions \eqref{eq:BNB1} and \eqref{eq:BNB2}. 

Let $B(\xi):V\rightarrow W'$ denote the parameter-dependent linear operator such that $\langle B(\xi) v,w\rangle = b(v,w;\xi)$ for all $(v,w)\in V\times W$. Problem \eqref{eq:problem} is therefore equivalent to the operator equation
\begin{align}
B(\xi) u(\xi) = f(\xi), \label{eq:Bu=f}
\end{align}
where 
Assumptions \eqref{eq:continuity}, \eqref{eq:BNB1} and \eqref{eq:BNB2}  are necessary and sufficient conditions for  $B(\xi)$ 
to be an isomorphism from $V$ to $W'$ which satisfies 
\begin{align}
\alpha \Vert v \Vert_V \le \Vert B(\xi) v\Vert_{W'} \le \beta \Vert v\Vert_V, \quad \forall v\in V, \label{eq:norm-isomorphism}
\end{align} 
for almost all $\xi\in \Xi$. 
Let  $B(\xi)^*:W\rightarrow V'$ denote the adjoint of $B(\xi)$, defined by $\langle B(\xi)v,w\rangle = \langle v , B(\xi)^*w\rangle$. Property \eqref{eq:BNB2} is equivalent to $\Vert B(\xi)^* w\Vert_{V'} >0$ for all $w\neq 0$. Problem \eqref{eq:problem}  admits a unique solution $u(\xi)$ satisfying 
\begin{align}
\Vert u(\xi)\Vert_V \le \frac{1}{\alpha} \Vert f(\xi)\Vert_{W'}. \label{eq:data_to_solution_map_continuous}
\end{align}
From \eqref{eq:data_to_solution_map_continuous}, it can be deduced that if $f\in L^p_\mu(\Xi;W')$
 for a  certain $p>0$, then the solution $u \in L_\mu^{p'}(\Xi;V)$ for any $p'\le p$. 
 
 \begin{rem}
Note that the above presentation includes the case of parameter-dependent algebraic equations, for which $V = W = \Rbb^N$, $B(\xi)$ is a matrix in $\Rbb^{N\times N}$ and $f(\xi)$ is a vector in $\Rbb^N$.
 \end{rem}

\subsubsection{Example 1: elliptic diffusion equation with random coefficients.}  \label{sec:example1}
Let $D $ be an open bounded domain of $\Rbb^m$ with Lipschitz boundary $\partial \Omega$. Let $\kappa$ be a random field indexed by $x\in D$ defined on a probability space $(\Omega,\Bc,\Pbb)$ and such that it can be expressed as a function 
of random variables $\xi:\Omega\rightarrow \Xi\subset \Rbb^s$, that means $\kappa= \kappa(x,\xi)$. We  consider the following boundary value problem:
\begin{align*}
-\nabla \cdot (\kappa(\cdot,\xi) \nabla u(\xi)) = g(\cdot,\xi) \quad \text{on} \quad  D, \quad u = 0 \quad \text{on} \quad \partial D,
\end{align*}
with $g(\cdot,\xi) \in L^2(D)$. Let $V$ be an approximation space in $H^1_0(D)$, $W=V$, and 
$\Vert v\Vert_V =\Vert v\Vert_W=\left(\int_{D} \vert \nabla v \vert^2\right)^{1/2}$. A Galerkin approximation of the solution, still denoted $u(\xi)\in V$, is the solution of  \eqref{eq:problem} where $b(\cdot,\cdot;\xi) : V \times V \rightarrow \Rbb$ and $f(\xi)\in V'$ are bilinear and linear forms defined by
\begin{align*}
&b(v,w;\xi) = \int_D \kappa(\cdot,\xi) \nabla v \cdot  \nabla w,\quad \text{and} \quad  \langle f(\xi) , w\rangle = \int_{D} g(\cdot,\xi) w.
 \end{align*}
 If $\kappa$ satisfies almost surely and almost everywhere 
\begin{align}
\alpha \le \kappa(x,\xi) \le \beta,\label{eq:kappabounds}
\end{align}
then properties \eqref{eq:continuity} and \eqref{eq:uniformly_coercive} are satisfied.
Let us consider a classical situation where  $\kappa$ admits the following representation
\begin{align}
\kappa(x,\xi) = \kappa_0(x) +  \sum_{i=1}^{N} \kappa_i(x) \lambda_i(\xi), \label{eq:kappa_decomp}
\end{align}
yielding the following decomposition of the parameter-dependent bilinear form $b$:
$$
b(v,w;\xi) = \int_D \kappa_0 \nabla v \cdot\nabla w + \sum_{i=1}^N \left(\int_D \kappa_i \nabla v \cdot\nabla w \right) \lambda_i(\xi).
$$
For a spatially correlated second order random field $\kappa$, the representation \eqref{eq:kappa_decomp} can be obtained by using truncated Karhunen-Lo\`eve decomposition and truncated polynomial chaos expansions (see e.g. \cite{nouy2014random}). This representation also holds in the case where $\kappa_0$ is a mean diffusion field and the $\lambda_i$ represent random fluctuations of the diffusion coefficient in subdomains $D_i \subset D$ characterized by their indicator functions $\kappa_i(x) = I_{D_i}(x)$. This problem has been extensively analyzed, see e.g. \cite{BAB05,FRA05,MAT05}. 
\begin{rem}
For some problems of interest, the random field $\kappa$ may not be uniformly bounded (e.g. when considering log-normal random fields)  and may only satisfy 
 $0<\alpha(\xi) \le \kappa(x,\xi) \le \beta(\xi)<+\infty$ where $\alpha$ and $\beta$ possibly depend on $\xi$. For the mathematical analysis of such stochastic problems, we refer the reader to  \cite{mugler2013convergence,Charrier:2012,Charrier:2013,Charrier:2013fk,nouy2014random}.  
\end{rem} 
  
\subsubsection{{Example 2: evolution equation}}\label{sec:example2}
Let $D$ denote a bounded domain of $\Rbb^m$ with Lipschitz boundary $\partial \Omega$ and let $I=(0,T)$ denote a time interval. We consider the following evolution equation\begin{align*}
&\frac{\partial u}{\partial t} - \nabla \cdot(\kappa(\cdot,\xi)  \nabla u) = g(\cdot,\cdot,\xi)\quad \text{on}\; D\times I,
\end{align*}
with initial and boundary conditions,
\begin{align*}
u = u_0(\cdot,\xi) \quad \text{on}\; D\times \{0\}, \quad \text{and} \quad 
u = 0\quad \text{on} \; \partial D \times I.&
\end{align*}
We assume that $\kappa$ satisfies the same properties as in Example 1,
$g(\cdot,\cdot,\xi)\in L^2(D\times I)$,  and $u_0(\cdot,\xi)\in L^2(D)$.
A space-time Galerkin approximation of the solution, still denoted $u(\xi)$, can be defined by introducing an approximation space 
$$V \subset  L^2(I;H^1_0(D)) \cap H^1(I;L^2(D)):=\Vc, $$
equipped with the norm $\Vert\cdot\Vert_V$ such that $\Vert v \Vert_V^2 =  \Vert v \Vert^2_{ L^2(I;H^1_0(D))} + \Vert v \Vert^2_{ H^1(I;L^2(\Omega))}$, and a test space
 $$W = W_1 \times W_2 \subset   L^2(I;H^1_0(D)) \times L^2(D) := \Wc,$$ equipped with the norm $\Vert \cdot\Vert_W $ such that for $w=(w_1,w_2)\in W$, $\Vert w\Vert_W^2 = \Vert w_1\Vert_{L^2(I;H^1_0(D))}^2 + \Vert w_2\Vert_{L^2(D)}^2$. Then, the Galerkin approximation $u(\xi)\in V$ is defined by equation \eqref{eq:problem} where  the parameter-dependent bilinear form $b(\cdot,\cdot;\xi):V\times W \rightarrow \Rbb$ and the parameter-dependent linear form $f(\xi) : W \rightarrow \Rbb$ are defined for $v\in V$ and $w=(w_1,w_2) \in W$ by
\begin{align*}
&b(v,w;\xi) = \int_{D\times I} \frac{\partial v}{\partial t} w_1 + \int_{D\times I} \kappa(\cdot,\xi) \nabla v \cdot \nabla w_1   + \int_{D} v(\cdot,0) w_2 , \quad \text{and}
\\
&\langle f(\xi) , w \rangle = \int_{D\times I} g(\cdot,\cdot,\xi) w_1 + \int_D u_0(\cdot,\xi) w_2.
\end{align*}
For the analysis of this formulation, see \cite{Schwab-Stevenson:2009}.
 
\begin{rem}\label{rem:tensor-structure-space-time-model-example2}
$L^2(I;H^1_0(D))$ and $H^1(I;L^2(D))$ are identified with tensor Hilbert spaces 
$\overline{L^2(I) \otimes H^{1}_0(D)}^{\Vert\cdot\Vert_{L^2(I;H^1_0(D))}}$ and $  \overline{H^1(I) \otimes L^2(D)}^{\Vert\cdot\Vert_{H^1(I;L^2(D))}}$ respectively, so that 
the space $\Vc = L^2(I;H^1_0(D)) \cap H^1(I;L^2(D))$ is an intersection tensor Hilbert space which coincides with $ \overline{H^1(I) \otimes H_0^1(D)}^{\Vert\cdot\Vert_V}$ (see \cite[Section 4.3.6]{HAC12}).  Approximation spaces $V$ in $\Vc$ can be chosen of the form $V= V(I)\otimes V(D)$ in the algebraic tensor space $H^1(I) \otimes H_0^1(D)$, with approximation spaces (e.g. finite element spaces) $ V(I) \subset H^1(I)$ and $V(D)\subset H_0^1(D)$. Low-rank methods can also exploit this tensor structure and provide approximations of elements $v\in V$ under the form $v(x,t) = \sum_{i=1}^r  a_i(x) b_i(t)$, with $a_i\in V(D)$ and $b_i\in V(I)$. This is the basis of POD methods for evolution problems\citechap{ (see \refchapPOD{})} and of the first versions of Proper Generalized Decomposition methods which were introduced for solving evolution equations with variational formulations in time  \cite{LAD99b,NOU04,LAD10,NOU10b}.
\end{rem}
 
\subsection{Tensor-structured equations}\label{sec:Tensor structured equations}

Let us assume that $f\in L^2_\mu(\Xi;W'),$ so that the solution of \eqref{eq:problem} is in $L^2_\mu(\Xi;V)$. 
In this section, we use the notations $\boldsymbol{V} = L^2_\mu(\Xi;V)$ and $ \boldsymbol{W}=L^2_\mu(\Xi;W)$.
The solution $u \in\boldsymbol{V} $ satisfies 
\begin{align}
a(u,w) = F(w), \quad \forall w \in \boldsymbol{W},\label{eq:varprob_L2}
\end{align}
where  $a : \boldsymbol{V}\times \boldsymbol{W} \rightarrow \Rbb$ is the bilinear form defined by
\begin{align}
a(v,w) =  
\int_{\Xi} b(v(y),w(y);y) \mu(dy), \label{eq:bilinear_form_a}
\end{align}
and $F: \boldsymbol{W} \rightarrow \Rbb$ is the continuous linear form defined by 
$$
F(w) = 
 \int_{\Xi} \langle f(y),w(y) \rangle \mu(dy).
$$
Under Assumptions \eqref{eq:continuity}, \eqref{eq:BNB1} and \eqref{eq:BNB2}, it can be proved that $a$ satisfies 
 \begin{align}
&\sup_{v\in \boldsymbol{V}}\sup_{w\in \boldsymbol{W}} \frac{a(v,w)}{\Vert v\Vert_{\boldsymbol{V}} \Vert w\Vert_{\boldsymbol{W}}} \le \beta<\infty, \label{eq:weak_continuity}
\\
& \adjustlimits \inf_{v\in \boldsymbol{V}} \sup_{w\in \boldsymbol{W}} \frac{a(v,w)}{\Vert v\Vert_{\boldsymbol{V}} \Vert w\Vert_{\boldsymbol{W}}} \ge \alpha >0, \label{eq:weak_BNB1}
\end{align} 
and for all $0\neq w\in \boldsymbol{W}$,
\begin{align}
\sup_{v\in \boldsymbol{V}} a(v,w) >0. \label{eq:weak_BNB2}
\end{align} 
Equation \eqref{eq:varprob_L2} can be equivalently rewritten as an operator equation 
\begin{align}
A u =F, \label{eq:stochastic_operator_equation_L2}
\end{align}
where $A : \boldsymbol{V} \rightarrow \boldsymbol{W}'$ is the continuous linear operator associated with $a$, such that \begin{align}
\langle A v , w\rangle = a(v,w) \quad \text{for all } (v,w)\in \boldsymbol{V}\times \boldsymbol{W}.\label{def_operator_A}
\end{align} 
Properties \eqref{eq:weak_continuity}, \eqref{eq:weak_BNB1} and \eqref{eq:weak_BNB2}  imply that $A$ is an isomorphism from $\boldsymbol{V}$ to $\boldsymbol{W}'$ such that for all $v\in \boldsymbol{V}$,
\begin{align}
\alpha \Vert v\Vert_{\boldsymbol{V}} \le \Vert A v\Vert_{\boldsymbol{W}'} \le \beta \Vert v \Vert_{\boldsymbol{V}}. \label{eq:norm_iso_A}
\end{align}
Problem \eqref{eq:varprob_L2} therefore admits a unique solution such that $\Vert u \Vert_{\boldsymbol{V}} \le \frac{1}{\alpha} \Vert F\Vert_{\boldsymbol{W}'}$.

\subsubsection{Order-two tensor structure}\label{sec:Order-two-tensor-structure}
$u$ (resp. $f$), as an element of the Bochner space $L^2_\mu(\Xi;V)$ (resp. $L^2_\mu(\Xi;W')$), can be identified with a tensor in 
$L^2_\mu(\Xi) \otimes_{\Vert\cdot\Vert_2} V$ (resp. $ L^2_\mu(\Xi) \otimes_{\Vert\cdot\Vert_2} W'$).
Let us further assume that $f \in L^2_\mu(\Xi) \otimes W'$ admits the following representation 
\begin{align}
f(\xi)= \sum_{i=1}^L \gamma_i(\xi) f_i ,\label{eq:decomp_f}
\end{align}
 with $f_i\in W'$ and $\gamma_i\in L^2_\mu(\Xi)$. Then $F$ is identified with the finite-rank tensor 
 \begin{align}
 F=\sum_{i=1}^L \gamma_i \otimes f_i  .
 \label{order-two-tensor-structure-F}
 \end{align}
 Let us now assume that the parameter-dependent operator $B(\xi):V\rightarrow W'$ associated with the parameter-dependent bilinear form $b(\cdot,\cdot;\xi)$ admits the following representation (so called \emph{affine representation} in the context of Reduced Basis methods)
 \begin{align}
B(\xi) = \sum_{i=1}^{R}  \lambda_i(\xi)B_i,\label{eq:decomp_B}
\end{align}
where the $B_i:V\rightarrow W'$ are parameter-independent operators associated with parameter-independent bilinear forms 
$b_i$,  and where the $\lambda_i$ are real-valued functions defined on $\Xi$.
\begin{rem}
Let us assume that $\lambda_i\in L^{\infty}_\mu(\Xi)$, $1\le i\le R$, and $\lambda_1\ge 1$.
Let us denote by $\alpha_i$ and $\beta_i$ the constants such that 
 $$\alpha_i \Vert  v\Vert_V \le \Vert B_i v \Vert_{W'} 
 \le \beta_i \Vert v\Vert_V.$$
Property \eqref{eq:norm-isomorphism} is satisfied with $\beta =  \sum_{i=1}^R \beta_i\Vert \lambda_i\Vert_{\infty}$ and with 
$\alpha = \alpha_1 - \sum_{i=2}^R \beta_i \Vert \lambda_i\Vert_{\infty}$ if $\alpha_1> \sum_{i=2}^R \beta_i \Vert \lambda_i\Vert_{\infty}$.
In the case where  $V=W$, if all the $B_i$ satisfy $\inf_{v\in V} \frac{\langle B_i v,v\rangle}{\Vert v\Vert_V^2} \ge \alpha_i > -\infty$, then Property  \eqref{eq:norm-isomorphism} is satisfied 
with either $\alpha = \alpha_1 - \sum_{i=2}^R \alpha_i \Vert \lambda_i\Vert_{\infty}$ if   $\alpha_1> \sum_{i=2}^R \alpha_i \Vert \lambda_i\Vert_{\infty}$, or $\alpha = \alpha_1$ if $\alpha_1>0$ and $\alpha_i\ge 0$ and $\lambda_i\ge 0$ for all $i\ge 2$.
 \end{rem}
 
 \begin{rem}
 If the parameter-dependent operator $B(\xi)$ (resp. right-hand side $f(\xi)$)  does not admit an affine representation of the form \eqref{eq:decomp_B} (resp. \eqref{eq:decomp_f}), or if the initial affine representation contains
a high number of terms,  low-rank approximation methods can be used in order to obtain an affine representation with a small number of terms. For that purpose, one can rely on SVD or on the Empirical Interpolation Method\citechap{ (see \refchapEIM{})}, the latter approach being commonly used in the context of Reduced Basis Methods.
 \end{rem}
Assuming that $\lambda_i\in L^{\infty}_\mu(\Xi)$,  
the operator $A : \boldsymbol{V}  
\rightarrow \boldsymbol{W}' 
$ admits the following representation\footnote{$A$ is a finite-rank tensor in $\Lc(L^2_\mu(\Xi),L^2_\mu(\Xi))\otimes \Lc(V,W')$.}
\begin{align}
A = \sum_{i=1}^{R} \Lambda_i \otimes B_i ,\label{order-two-tensor-structure-A}
\end{align}
where $\Lambda_i : L^2_\mu(\Xi) \to L^2_\mu(\Xi)$ is a continuous linear operator associated with $\lambda_i$ such that for $\psi \in L^2_\mu(\Xi)$, $\Lambda_i \psi $ is defined by
$$\langle \Lambda_i \psi , \phi \rangle = 
\int_{\Xi} \lambda_i(y) \psi(y)\phi(y) \mu(dy)
\quad \text{ for all }  \phi \in 
 L^2_\mu(\Xi).$$
  Therefore, Equation \eqref{eq:stochastic_operator_equation_L2} can be written as a tensor-structured equation
  \begin{align}
\Big(\sum_{i=1}^{R} \Lambda_i \otimes B_i\Big) u =\sum_{i=1}^L \gamma_i \otimes f_i  . \label{eq:stochastic_tensor_format}
\end{align}

\subsubsection{Higher-order tensor structure}\label{sec:higher-order-tensor-structure}
Let us assume that $\mu = \mu_1\otimes \hdots \otimes \mu_d$ is a product measure on $\Xi=\Xi_1\times \hdots \times \Xi_d$, with  $\mu_\nu$ being a measure on $\Xi_\nu \subset \Rbb^{s_\nu}$, $1\le \nu\le d$, with  $s=\sum_{\nu=1}^d s_\nu$. Then $L^2_\mu(\Xi)= {}_{\Vert\cdot\Vert_2} \bigotimes_{\nu=1}^d L^2_{\mu_\nu}(\Xi_\nu)$ (see Section \ref{sec:Lp_spaces}).  
In a probabilistic context,  $\mu$ would be the measure induced by $\xi = (\xi_1,\hdots,\xi_d)$, where the $\xi_\nu$ are independent random variables with values in $\Xi_\nu$ and probability law $\mu_\nu$.

Let us assume that the functions $\gamma_i$ in \eqref{eq:decomp_f}, $1\le i\le L$,  are such that 
\begin{align}
\gamma_i(\xi) =  \gamma_i^{(1)}(\xi_1) \hdots  \gamma_i^{(d)}(\xi_d),\label{ass:gammai}
\end{align}
with $\gamma^{(\nu)}_i \in L^2_{\mu_\nu}(\Xi_\nu) $. Then $f$ is an element of 
 $L^2_{\mu_1}(\Xi_1)\otimes \hdots \otimes  L^2_{\mu_d}(\Xi_d) \otimes W' $ and $F$ admits the following representation
\begin{align}
F = \sum_{i=1}^L  \gamma_i^{(1)} \otimes \hdots \otimes  \gamma_{i}^{(d)}
 \otimes f_i . \label{higher-order-tensor-structure-F}
 \end{align} 
Let us assume that  in the representation \eqref{eq:decomp_B} of $B(\xi)$, the  functions $\lambda_i$, $1\le i \le R$, are such that 
\begin{align}
\lambda_i(\xi) =  \lambda^{(1)}_i(\xi_1)\hdots \lambda^{(d)}_i(\xi_d).\label{ass:lambdai}
\end{align}
Assuming that $\lambda^{(\nu)}_i \in L^\infty_{\mu_\nu}(\Xi_\nu)$,  $\lambda_i^{(\nu)}$ can be identified with an operator $\Lambda_i^{(\nu)} : S_\nu
 \rightarrow \tilde S_\nu'$, where for $\psi\in S_\nu$, $\Lambda^{(\nu)}_i \psi$ is defined by 
 $$\langle \Lambda^{(\nu)}_i \psi , \phi\rangle = \int_{\Xi_\nu}  \lambda^{(\nu)}_i(y_\nu) \psi(y_\nu)\phi(y_\nu)\mu_\nu(dy_\nu)
 \quad \text{ for all } \phi\in \tilde S_\nu.$$ Then, $\lambda_i$ also defines an operator $\Lambda_i : S \rightarrow \tilde S'$ such that 
$$
\Lambda_i=   \Lambda^{(1)}_i \otimes \hdots \otimes \Lambda^{(d)}_i.
$$
Then the operator $A$, as an operator from $L^2_{\mu_1}(\Xi_1)\otimes \hdots \otimes L^2_{\mu_d}(\Xi_d) \otimes V$ to  
$(L^2_{\mu_1}(\Xi_1)\otimes \hdots \otimes L^2_{\mu_d}(\Xi_d) \otimes W)'$,
admits the following  decomposition\footnote{$A$ is a finite-rank tensor in $\Lc(L^2_{\mu_1}(\Xi_1),L^2_{\mu_1}(\Xi_1))\otimes \hdots \otimes \Lc(L^2_{\mu_d}(\Xi_d),L^2_{\mu_d}(\Xi_d))\otimes \Lc(V,W')$.}
\begin{align}
A = \sum_{i=1}^{R} \Lambda^{(1)}_i \otimes \hdots \otimes \Lambda^{(d)}_i \otimes B_i.
\label{higher-order-tensor-structure-A}
 \end{align} 
Therefore, Equation \eqref{eq:stochastic_operator_equation_L2} can be written as a tensor-structured equation
\begin{equation}
\Big(  \sum_{i=1}^{R} \Lambda^{(1)}_i \otimes \hdots \otimes \Lambda^{(d)}_i \otimes B_i \Big) u =
 \sum_{i=1}^L  \gamma_i^{(1)} \otimes \hdots \otimes  \gamma_{i}^{(d)} \otimes f_i.
\label{eq:stochastic_tensor_format_high_order}
\end{equation}

\subsection{Galerkin approximations}\label{sec:galerkin}

Here, we present Galerkin methods for the approximation of the solution of {\eqref{eq:problem}} in a subspace $S \otimes  V $ of $L^2_\mu(\Xi;V)$, where $S$ is a finite-dimensional subspace in $L^2_\mu(\Xi)$. 
In this section, $\boldsymbol{V}= S \otimes  V$ denotes the approximation space in $L^2_\mu(\Xi;V)$, which is 
equipped with the natural norm in $L^2_\mu(\Xi;V)$, denoted $\Vert\cdot\Vert_{\boldsymbol{V}}$.

 \subsubsection{Petrov-Galerkin approximation}\label{sec:petrov-galerkin}
 
Let us introduce a finite-dimensional subspace $\tilde S$ in $L^2_\mu(\Xi)$, with $\dim(S)=\dim(\tilde S)$, and let us introduce the tensor space 
$\boldsymbol{W}=\tilde S \otimes  W \subset L^2_\mu(\Xi;W)$, equipped with the natural norm in $L^2_\mu(\Xi;W)$, denoted $\Vert\cdot\Vert_{\boldsymbol{W}}$.  
A Petrov-Galerkin approximation in $\boldsymbol{V}=  S\otimes V$ of the solution of Problem  \eqref{eq:problem}, denoted $u_G$, 
is defined by the equation 
\begin{align}
a(u_G,w) = F(w), \quad \forall w \in \boldsymbol{W},\label{eq:varprob_galerkin}
\end{align}
which can be equivalently rewritten as an operator equation 
\begin{align}
A u_G =F, \label{eq:stochastic_operator_equation_galerkin}
\end{align}
where $A : \boldsymbol{V} \rightarrow \boldsymbol{W}'$ is associated (through Equation  \eqref{def_operator_A}) with the bilinear form $a$.
\begin{rem}
Assuming that the approximation space $\boldsymbol{V}$ and the test space $\boldsymbol{W}$ are such that 
Properties \eqref{eq:weak_continuity}, \eqref{eq:weak_BNB1} and \eqref{eq:weak_BNB2} are satisfied,  then 
$A$ satisfies \eqref{eq:norm_iso_A} and Equation \eqref{eq:stochastic_operator_equation_galerkin} admits a unique solution $u_G$ which is a quasi-optimal approximation of $u$, with  
$\Vert u_G-u \Vert_{\boldsymbol{V}} \le (1+\frac{\beta}{\alpha}) \min_{v \in \boldsymbol{V}} \Vert u-v \Vert_{\boldsymbol{V}}$.
\end{rem}
\begin{rem}
Letting $\{\psi_i\}_{i=1}^{K}$ and $\{\phi_i\}_{i=1}^{K}$ be bases of $S$ and $\tilde S$ respectively, the solution $u_G$ of \eqref{eq:stochastic_operator_equation_galerkin} can be written $u_G = \sum_{i=1}^K  \psi_i \otimes u_i$, where the tuple $\{u_i\}_{i=1}^K \in V^K$ verifies the coupled system of equations
\begin{align}
\sum_{j=1}^P A_{ij} u_j = F_i, \quad 1\le i\le K, \label{eq:galerkin_block_system}
\end{align}
with 
$
A_{ij} = \int_\Xi B(y)\psi_j(y)\phi_i(y)\mu(dy) 
$ and $F_i = \int_\Xi f(y)\phi_i(y) \mu(dy)
$,
for $1\le i,j\le K$. The tuple   $\{u_i\}_{i=1}^K \in V^K$  can be identified with a tensor in $\Rbb^K \otimes V$. \end{rem}
\begin{rem}\label{Galerkin-quadrature}In practice, integrals over $\Xi$ with respect to the measure $\mu$ 
can be approximated by using a suitable quadrature rule $\{(y^k,\omega^k)\}_{k=1}^K$, therefore replacing \eqref{eq:varprob_galerkin} by 
$$
\sum_{k=1}^K \omega^k \langle B(y^k)u_G(y^k),w(y^k)\rangle = \sum_{k=1}^K \omega^k \langle f(y^k),w(y^k)\rangle.
$$
 \end{rem}
Under the assumptions of Section \ref{sec:Order-two-tensor-structure}, Equation \eqref{eq:stochastic_operator_equation_galerkin} can be written in the form of a tensor-structured equation \eqref{eq:stochastic_tensor_format}, where the functions $\gamma_i$ are now identified with elements of $\tilde S'$ such that 
$\langle \gamma_i, \psi \rangle =\int_\Xi \gamma_i(y) \psi(y) \mu(dy)
$ for all $\psi\in \tilde S$, and where the $\Lambda_i$ are now considered as operators from $S$ to $\tilde S'$ such that for $\psi \in \tilde S$, $\Lambda_i \psi$ is defined by 
$\langle \Lambda_i \psi , \phi\rangle = \int_\Xi \lambda_i(y)\psi(y)\phi(y) \mu(dy)
$ for all $\phi\in \tilde S$.

Under the stronger assumptions of Section \ref{sec:higher-order-tensor-structure},  \eqref{eq:stochastic_operator_equation_galerkin} can be written in the form of a tensor-structured equation \eqref{eq:stochastic_tensor_format_high_order}, where the functions $\gamma_i^{(\nu)}$ are now identified with elements of $\tilde S_\nu'$ such that 
$\langle \gamma_i^{(\nu)}, \psi \rangle = \int_{\Xi_\nu} \gamma_i(y_\nu) \psi(y_\nu) \mu_\nu(dy_\nu)
 $ for all $\psi\in \tilde S_\nu$, and where the $\Lambda_i^{(\nu)}$ are now considered as operators from $S_\nu$ to $\tilde S_\nu'$ such that for $\psi \in \tilde S_\nu$, $\Lambda_i^{(\nu)} \psi$ is defined by 
$\langle \Lambda_i^{(\nu)} \psi , \phi\rangle = \int_{\Xi_\nu} \lambda_i^{(\nu)}(y_\nu)\psi(y_\nu)\phi(y_\nu) \mu_\nu(dy_\nu)
$ for all $\phi\in \tilde S_\nu$.
 
\subsubsection{Minimal residual Galerkin approximation}
\label{sec:parametric-galerkin-minimal-residual}
Let $C(\xi): W' \rightarrow W$ be a symmetric operator that defines on $W'$ an inner product $\langle\cdot,\cdot\rangle_{C(\xi)}$ defined by 
$\langle g , h  \rangle_{C(\xi)} = \langle g, C(\xi) h\rangle=\langle C(\xi) g,  h\rangle$ for $g,h\in W'$. Let $\Vert \cdot \Vert_{C(\xi)}$ denote the associated norm on $W'$ and assume that 
\begin{equation} \alpha_C \Vert \cdot\Vert_{W'} \le \Vert \cdot \Vert_{C(\xi)} \le \beta_C \Vert \cdot\Vert_{W'}\label{assumptionsC}
\end{equation} 
for some constants $0<\alpha_C\le \beta_C<\infty.$
\begin{rem}
A natural choice for $C$ is to take the inverse of the (parameter-independent) Riesz map $R_W:W\to W'$, so that $\Vert h \Vert_C = \Vert h \Vert_{W'} = \Vert R_W^{-1} h \Vert_W$. When $V=W$ and $B(\xi)$ is coercive, another possible choice for $C(\xi)$ is to take the inverse of the symmetric part of $B(\xi)$.
\end{rem}	
A minimal residual Galerkin approximation in $\boldsymbol{V} =  S\otimes V$ of the solution of Problem  \eqref{eq:problem}, denoted $u_R$, can be defined by 
\begin{equation}
u_R = \arg\min_{v \in \boldsymbol{V}}\; \Ec(u,v),\label{minresgalerkin_optim}
\end{equation}
with
\begin{equation}
 \Ec(u,v)^2=  
 \int_\Xi\Vert B(y)v(y) - f(y) \Vert_{C(y)}^2 \mu(dy)
 .\label{minresgalerkin_optim_functional}
\end{equation} 
 Let $B(\xi)^*:W\rightarrow V'$ denote the adjoint of $B(\xi)$.
Then, we define the symmetric bilinear form $\tilde a: \boldsymbol{V} \times \boldsymbol{V} \rightarrow \Rbb$ such that
$$
\tilde a(v,w) = 
\int_\Xi \langle B(y) v(y), B(y) w(y) \rangle_{C(y)} \mu(dy) = 
\int_\Xi \langle \tilde B(y) v(y) , w(y) \rangle\mu(dy) ,
$$
with $\tilde B(\xi) = B(\xi)^* C(\xi) B(\xi)$,
and the linear form $\tilde F: \boldsymbol{V} \rightarrow \Rbb$ such that
$$
\tilde F(w) = 
\int_\Xi \langle f(y), B(y)w(y)\rangle_{C(y)} \mu(dy)=  \int_\Xi \langle \tilde f(y), w(y)\rangle \mu(dy),
$$
with $\tilde f(\xi) = B(\xi)^* C(\xi)  f(\xi)$.
The approximation $u_R \in \boldsymbol{V} =  S\otimes V$ defined by \eqref{minresgalerkin_optim}
is equivalently defined by
\begin{align}
\tilde a(u,v) =  \tilde F(v), \quad \forall v\in \boldsymbol{V}, \label{minresgalerkin}
\end{align}
which can be rewritten as an operator equation 
\begin{align}
\tilde A u_R = \tilde F, \label{eq:stochastic_operator_equation_tilde}
\end{align}
where $\tilde A : \boldsymbol{V} \to \boldsymbol{V}'$ is the operator associated with the bilinear form $\tilde a$. The approximation $u_R$ is the standard Galerkin approximation of the solution of the parameter-dependent equation 
\begin{align}
\tilde B(\xi) u(\xi) = \tilde f(\xi).\label{eq:parametric-minres}
\end{align}
\begin{rem}
Under assumptions  \eqref{eq:continuity}, \eqref{eq:BNB1}, \eqref{eq:BNB2} and \eqref{assumptionsC}, we have that 
 \begin{align}
&\sup_{v\in \boldsymbol{V}}\sup_{w\in \boldsymbol{V}} \frac{\tilde a(v,w)}{\Vert v\Vert_{\boldsymbol{V}} \Vert w\Vert_{\boldsymbol{V}}} \le \tilde \beta <\infty,  \quad 
\inf_{v\in \boldsymbol{V}}  \frac{\tilde a(v,v)}{\Vert v\Vert_{\boldsymbol{V}}^2} \ge \tilde\alpha  >0,
\end{align} 
with $\tilde \alpha =  \alpha_C^2 \alpha^2 $ and $\tilde \beta =  \beta_C^2 \beta^2$, and $u_R$ is a quasi-optimal approximation of $u$ in $\boldsymbol{V}$, with 
$\Vert u - u_R \Vert_{\boldsymbol{V}} \le \sqrt{\frac{\tilde \beta}{\tilde \alpha}}\min_{v \in \boldsymbol{V}} \Vert u - v \Vert_{\boldsymbol{V}}$.
\end{rem}
\begin{rem}
Letting $\{\psi_i\}_{i=1}^{K}$ be a  basis of $S$, the solution $u_R$ of \eqref{eq:stochastic_operator_equation_tilde} can be written $u_R = \sum_{i=1}^K  \psi_i \otimes u_i$, where the set of vectors $\{u_i\}_{i=1}^K \in V^K$ verifies the coupled system of equations \eqref{eq:galerkin_block_system} 
with $
 A_{ij} =  
\int_\Xi \tilde B(y)\psi_j(y)\psi_i(y) \mu(dy)
 $
  and 
 $ F_i = 
 \int_\Xi \tilde f(y)\psi_i(y)\mu(dy) 
 ,$
 for $1\le i,j\le K$. The tuple   $\{u_i\}_{i=1}^K \in V^K$  can be identified with a tensor in $\Rbb^K \otimes V$.
 \end{rem}
  \begin{rem} 
  \label{rem:minresgalerkin-quadrature}
 In practice, \eqref{minresgalerkin_optim_functional} can be replaced by 
\begin{equation}
 \Ec(u,v)= \sum_{k=1}^K \omega^k \Vert B(y^k)v(y^k) - f(y^k) \Vert_{C(y^k)}^2 = 
\sum_{k=1}^K \omega^k \Vert v(y^k) - u(y^k) \Vert_{\tilde B(y^k)}^2 ,
 \label{minresgalerkin_optim_quadrature}
\end{equation} 
where  $\{(y^k,\omega^k)\}_{k=1}^K$ is a suitable quadrature rule for  the integration over $\Xi$ with respect to the measure $\mu$. 
\end{rem}
Under the assumptions of Section \ref{sec:Order-two-tensor-structure} and if we assume that 
$C$ admits an affine representation $C(\xi) = \sum_{i=1}^{R_C} C_i \eta_i(\xi)$, then $\tilde B(\xi)$ and $\tilde f(\xi)$ admit affine representations $\tilde B(\xi) = \sum_{i=1}^{\tilde R} \tilde B_i \tilde \lambda_i(\xi)$ and $\tilde f(\xi) = \sum_{i=1}^{\tilde L} \tilde f_i \tilde \gamma_i(\xi)$ respectively. Therefore, Equation \eqref{eq:stochastic_operator_equation_tilde} can be written in the form of a tensor-structured equation \eqref{eq:stochastic_tensor_format}, where all the quantities are replaced by their tilded versions. 
Under the stronger assumptions of Section \ref{sec:higher-order-tensor-structure}, if we assume that $C$ admits a representation of the form $C(\xi) = \sum_{i=1}^{R_C} C_i \eta_i^{(1)}(\xi_1)\hdots \eta_i^{(d)}(\xi_d)$, 
 then $\tilde B(\xi)$ and $\tilde f(\xi)$ admit representations of the form $\tilde B(\xi) = \sum_{i=1}^{\tilde R}\tilde B_i \tilde \lambda_i^{(1)}(\xi_1) \hdots \tilde \lambda_i^{(d)}(\xi_d)$ and $\tilde f(\xi) = \sum_{i=1}^{\tilde L} \tilde f_i \tilde \gamma_i^{(1)}(\xi_1)\hdots \tilde \gamma_i^{(d)}(\xi_d)$. Therefore, Equation \eqref{eq:stochastic_operator_equation_tilde} can be written in the form of a tensor-structured equation \eqref{eq:stochastic_tensor_format_high_order}, where all the quantities are replaced by their tilded versions. 
 
\subsection{Interpolation (or collocation) method}\label{sec:interpolation-equation}
Let $\Xi_K = \{y_k\}_{k \in I} $ be a set of $K=\#I$ interpolation points in $\Xi$ and let $\{\phi_k\}_{k\in I}$ be an associated set of interpolation functions.  
An interpolation $\Ic_K(u)$
of the solution of \eqref{eq:problem} can then be written 
$$
\Ic_K(u)(y) = \sum_{k \in I} u(y^k) \phi_k(y),
$$
where $u(y^k)\in V$ is the solution of 
\begin{align}
B(y^k) u(y^k)= f(y^k), \quad k\in I, \label{sec:set-of-equations-interpolation}
\end{align}
which can be written as an operator equation
\begin{align}
A u_I= F, \label{sec:operator-equation-interpolation}
\end{align}
where $u_I = \{u(y^k)\}_{k \in I} \in V^K$,
$F = \{f(y^k)\}_{k \in I} \in (W')^K$ and $A:V^K \to (W')^K$.
\paragraph{Order-two tensor structure.}
The tuples $u_I$ and $F$ can be identified with tensors in 
$ \Rbb^K \otimes V$ and $\Rbb^K \otimes W'$ respectively. Also, the operator ${A}$ can be identified with a tensor  in $\Rbb^{K\times K} \otimes \Lc(V,W')$.
Under the assumption \eqref{eq:decomp_B} on $B(\xi)$, $A$ can be written in the form \eqref{order-two-tensor-structure-A}, where $\Lambda_i \in \Rbb^{K \times K}$ is  the diagonal matrix $\mathrm{diag}(\lambda_i(y^1),\hdots,\lambda_i(y^K))$. Also, under the assumption \eqref{eq:decomp_f} on $f$, $F$ can be written in the form 
\eqref{order-two-tensor-structure-F}, where  ${\gamma}_i = (\gamma_i(y^1),\hdots,\gamma_i(y^K)) \in \Rbb^K$.

\paragraph{Higher-order tensor structure.}

Now, using the notations of Section \ref{sec:interpolation}, we consider 
a tensorized interpolation grid 
$\Xi_K = \Xi_{1,K_1} \times \hdots \times \Xi_{d,K_d}$ and a corresponding set of interpolation functions 
 $\{\phi_k(y) = \phi_{k_1}^1(y_1)\hdots \phi_{k_d}^d(y_d)\}_{k\in I}$, with $I=\times_{k=1}^d \{1,\hdots,K_k\}
 $. 
 The tuples $u_I$ and $F$ can now be identified with tensors in $\Rbb^{K_1} \otimes \hdots \otimes \Rbb^{K_d} \otimes V$ and $\Rbb^{K_1} \otimes \hdots \otimes \Rbb^{K_d} \otimes W'$ respectively, and the operator $A$ can be identified with a tensor in $\Rbb^{K_1\times K_1} \otimes \hdots \otimes \Rbb^{K_d\times K_d} \otimes \Lc(V,W')$.
Under the assumptions of Section \ref{sec:higher-order-tensor-structure}, $A$ can be written in the form 
\eqref{higher-order-tensor-structure-A}, with $\Lambda_i^{(\nu)} =\mathrm{diag}(\lambda_i^{(\nu)}(y_\nu^1),\hdots,\lambda_i^{(\nu)}(y_\nu^{K_\nu})) \in \Rbb^{K_\nu \times K_\nu}$, 
and $F$ can be written in the form 
\eqref{higher-order-tensor-structure-F} with $\gamma_i^{(\nu)} = (\gamma_i^{(\nu)}(y_\nu^1),\hdots,\gamma_i^{(\nu)}(y_\nu^{K_\nu})) \in \Rbb^{K_\nu}$.

\begin{rem}
Note that the interpolation (or collocation) method provides an approximation $u_I$ in $S\otimes V$, with $S=\mathrm{span}\{\phi_k\}_{k\in I}$, which coincides with the approximation obtained by a ``pseudo-spectral'' Galerkin method where the integrals over $\Xi$ are approximated using a numerical quadrature with $\{y_k\}_{k \in I}$ as the set of integration points (see Remarks \ref{Galerkin-quadrature} and \ref{rem:minresgalerkin-quadrature}).
\end{rem}

\subsection{Low-rank structures of the solution of parameter-dependent or stochastic equations}
When solving parameter-dependent or stochastic equations with low-rank tensor methods, the first question that should be asked is: does the solution $u $ present a low-rank structure or admit an accurate approximation with low rank ? Unfortunately, there are only a few quantitative answers to this question.

When $u \in L^p_\mu(\Xi;V)$ is 
seen as an order-two tensor in $\overline{L^p_\mu(\Xi)\otimes  V}^{\Vert\cdot\Vert_p}$ (see Section \ref{sec:low-rank-bochner}), there exist some results about the convergence of best rank-$r$ approximations  for some classes of functions. For $p=2$, these results are related to the decay of singular values of $u$  (or equivalently of the compact operator associated with $u$). For $p=\infty$, these results are related to the convergence of the Kolmogorov widths of the set of solutions $u(\Xi)$.  Exploiting these results requires a fine analysis of parameter-dependent  (or stochastic) equations in order to precise the class of their solutions. 
The question is more difficult when looking at $u$ as a higher-order tensor in 
$\overline{L^p_{\mu_1}(\Xi_1) \otimes \hdots \otimes L^p_{\mu_d}(\Xi_d) \otimes V}^{\Vert\cdot\Vert_p}$, in particular because of the combinatorial nature of the definition of rank.
Some results are available for the convergence of best rank-$r$ approximations of some general classes of functions, for canonical or tree-based Tucker formats \cite{Temlyakov:1989,Schneider201456,Hackbusch:2013vn}. However, these results usually exploit some global regularity and do not exploit specific structures of the solution (such as anisotropy), which would again require a detailed analysis of the parameter-dependent equations. 
\begin{ex}
Low-rank structures can be induced by particular parametrizations of operators and right-hand sides. As an example,  consider the equation 
$B(\xi_1) u(\xi_1,\xi_2) = f(\xi_2),
$ where $u$ and $f$ are considered as tensors in $L^2_{\mu_1}(\Xi_1) \otimes L^2_{\mu_2}(\Xi_2) \otimes V$ and $L^2_{\mu_1}(\Xi_1) \otimes L^2_{\mu_2}(\Xi_2) \otimes W'$ respectively. 
Here, $\rank_1(f)=1$ and $\rank_2(f) = \rank_3(f) $. Then, $\rank_2(u)=\rank_2(f)$.
\end{ex}
\begin{ex}{There are specific structures that are particular cases of low-rank structures and that can therefore be captured by low-rank methods, such as low effective dimensionality or low-order interactions. For example, a function $u(\xi_1,\hdots,\xi_d)$ which can be well approximated by a low-dimensional function $\tilde u(\xi_{\beta_1},\hdots,\xi_{\beta_k})$, with $\{\beta_1,\hdots,\beta_k\}:=\beta \subset \{1,\hdots,d\}$, can therefore be approximated with a tensor $\tilde u$ with $\rank_\beta(\tilde u)=1$ and $\rank_\alpha(\tilde u)=1$ for any $\alpha \subset \{1,\hdots,d\} \setminus \beta$. When using tree-based tensor formats, the tree should be adapted in order to reveal these low-rank structures. 
}\end{ex}

Although only a few a priori results are available, it is observed in many applications that the solutions of parameter-dependent (or stochastic) equations can be well approximated using low-rank tensor formats. 
However, there is a need for a rigorous classification of problems in terms of the expected accuracy of the low-rank approximations of their solutions.  
In the following section, we let apart the discussion about the good approximability of functions by low-rank tensors and focus on numerical methods for computing low-rank approximations. 
\begin{rem}
Note that results on the convergence of best $r$-term approximations on a polynomial basis for particular classes of  parameter-dependent equations \cite{Cohen:2011,KHO11,Chkifa2014}, which exploit the anisotropy of the  solution map $u:\Xi\to V$, provide as a by-product upper bounds for the convergence of low-rank approximations.
\end{rem}

 \section{Low-rank approximation for equations in tensor format}\label{sec:low-rank-equations}

In this section, we present algorithms for the approximation in low-rank formats of the solution of an operator equation
\begin{equation}
A u = F, \label{eq:operator_equation} 
\end{equation}
where ${u}$ is an element of a finite-dimensional tensor space $\boldsymbol{V}=V_1\otimes \hdots \otimes V_{D}$, and $A $ is an operator from $\boldsymbol{V}$ to $\boldsymbol{W}'$, with $\boldsymbol{W}=W_1\otimes \hdots \otimes W_{D}$.
We can distinguish two main families of algorithms.  
The first family of algorithms relies on the use of classical iterative methods  with efficient low-rank truncations of iterates.
The second family of algorithms directly computes a low-rank approximation of the solution based on the minimization of a certain distance between the solution $u$ and its low-rank approximation, using either a direct minimization in subsets of low-rank tensors or suboptimal but constructive (greedy) algorithms.

The algorithms are presented in a general setting and will be detailed for the particular case of
parameter-dependent equations presented in Section \ref{sec:parametric_problems}. In this particular case, 
the space $\boldsymbol{V}$ can be considered as a space $S \otimes V$ of order-two tensors ($D=2$), where $S$ is either a $K$-dimensional approximation space in $L^2_\mu(\Xi)$ (for Galerkin methods) or $\Rbb^K$ (for interpolation or collocation methods). The space $\boldsymbol{V}$ can also be considered as a space $S_1\otimes \hdots \otimes S_d \otimes V$ of higher-order tensors ($D=d+1$), where $S_\nu$ is either a $K_\nu$-dimensional subspace of $L^2_{\mu_\nu}(\Xi_\nu)$ (for Galerkin methods) or $\Rbb^{K_\nu}$ (for interpolation or collocation methods). 

In practice,  bases of finite-dimensional tensor spaces $\boldsymbol{V}$ and $\boldsymbol{W}$ are introduced,  
so that \eqref{eq:operator_equation}  can be equivalently rewritten
\begin{equation}
\mathbf{A} \mathbf{u}= \mathbf{F}, \label{eq:operator_equation_discrete} 
\end{equation}
where $\mathbf{u} \in \mathbf{X} = \Rbb^{N_1}\otimes \hdots \otimes \Rbb^{N_D}$ is the set of coefficients of $u$ in the chosen basis of $\boldsymbol{V}$, and where $\mathbf{A} \mathbf{u}$ and $\mathbf{F}$ are respectively the coefficients of $Au$ and $F$ in the dual basis of the chosen basis  in $\boldsymbol{W}$. Here 
$\mathbf{A}$ is an operator from $\mathbf{X}$ to $\mathbf{X}$.

\subsection{Classical iterative methods using low-rank truncations}\label{sec:iterative_methods}

 Simple iterative algorithms  (e.g. Richardson iterations, Gradient algorithm...)  
  take the form ${u}^{i+1} = M({u}^i)$, where $M$ is an iteration map which involves simple algebraic operations (additions, multiplications...) between tensors.  Low-rank truncation methods can be systematically used to reduce the storage and computational complexities of these algebraic operations. This results in approximate iterations 
  $${u}^{i+1} \approx M({u}^i),$$ 
  where the iterates $\{{u}^{i}\}_{i\ge 1}$ are in low-rank format. 
 The resulting algorithm can be analyzed as a perturbed version of the initial algorithm (see e.g. \cite{HAC08}). The reader is referred to \cite{Kressner:2011ys,KHO11,Ballani:2013vn} for a detailed introduction to these techniques in a general algebraic setting, and to \cite{MAT12} for an application to parameter-dependent equations. Note that these iterative methods usually require the construction of good preconditioners in low-rank tensor formats (see \cite{Touzene:2008,Khoromskij:2009fk,Oseledets:2012,Giraldi:2014}).

As an example, let us consider  simple 
Richardson iterations for solving \eqref{eq:operator_equation_discrete},
where $M(\mathbf{u}) = \mathbf{u}  - \alpha (\mathbf{A} \mathbf{u} -\mathbf{F})$. Approximate iterations can take the form 
 $$\mathbf{u}^{i+1} = \Pi_\epsilon( \mathbf{u}^i - \alpha (\mathbf{A}\mathbf{u}^i -\mathbf{F}) ) ,$$ where $\Pi_\epsilon$ is an operator which associates to a tensor $\mathbf{u}$ an approximation $ \Pi_\epsilon(\mathbf{u})$ in low-rank format with a certain precision $\epsilon$.  
\paragraph{Low-rank truncation controlled in $2$-norm.}
For a control of the error  in the $2$-norm, the operator $\Pi_\epsilon$ provides an approximation 
 $\Pi_\epsilon(\mathbf{w})$ of a tensor $\mathbf{w}$ such that  $\Vert \mathbf{w} - \Pi_\epsilon(\mathbf{w})\Vert_2 \le \epsilon \Vert \mathbf{w}\Vert_2$, where $\Vert \mathbf{w} \Vert_2 = (\sum_{i_1,\hdots,i_D} \vert\mathbf{w}_{i_1,\hdots,i_D}\vert^2)^{1/2}$. For order-two tensors, $\Pi_\epsilon(\mathbf{w})$ is obtained by truncating the SVD of $\mathbf{w}$, which can be computed with standard and efficient algorithms.  For higher-order tensors,  low-rank truncations (in tree-based Tucker format) can be obtained using efficient higher-order SVD algorithms (also allowing a control of the error, see Section \ref{sec:hosvd}) or other optimization algorithms in subsets of low-rank tensors. 
  
 \paragraph{Low-rank truncation controlled in $\infty$-norm.}
 For a control of the error in the $\infty$-norm, the operator $\Pi_\epsilon$ should provide an approximation 
 $\Pi_\epsilon(\mathbf{w})$ of a tensor $\mathbf{w}$ such that $\Vert \mathbf{w} - \Pi_\epsilon(\mathbf{w})\Vert_\infty \le \epsilon \Vert \mathbf{w}\Vert_\infty$, where $\Vert \mathbf{w} \Vert_\infty = \max_{i_1,\hdots,i_D} \vert \mathbf{w}_{i_1,\hdots,i_D} \vert$. A practical implementation of the 
truncation operator $\Pi_\epsilon$ can be based on the Empirical Interpolation Method \cite{BAR02,MAD09} or higher-order extensions of Adaptive Cross Approximation (ACA) algorithms \cite{OSE10}.
 
 \begin{rem}
For the solution of parameter-dependent equations with interpolation methods, simple iterative algorithms take the form 
$$u^{i+1}(\xi) = M(\xi)(u^{i}(\xi)), \quad \xi \in \Xi_K,$$
where $M(\xi)$ is a parameter-dependent iteration map and $\Xi_K$ is a discrete parameter set. For the example of Richarson iterations, exact iterations are 
 $u^{i+1}(\xi) = u^{i}(\xi) - \alpha(B(\xi)u^i(\xi)-f(\xi)),$ $\xi \in \Xi_K.$
 Low-rank truncations of the iterates $\{u^{i}(\xi)\}_{\xi\in \Xi_K}$ should therefore be controlled in the $2$-norm (resp. $\infty$-norm) if one is interested in a mean-square (resp. uniform) control of the error over the discrete parameter set $\Xi_K$. However, note that  under some assumptions on the regularity of a function, controlling the 
 $2$-norm may be sufficient for controlling the $\infty$-norm. 
  \end{rem}
  \begin{rem}
  Note that one could be interested in controlling the error with respect to other norms, such as the norm $\Vert \mathbf{w} \Vert_{(\infty,2)} = \max_{i_1} (\sum_{i_2} \vert \mathbf{w}_{i_1,i_2}\vert^2)^{1/2}$ for an order-two tensor $\mathbf{w}$. For the solution of parameter-dependent equations with interpolation methods, where $\mathbf{w} \in \Rbb^K \otimes \Rbb^{\dim(V)}$ represents the components on an orthonormal basis of $V$ of samples $\{w(\xi)\}_{\xi\in \Xi_K}\in V^K$ of a function $w$ on a discrete parameter set $\Xi_K$, we have $\Vert \mathbf{w} \Vert_{(\infty,2)}= \max_{\xi \in \Xi_K} \Vert w(\xi) \Vert_V$. This allows a uniform control over $\Xi_K$ of the error measured in the $V$-norm. A practical implementation of the 
truncation operator with a control in norm $\Vert \cdot \Vert_{(\infty,2)}$ can be based on the Generalized Empirical Interpolation Method \cite{MAD09}.
 \end{rem}
 \subsection{Minimization of a residual-based distance to the solution} \label{sec:direc-min-residual}

A distance $\Ec({u},{w})$ from $w$ to the solution ${u}$ of \eqref{eq:operator_equation} can be defined by using a  residual norm,
\begin{align}
\Ec(u,w) = \Vert Aw-F \Vert_{D}, \label{eq:minres_metricD}
\end{align}
where $D : \boldsymbol{W}' \rightarrow \boldsymbol{W} $ is an operator which defines on 
$\boldsymbol{W}'$ an inner product norm $\Vert \cdot\Vert_{D} =\sqrt{ \langle D \cdot,\cdot\rangle}$. The operator $D$ plays the role of a preconditioner. 
It can be chosen such that $\Vert \cdot\Vert_{D} = \Vert \cdot\Vert_{\boldsymbol{W}'}$, but it can also be defined in a different way.
For a linear operator $A$, $w\mapsto \Ec(u,w)$ is a quadratic functional. 
\begin{rem}
Note that the distance $\Ec(u,w)$ between the tensors $u$ and $w$ in $\boldsymbol{V}$ corresponds to a distance $\boldsymbol{\Ec}(\mathbf{u},\mathbf{w})$ between the associated tensors $\mathbf{u}$ and $\mathbf{w}$ in $\Rbb^{N_1}\otimes \hdots \otimes \Rbb^{N_D}$, with
$
\boldsymbol{\Ec}(\mathbf{u},\mathbf{w}) = \Vert \mathbf{A}\mathbf{w}-\mathbf{F} \Vert_{\mathbf{D}},
$ for some operator $\mathbf{D} : \mathbf{X} \to \mathbf{X}$.
\end{rem}

Let $\Sc_r$ denote a subset of tensors in $\boldsymbol{V}$ with bounded rank $r$. 
Let $u_r$ denote the minimizer of $w\mapsto  \Ec(u,w)$ over $\Sc_r$, i.e. 
\begin{align}
\Ec(u,u_r) = \min_{w\in \Sc_r}   \Ec(u,w).\label{eq:minres-ur-euv}
\end{align}
 If $A$ is a linear operator such that $\alpha \Vert v\Vert_{\boldsymbol{V}} \le \Vert A v\Vert_{\boldsymbol{W}'} \le \beta \Vert v \Vert_{\boldsymbol{V}}$ and if the operator $D$ is such that $\alpha_D \Vert \cdot\Vert_{\boldsymbol{W}'} \le \Vert \cdot\Vert_{D}\le \beta_D \Vert \cdot\Vert_{\boldsymbol{W}'}$,   then 
 \begin{align}
 \tilde \alpha  \Vert u-w \Vert_{\boldsymbol{V}} \le \Ec(u,v) \le \tilde  \beta \Vert u-w \Vert_{\boldsymbol{V}},
\label{distance-equivalent-V} \end{align} with 
 $\tilde \beta = \beta_D \beta$ and $\tilde \alpha = \alpha_D \alpha$, and 
the solution $u_r$ of \eqref{eq:minres-ur-euv} is a quasi-best approximation in $\Sc_r$ with respect to the norm $\Vert\cdot\Vert_{\boldsymbol{V}}$, with 
 \begin{align}
 \Vert u-u_r\Vert_{\boldsymbol{V}} \le {\frac{\tilde \beta}{\tilde \alpha}} \inf_{w\in \Sc_r}  \Vert u-w \Vert_{\boldsymbol{V}}.\label{eq:quasi-optim-min-res-lowrank}
\end{align}
In practice, an approximation $u_r$ in a certain low-rank format can be obtained by directly solving the optimization problem \eqref{eq:minres-ur-euv} over a subset $\Sc_r$ of low-rank tensors. Constructive algorithms presented in Section \ref{sec:greedy} (which provide only suboptimal approximations) can also be applied and should be preferred 
 when dealing with complex numerical models.  

\begin{rem}
Equation \eqref{eq:quasi-optim-min-res-lowrank} highlights the utility of working with well chosen norms, such that $ {\tilde \beta}/{\tilde \alpha} \approx 1$ if one is interested in minimizing the error in the norm $\Vert \cdot \Vert_{\boldsymbol{V}}$.
In \cite{Cohen:2012,MZA:9370830}, the authors introduce a norm on $\boldsymbol{W}$ such that the residual norm $\Vert Aw-F\Vert_{\boldsymbol{W'}}$ coincides with the error $\Vert w-u\Vert_{\boldsymbol{V}}$, where $\Vert \cdot\Vert_{\boldsymbol{V}}$ is a  norm of interest. Quasi-best approximations are then computed using an iterative algorithm.
\end{rem} 
\begin{rem}
For linear problems, a necessary condition of optimality for problem \eqref{eq:minres-ur-euv} writes\footnote{This  stationarity condition reveals the importance of the analysis of the manifold structure of subsets of tensors with bounded rank (see e.g. \cite{USC12,falco2015geometric}).}
\begin{align}
\langle A u_r-F  , D A \delta w\rangle = 0 \quad \text{for all }\delta w \in T_{u_r} \Sc_r, \label{eq:stationarity_min_res}
\end{align}
where $T_{u_r} \Sc_r$ is the tangent space to the manifold $\Sc_r$ at $u_r$. Since $\Sc_r$ is not a linear space nor a convex set, the condition \eqref{eq:stationarity_min_res} is not a sufficient condition for $u_r$ to be a solution of \eqref{eq:minres-ur-euv}.
\end{rem}
\begin{rem}
For linear symmetric coercive problems, where $\boldsymbol{V}=\boldsymbol{W}$, $A^{-1}$ defines a norm  $\Vert \cdot\Vert_{A^{-1}}$
on $\boldsymbol{W}'$ such that $\Vert F\Vert_{A^{-1}} = \langle F, A^{-1} F \rangle$. Then, letting $D=A^{-1}$, the distance to the solution can be  chosen as 
\begin{align}
\Ec(u,w) =  \Vert Aw-F\Vert_{A^{-1}} =  \Vert w- u\Vert_{A} , \label{eq:sym-coercive-distance}
\end{align}
where $\Vert w \Vert_A^2 = \langle A w,w \rangle$. In this case, the minimization of $w\mapsto \Ec(u,w)$ on a subset $\Sc_r$ provides a best approximation of $u$ in $\Sc_r$ with respect to the operator norm $\Vert\cdot\Vert_A$.
Denoting by $J(w) =  \langle A w,w\rangle - 2\langle F,w\rangle$, we have $\Ec(u,w)^2 = J(w) - J(u)$, so that minimizing $\Ec(u,w)$ is equivalent to minimizing the functional $J(w)$, which is a strongly convex quadratic functional.
\end{rem}

\paragraph{Parameter-dependent equations.}   
 We now consider the particular case of the solution of parameter-dependent equations using Galerkin or interpolation methods. 
For Petrov-Galerkin methods (see Section \ref{sec:galerkin}), the distance $\Ec(u,w)$ can be chosen as in \eqref{eq:minres_metricD} with an operator $D$ such that $\Vert Au \Vert_D = \Vert Au \Vert_{\boldsymbol{W}'}$ or simply  $\Vert Au \Vert_D = \Vert \mathbf{A} \mathbf{u} \Vert_2$. For 
minimal residual Galerkin methods (see Section \ref{sec:parametric-galerkin-minimal-residual}), 
the distance $\Ec(u,w)$ can be chosen such that
\begin{align}
\Ec(u,w)^2 =   \int_{\Xi} \Vert B(y) w(y) - f(y) \Vert_{C(y)}^2 \mu(dy).
 \label{minresgalerkin_optim_functional_innumericalsolution}
\end{align}
which corresponds to $\Ec(u,w) = \Vert w-u \Vert_{\tilde A}^2=\Vert \tilde Aw-\tilde F \Vert_{\tilde A^{-1}}^2$. %
  In the case of interpolation (or collocation) methods (see Section \ref{sec:interpolation-equation}), the distance $\Ec(u,w)$ can be chosen such that 
\begin{align}
\Ec(u,w)^2 = \sum_{k=1}^K \omega^k  \Vert B(y^k) w(y^k) - f(y^k) \Vert_{C(y^k)}^2,\label{minresinterpolation}
\end{align}
with suitable weights $\{\omega^k\}_{k=1}^K$
\footnote{By identifying an element $F\in (\Rbb^K\otimes W)'$ with an element $\{f_k\}_{k=1}^K \in (W')^K$, \eqref{minresinterpolation} corresponds to \eqref{eq:minres_metricD} with an operator $D$ such that $\Vert F \Vert_D^2 =  \sum_{k=1}^K \omega^k \Vert f_k \Vert_{C(y^k)}^2$, i.e. $D = \sum_{k=1}^K \Omega_k \otimes C(y^k)$, with $\Omega_k \in \Rbb^{K\times K}$ such that 
$(\Omega^k)_{ij} = \delta_{ik} \delta_{jk} \omega^k$. For $C(y)=C$ independent of $y$, $D = \Omega \otimes C$, with $\Omega = \mathrm{diag}(\omega^1,\hdots,\omega^k) \in \Rbb^{K\times K}$.}.
In both cases \eqref{minresgalerkin_optim_functional_innumericalsolution} and \eqref{minresinterpolation}, with a linear operator $B(\xi)$ satisfying the 
 assumptions of Section \ref{sec:parametric-galerkin-minimal-residual},   
Property \eqref{distance-equivalent-V} is satisfied, where in the case of interpolation methods, 
  $\Vert \cdot\Vert_{\boldsymbol{V}}$ coincides with the norm $\Vert \cdot\Vert_{2,K}$ defined by \eqref{p-K-norm}.

\begin{rem} 
The distance could also be chosen as
$$
\Ec(u,w) = \sup_{1\le k \le K}  \Vert B(y^k)w(y^k) - f(y^k) \Vert_{C(y^k)},$$ 
therefore moving from a Hilbert setting to a Banach setting. This is the classical framework for the so-called Reduced Basis methods. With a linear operator $B(\xi)$ satisfying the assumptions of  
 Section \ref{sec:parametric-galerkin-minimal-residual},  Property \eqref{distance-equivalent-V} is satisfied,
where  $\Vert \cdot\Vert_{\boldsymbol{V}}$ coincides with the norm $\Vert \cdot\Vert_{\infty,K}$ defined by \eqref{infinite-K-norm}. The optimal rank-$r$ approximation $u_r$ in $\Rbb^K \otimes V$ such that $\Ec(u,u_r) = \min_{w\in \Rc_r} \Ec(u,w)$ satisfies $$\Vert u-u_r \Vert_{\infty,K} \le \frac{\tilde \beta}{\tilde \alpha} \min_{v\in \Rc_r} \Vert u-w\Vert_{\infty,K} =  \frac{\tilde \beta}{\tilde \alpha} d_r(u(\Xi_K))_V,$$
where $d_r(u(\Xi_K))_V$ is the Kolmogorov $r$-width of the discrete set of solutions $u(\Xi_K)$ in $V$. In practice, one can rely on an algorithm based on a greedy construction of subspaces $V_r \subset  V$, such as presented  in section 
\ref{sec:low-rank-sample-norm} (replacing $\Vert w(y) -u(y) \Vert_V$ by $\Vert B(y) w(y) - f(y) \Vert_{C(y)}$). This is the so-called  offline phase of Reduced Basis methods\citechap{ (see \refchapRB{})}, and convergence results 
for this algorithm can be found in \cite{Binev:2011,buffa:2012,DeVore:2013fk}, where   
the error $\Vert u-u_r \Vert_{\infty,K} = \sup_{y\in \Xi_K} \Vert u(y) - P_{V_r} u(y)\Vert_V$ is compared with the best rank-$r$ approximation error 
$\rho^{(\infty,K)}_r(u) = d_r(u(\Xi_K))_V$\citechap{ (see \refchapRBtheory{})}. 
\end{rem}

\subsection{Coupling iterative methods and residual norm minimizations}

Methods presented in Sections \ref{sec:iterative_methods} and \ref{sec:direc-min-residual} can be combined.
It allows the use of a larger class of iterative solvers for which one iteration takes the form $u^{i+1} = M(u^i)$, with 
$M(u^i) = C_i^{-1} G(u^i)$, where $C_i$ is an operator given in low-rank format whose inverse $C_i^{-1}$ is not computable explicitly. At iteration $i$, a low-rank approximation 
 $u^{i+1}$ of $C_i^{-1} G(u^i)$ can be obtained by minimizing the functional $w\mapsto \Ec(C_i^{-1} G(u^i),w) = \Vert C_i w - G(u^i) \Vert_{\star} $, where $\Vert \cdot\Vert_\star$ is some computable residual norm, either by a direct optimization in subsets of low-rank tensors or by using greedy algorithms. 
 
 The above iterations can be associated with an advanced iterative method for solving the linear system $Au=F$  ($C_i$ could be the inverse of a known preconditioner of the operator $A$, or a piece of the operator $A$ in a method based on operator splitting), or with a nonlinear iterative solver for solving a nonlinear equation $A(u) = F$, with $A$ being a nonlinear map. For example, for a Newton solver, $C_i$ would be the differential of $A$ (tangent operator) at $u^i$.

\subsection{Galerkin approaches for low-rank approximations}\label{sec:galerkin-low-rank}

The minimal residual-based approaches presented in Section \ref{sec:direc-min-residual} are robust approaches that guarantee the convergence of low-rank approximations. However, when an approximation $u_r$ is defined as the minimizer of the residual norm $\Vert Au_r - F\Vert_D= \langle u_r - u , A^* D A (u_r -u) \rangle^{1/2}$ with a certain operator $D$, these approaches 
 may suffer 
 from ill-conditioning and they may induce high computational costs since they require  operations between objects with a possibly high rank (operator $A^*DA$ and right-hand side $A^*D Au = A^*D F$). Moreover, in the context of parameter-dependent (or stochastic) equations, they require the solution of problems that have not in general the structure of standard parameter-independent problems, and therefore, they cannot rely on standard parameter-independent (or deterministic) solvers (see Section 
 \ref{sec:Proper Generalized Decompositions for parametric equations in tensor format}). In order to address these issues, low-rank approximations can also be defined using other Galerkin orthogonality criteria (see \cite{NOU10b,cances2013greedy}).

\subsubsection{Galerkin orthogonality}
Let us assume that $\boldsymbol{V}=\boldsymbol{W}$ .
Let $\Sc_r$ denote a subset of tensors in $\boldsymbol{V}$  with a rank bounded by $r$.
An approximation $u_r$ in $\Sc_r$  can be searched such that 
the residual $Au_r-F$ is orthogonal to the tangent space $T_{u_r} \Sc_r$  to the manifold $\Sc_r$ at $u_r$, i.e. such that
\begin{align}
\langle A u_r - F , \delta w\rangle = 0  \quad \text{for all } \delta w\in T_{u_r} \Sc_r.\label{eq:ortho-galerkin}
\end{align}
\begin{rem}
If $\Sc_r$ admits a multilinear parametrization of the form \eqref{multilinear-param-Sr}, then $u_r \in \Sc_r$ can be written $u_r = F_{\Sc_r}(p_1,\hdots,p_M)$. Assuming that the parameters are in vector spaces $P_i$,  \eqref{eq:ortho-galerkin} is equivalent to a set of $M$ coupled equations on the parameters $(p_1,\hdots,p_M)\in P_1\times \hdots \times P_M$:
\begin{align}
\langle A F_{\Sc_r}(p_1,\hdots,p_M) - F , F_{\Sc_r}(p_1,\hdots,\delta p_i,\hdots,p_M)\rangle = 0 \quad \forall \delta p_i\in P_i, \label{eq:ortho-galerkin-parameters}
\end{align}
$1\le i\le M$.
\end{rem}

In order to simplify the presentation of this formulation, let us consider the particular case where an approximation $u_{r-1}$ of rank $r-1$ in $S\otimes V$ is given and let us define $\Sc_r = u_{r-1} + \Rc_1$. Then $u_r \in \Rc_r$ is searched under the form  $u_r= u_{r-1} + w_r$, where $w_r = s_r\otimes v_r \in \Rc_1$ is a rank-one correction of $u_{r-1}$ which must satisfy the Galerkin orthogonality condition \eqref{eq:ortho-galerkin}. Since $T_{u_r} \Sc_r = T_{w_r} \Rc_1 = \{\delta w = s \otimes v_r + s_r\otimes v : s\in S, v\in V \}$, the 
 condition \eqref{eq:ortho-galerkin} becomes 
\begin{align}
\langle A w_r  - (F-Au_{r-1}) , \delta w\rangle = 0 \quad \text{for all } \delta w\in T_{w_r} \Rc_1,\label{eq:ortho-galerkin-R1}
\end{align}
or equivalently 
\begin{subequations}
\begin{align}
&\langle A s_r\otimes v_r  - (F-Au_{r-1}) , s \otimes v_r \rangle =  0 \quad \text{for all } s \in  S,\label{eq:ortho-galerkin-R1-s}
\\
&\langle A s_r\otimes v_r  - (F-Au_{r-1}) ,  s_r \otimes v\rangle =  0 \quad \text{for all } v \in  V.
\label{eq:ortho-galerkin-R1-v}
\end{align}
\end{subequations}
Equation \eqref{eq:ortho-galerkin-R1} may not have any solution $w_r$ or may have many solutions (possibly infinitely many for problems formulated in infinite-dimensional spaces), with particular solutions that are not relevant for the approximation of $u$. 
In practice, heuristic algorithms are used to solve equation \eqref{eq:ortho-galerkin-R1}, such as alternating direction algorithms. This consists in solving successively equation \eqref{eq:ortho-galerkin-R1-s} with a fixed $v_r$ and 
equation  \eqref{eq:ortho-galerkin-R1-v} with a fixed $s_r$.
It has to be noted that  when  \eqref{eq:ortho-galerkin-R1} admits solutions, this heuristic algorithm selects particular solutions $u_r$ that are usually relevant. It can be understood in the case where $A$ is symmetric and defines a norm, since the alternating direction algorithm coincides with an alternating minimization algorithm. This explains why solutions that minimize a certain residual norm are selected.
\begin{rem}
As an illustration, let us consider the case where 
$A$ is such that $\langle A v,w\rangle = \langle v,w\rangle$, with $\langle \cdot,\cdot\rangle$ the canonical norm on $S\otimes V$. For $r=1$, equation \eqref{eq:ortho-galerkin} writes $\langle  s_1\otimes v_1 -u , s \otimes v_1 + s_1\otimes v  \rangle = 0
$ for all $(s,v) \in S\times V$. This implies that $u(v_1) = \Vert v_1\Vert_V^2 s_1$ and $u^*(s_1) = \Vert s_1\Vert_S^2 v_1$, that means that the solutions of  equation \eqref{eq:ortho-galerkin} are the tensors $s_1 \otimes v_1$ where $v_1$ and $s_1$ are right and left singular vectors of $u$ associated with a certain non-zero singular value $\Vert v_1\Vert_V^2 = \Vert s_1\Vert^2_S $. In this case, the alternating direction algorithm corresponds to a power method for finding the dominant eigenvector of $u^*\circ u$. That means that the algorithm allows selecting the optimal rank-one approximation with respect to the canonical norm $\Vert\cdot\Vert$ among possible solutions of equation \eqref{eq:ortho-galerkin-R1}. 
See  \cite{NOU08b,NOU10b} for further details on algorithms and their interpretation as algorithms for solving invariant subspace problems associated with generalization of singular value decompositions.
\end{rem}
In many applications, Galerkin orthogonality criteria (in conjunction with suitable algorithms for solving \eqref{eq:ortho-galerkin}) provide rather good low-rank approximations, although they are not  associated with the minimization of a certain distance to the solution, and therefore do not guarantee the convergence with the rank. Note that these Galerkin approaches have also been applied successfully to some nonlinear problems ($A$ being a nonlinear map), see \cite{NOU09b,Tamellini:2014}.
\paragraph{Parameter-dependent equations.}
For the case of parameter-dependent equations, the equation \eqref{eq:ortho-galerkin-R1-v} for $v_r \in V$ (with fixed $s_r$) is a parameter-independent equation of the form $$\widehat B_{r,r} v_r= \widehat f_r - \sum_{i=1}^{r-1} \widehat B_{r,i} v_i,$$ where $\widehat B_{r,i} =  \int_\Xi B(y)s_r(y)s_i(y) \mu(dy)$ and $\widehat f_r = \int_\Xi f(y)s_r(y)\mu(dy)$ for Galerkin methods, or $\widehat B_{r,i} =  \sum_{k=1}^K \omega^k B(y^k)s_r(y^k)s_i(y^k) $ and $\widehat f_r =  \sum_{k=1}^K \omega^k f(y^k)s_r(y^k)$ for interpolation (or collocation) methods. When $B(\xi)$ and $f(\xi)$ admit affine representations of the form \eqref{eq:decomp_B} and \eqref{eq:decomp_f} respectively, then $\widehat B_{r,i}$ and $\widehat{f}_r$ can be interpreted as evaluations of $B$ and $f$ for particular values of parameter-dependent functions $\lambda_i$ and $\gamma_i$. Therefore, equation \eqref{eq:ortho-galerkin-R1-v} can be usually solved with a standard solver for parameter-independent or deterministic  models   (see Section \ref{sec:Proper Generalized Decompositions for parametric equations in tensor format} for further discussion and illustration on model example 1).

\subsubsection{Petrov-Galerkin orthogonality}
For non-symmetric problems,  possibly with $\boldsymbol V \neq \boldsymbol W$, it has been proposed in \cite{NOU10b} an alternative construction based on a Petrov-Galerkin orthogonality criteria. At step $r$, assuming that $u_{r-1} \in \Rc_{r-1} \subset \boldsymbol V$ is known,  a rank-one correction $w_r=s_r\otimes v_r \in \boldsymbol V$ and an auxiliary rank-one element $\tilde w_r = \tilde s_r \otimes \tilde v_r \in \boldsymbol{W} $ are constructed such that 
\begin{align*}
&\langle A w_r  - (F-Au_{r-1}) , \delta \tilde w\rangle = 0, \quad \text{for all } \delta \tilde w\in T_{\tilde w_r} \Rc_1 \subset \boldsymbol{W}, 
\\
&\langle A \delta w , \tilde w_r\rangle = \langle \delta w, w_r\rangle_{\boldsymbol{V}}, \quad \text{for all } \delta w\in T_{w_r} \Rc_1 \subset \boldsymbol{V}, 
\end{align*}
or equivalently
\begin{subequations}
\begin{align}
&\langle A s_r\otimes v_r  - (F-Au_{r-1}) , \tilde s \otimes \tilde v_r \rangle =  0, \quad \text{for all } \tilde s \in  \tilde S,\label{eq:ortho-petrov-galerkin-R1-s}
\\
&\langle A s\otimes v_r   , \tilde s_r \otimes \tilde v_r \rangle =  \langle s\otimes v_r  , s_r\otimes v_r\rangle_{\boldsymbol{V}}, \quad \text{for all }  s \in   S,\label{eq:ortho-petrov-galerkin-R1-tildes}
\\
&\langle A s_r\otimes v_r  - (F-Au_{r-1}) , \tilde s_r \otimes \tilde v\rangle =  0, \quad\text{for all } \tilde v \in  W,
\label{eq:ortho-petrov-galerkin-R1-v}\\
&\langle A s_r\otimes v  , \tilde s_r \otimes \tilde v_r\rangle =   \langle s_r\otimes v , s_r\otimes v_r\rangle_{\boldsymbol{V}}, \quad \text{for all }  v \in  V.
\label{eq:ortho-petrov-galerkin-R1-tildev}
\end{align}
\end{subequations}

A heuristic alternating direction algorithm can be used, which consists in solving successively equations \eqref{eq:ortho-petrov-galerkin-R1-s} to \eqref{eq:ortho-petrov-galerkin-R1-tildev}  respectively for $s_r$, $\tilde s_r$, $v_r$, $\tilde v_r$ (see \cite{NOU10b} for the application to evolution problems and \cite{cances2013greedy} for the application to parameter-dependent equations). 
\paragraph{Parameter-dependent equations.}
For the case of parameter-dependent equations, we note that problems \eqref{eq:ortho-petrov-galerkin-R1-v} and \eqref{eq:ortho-petrov-galerkin-R1-tildev} are parameter-independent equations respectively of the form $$\widehat B_{r,r}v_r= \hat f_r - \sum_{i=1}^{r-1} \widehat B_{r,i} v_i, \quad \text{and} \quad \widehat B_{r,r}^*\tilde v_r =\hat s_r  R_V v_r, $$
where $R_V : V\rightarrow V'$ is such that $\langle R_V v , \hat v\rangle = \langle v, \hat v\rangle_V$, and 
where $\widehat B_{r,i} = \int_\Xi B(y)s_i(y) \tilde s_r(y) \mu(dy)$, $\widehat f_r = \int_\Xi f(y)\tilde s_r(y)\mu(dy)$ and $  \hat s_r = \int_{\Xi} s_r(y)^2\mu(dy)$ for the case of Galerkin methods, or 
$\widehat B_{r,i} = \sum_{k=1}^K \omega^k B(y^k)s_i(y^k) \tilde s_r(y^k)$, $\widehat f_r = \sum_{k=1}^K \omega^k  f(y^k)\tilde s_r(y^k)$ and $\hat s_r = \sum_{k=1}^K \omega^k s_r(y^k)^2$ for the case of interpolation (or collocation) methods. When $B(\xi)$ and $f(\xi)$ admit affine representations of the form \eqref{eq:decomp_B} and \eqref{eq:decomp_f} respectively, then $\widehat B_{r,i}$ and $\widehat{f}_r$ can be interpreted as evaluations of $B$ and $f$ for particular values of parameter-dependent functions $\lambda_i$ and $\gamma_i$. Therefore, equations \eqref{eq:ortho-petrov-galerkin-R1-v} and \eqref{eq:ortho-petrov-galerkin-R1-tildev} can be usually solved with a standard solver for parameter-independent or deterministic models\footnote{For example, when applied to model example 2, equation \eqref{eq:ortho-petrov-galerkin-R1-v} is a weak form of a deterministic evolution equation 
$\hat \alpha \frac{\partial v_r}{\partial t} - \nabla \cdot (\hat \kappa  \nabla v_r) = \hat g$, with initial condition
$v_r(\cdot,0) = \Ebb_\mu(u_0(\cdot,\xi)\tilde s_r)$, and where $\hat \alpha =\Ebb_\mu(s_i(\xi)\tilde s_r(\xi))$ and $\hat \kappa(\cdot) = \Ebb_\mu(\kappa(\cdot,\xi) s_r(\xi)\tilde s_r(\xi))$.}.
 
\subsection{Greedy construction of subspaces for parameter-dependent equations}\label{sec:Proper Generalized Decompositions for parametric equations in tensor format}
Any of the algorithms presented in  Section \ref{sec:greedy} can be applied to construct a low-rank approximation $w$ of the solution $u$ when a suitable measure $\Ec(u,w)$ of the error has been defined. However, the variants based on the progressive construction of reduced spaces  
 (see Section \ref{sec:partial_greedy_subspace}) are particularly pertinent in the context of parameter-dependent problems since they  only involve the solution of a sequence of problems with the complexity of a parameter-independent problem, and of reduced order parameter-dependent models which are the projections of the initial model on the reduced spaces $V_r$.
Moreover, specific algorithms can take advantage of the particular structure of the parameter-dependent model, so that parameter-independent equations have the structure of standard problems which can be solved with available solution codes for parameter-independent models. The reader is referred to \cite{NOU07,NOU08b,NOU10,CHE12} for practical implementations and illustrations of the behavior of these algorithms.
\\
 
 Here, we illustrate the application of the algorithm presented in Section  \ref{sec:partial_greedy_subspace} for the computation of a sequence of low-rank approximations in $S\otimes V$. The algorithm relies on the construction of optimal nested subspaces $V_r \subset V$. 
For simplicity, we only consider the case of a symmetric coercive problem (e.g. model example 1 described in Section \ref{sec:example1}), with a residual-based error $\Ec(u,w) $ defined by \eqref{eq:sym-coercive-distance}, which corresponds to $\Ec(u,w)^2 = J(w) - J(u) $ with $J(w) = \langle A w,w \rangle - 2\langle F,w\rangle$. We have 
$$
J(w) = \int_\Xi \big( \langle B(y) w(y) , w(y) \rangle - 2\langle f(y),w(y) \rangle\big)\mu(dy)
$$
for the case of Galerkin methods, or 
$$
J(w) = \sum_{k=1}^K \omega^k \big( \langle B(y^k) w(y^k) , w(y^k) \rangle - 2\langle f(y^k),w(y^k) \rangle\big)
$$
for the case of interpolation (or collocation) methods. In order to simplify the presentation, $\Ebb_\mu(g(\xi))$ will  denote either $\int_\Xi g(y) \mu(dy) $ in the case of Galerkin methods, or $\sum_{k=1}^K \omega^k g(y^k)$ in the case of interpolation (or collocation) methods. With this notation, we have 
$$J(w) = \Ebb_\mu(\langle B(\xi) w(\xi) , w(\xi) \rangle - 2\langle f(\xi),w(\xi) \rangle).$$

\begin{rem}\label{rem:greedy-subspace-minres}
For non-symmetric problems, such as model example 2 described in Section \ref{sec:example2}, one can adopt the formulation presented in Section \ref{sec:parametric-galerkin-minimal-residual} and use the expression  \eqref{minresgalerkin_optim_functional_innumericalsolution} (or \eqref{minresinterpolation})
for $\Ec(u,w)$. The application of the algorithm follows the same lines, where we simply 
replace  operators and right-hand sides ($A$, $F$,  $B(\xi)$, $f(\xi)$) by their tilded versions ($\tilde A$, $\tilde F$, $\tilde B(\xi) = B(\xi)^*C(\xi) B(\xi)$, $\tilde f(\xi) = B(\xi)^*C(\xi) f(\xi)$), and where the approximation is searched as the minimizer of the functional 
$
J(w) = \langle \tilde A w,w \rangle - 2\langle \tilde F,w\rangle. 
$
\end{rem}

 The algorithm is defined by \eqref{eq:optimal-partial-greedy}. At iteration $r$, the $(r-1)$-dimensional reduced basis $\{v_1,\hdots,v_{r-1}\} $ of the subspace $V_{r-1}\subset V$ being given, the rank-$r$ approximation $u_r = \sum_{i=1}^{r} s_i^{(r)} \otimes v_i$ is defined by
$$
J(u_r) = \min_{\substack{V_r\in \Gbb_r(V) \\ V_r \supset V_{r-1}}} \min_{w\in S\otimes V_r}  J(w) = \min_{v_r\in V_r}  \min_{\{s_i\}_{i=1}^r \in S^r} J( \sum_{i=1}^{r} s_i \otimes v_i).
$$
For solving this optimization problem, we can use an alternating minimization algorithm, solving alternatively 
\begin{subequations}\begin{align}
\min_{v_r\in V}  J( \sum_{i=1}^{r} s_i  \otimes v_i),\label{eq:alternating-vr}
\\
\min_{\{s_i\}_{i=1}^r \in S^r} J( \sum_{i=1}^{r} s_i  \otimes v_i). \label{eq:alternating-si}
\end{align}
\end{subequations}
\paragraph{Solution of Problem \eqref{eq:alternating-vr} (a parameter-independent equation).}
Problem \eqref{eq:alternating-vr} is equivalent to solving the equation
\begin{align}
\widehat B_{r,r} v_r = \widehat{f}_r - \sum_{i=1}^{r-1} \widehat B_{r,i} v_i,\label{eq:alternating-vr-operator-equation}
\end{align}
where the operators $\widehat{B}_{r,i} : V\rightarrow W'$ and the vector $\widehat{f}_r \in W'$ are defined by
$$
\widehat{B}_{r,i} = \Ebb_\mu(B(\xi) s_r(\xi) s_i(\xi)) \quad \text{and} \quad \widehat{f}_r = \Ebb_\mu(f(\xi)s_r(\xi) ).$$ When $B(\xi)$ and $f(\xi)$ admit affine representations of the form \eqref{eq:decomp_B} and \eqref{eq:decomp_f} respectively, 
then $\widehat B_{r,i}$ and $\widehat{f}_r$ take the form 
\begin{align*}
&\widehat{B}_{r,i} = \sum_{l=1}^R B_l \widehat{\lambda}_{l,r,i} \quad \text{and} \quad \widehat{f}_r = \sum_{l=1}^L f_l \widehat{\gamma}_{l,r},
\end{align*}
where 
\begin{align*}
\widehat{\lambda}_{l,r,i} = \Ebb_\mu( \lambda_l(\xi)s_r(\xi)s_i(\xi)) \quad \text{and} \quad 
 \widehat{\gamma}_{l,r} = \Ebb_\mu( \gamma_l(\xi)s_r(\xi)).
\end{align*}
Let us emphasize that the operator $\widehat{B}_{r,i}= \sum_{l=1}^R B_l \widehat{\lambda}_{l,r,i}$ has the same structure as the parameter-dependent operator $B(\xi)= \sum_{l=1}^R B_l {\lambda}_l(\xi)$, but  $\widehat{\lambda}_{l,r,i}$ does not correspond to an evaluation of the function $\lambda_l(\xi)$ at some particular values of $\xi$. However, looking at $B$ as a family of operators parametrized by the $\lambda_l$, then 
$\widehat{B}_{r,i}$ corresponds to an evaluation of $B$ at some given values $\widehat{\lambda}_{l,r,i}$ of the parameters $\lambda_l$. In practical applications, that means that this problem can be solved with standard solvers (for parameter-independent  or deterministic models).
\begin{ex}
When this algorithm is applied to model example 1 (see Section \ref{sec:example1}), $\widehat f_r$ and $\widehat B_{r,i}$  are such that 
$$
\langle \widehat f_r,w \rangle = \int_D \widehat g_r \, w \quad \text{and} \quad  \langle \widehat B_{r,i} v , w \rangle = \int_{D}  \nabla  w \cdot \widehat \kappa_{r,i} \cdot \nabla v,
$$
with $\widehat g_r(\cdot) =  \Ebb_{\mu}(g(\cdot,\xi) s_r(\xi))\in L^2(D)$ and $\kappa_{r,i}(\cdot) = \Ebb_\mu(\kappa(\cdot,\xi) s_r(\xi)s_i(\xi))$. Problem  \eqref{eq:alternating-vr} therefore corresponds to the solution of  the deterministic diffusion equation 
$
-\nabla \cdot ( \widehat \kappa_{r,r} \nabla v_r) = \hat g_r + \sum_{i=1}^{r-1} \nabla \cdot (\widehat \kappa_{r,i} \nabla v_i)$.
\end{ex}

\paragraph{Solution of Problem \eqref{eq:alternating-si} (a reduced order parameter-dependent equation).}
Problem \eqref{eq:alternating-si} is equivalent to computing an approximation of the solution in $S\otimes V_r$, where $V_r$ is a reduced space with basis 
$\{v_1,\hdots,v_r\}$. Let $\mathbf{B}(\xi) \in \Rbb^{r\times r}$ be the parameter-dependent matrix   defined by $\mathbf{B}(\xi) = (\langle B(\xi) v_j , v_i\rangle)_{i,j=1}^r =  (b( v_j , v_i;\xi))_{i,j=1}^r$, and let $\mathbf{f}(\xi) \in \Rbb^r$ be the parameter-dependent vector   defined by $\mathbf{f}(\xi) = (\langle f(\xi),v_i \rangle)_{i=1}^r$. 
If $B(\xi)$ and $f(\xi)$ admit affine representations of the form \eqref{eq:decomp_B} and \eqref{eq:decomp_f} respectively, then 
$\mathbf{B}(\xi) = \sum_{l=1}^R \mathbf{B}_l \lambda_l(\xi)$ and $\mathbf{f}(\xi) = \sum_{l=1}^L \mathbf{f}_l\gamma_l(\xi)$, where the matrices $\mathbf{B}_l \in \Rbb^{r\times r}$ and vectors $\mathbf{f}_l \in \Rbb^r$ are associated with projections on the reduced spaces $V_r$ of operators $B_l$ and vectors $f_l$ respectively. Then, denoting $\mathbf{s} = (s_i)_{i=1}^r \in S^r = S\otimes \Rbb^r$, Problem \eqref{eq:alternating-si} is equivalent to
\begin{align}
\Ebb_\mu( \mathbf{t}(\xi)^T \mathbf{B}(\xi) \mathbf{s}(\xi)) = \Ebb_\mu(\mathbf{t}(\xi)^T \mathbf{f}(\xi)) \quad \forall \mathbf{t} \in S\otimes \Rbb^r,
\label{eq:reduced-equation-s}
\end{align}
which requires the solution of a system of $\dim(S)\times r$ equations. When $\dim(S)$ is large, order reduction methods can also be used at this step in order to obtain a reduced order approximation of the solution $\mathbf{s}$. For example, for the case of high-dimensional parameter-dependent models with a projection on a tensor-structured approximation space $S$ or with an interpolation on a tensor-structured grid, we can rely on sparse approximation methods or higher-order low-rank methods presented in Section \ref{sec:Low-rank approximation of higher-order tensors}.
Note that in the case of Galerkin methods, 
Equation \eqref{eq:reduced-equation-s} defines the Galerkin approximation 
of the reduced-order parameter-dependent equation 
\begin{align}
\mathbf{B}(\xi) \mathbf{s}(\xi) = \mathbf{f}(\xi),\label{eq:reduced-equation-s-strong}
\end{align}
so that an approximation of $\mathbf{s}$ can also be obtained by sampling-based approaches, based on many sample evaluations $\mathbf{s}(y^k) = \mathbf{B}(y^k)^{-1}\mathbf{f}(y^k) $ (only requiring the solution of reduced systems of equations).

\begin{rem}
As mentioned in Remark \ref{rem:greedy-subspace-minres}, this algorithm can be applied to non symmetric problems such as model example 2 (described in section \ref{sec:example2}) by using minimal residual formulations. However, 
when applied to this evolution problem, the algorithm requires the solution of parameter-independent problems of the form \eqref{eq:alternating-vr-operator-equation} which are global over space-time domain (time stepping methods cannot be used) and may be computationally intractable. An additional order reduction can be introduced by also exploiting the  tensor structure  of the space $V=V(D)\otimes V(I)$ of space-time functions (see Remark \ref{rem:tensor-structure-space-time-model-example2}). Low-rank methods that exploit this structure allow the  complexity of the representations of space-time functions to be reduced.
\end{rem}




\end{document}